\documentclass[a4paper,12pt,amsart,frenchb]{article}

\usepackage{amsmath,amsbsy,amsfonts,amssymb}
\usepackage[french]{babel}
\newcommand{\ass}[2]{\vskip0.3cm\noindent
{\bf {#1}}. { \sl {#2}}\vskip0.3cm\noindent
}

\begin{document}

  \title{ Repr\'esentations   de r\'eduction unipotente pour $SO(2n+1)$, III: exemples de fronts d'onde}
\author{J.-L. Waldspurger}
\date{2 d\'ecembre  2016}
\maketitle

\bigskip

{\bf Introduction}

\bigskip

Cet article est la suite de \cite{W3} et \cite{W4}.  Le corps de base $F$ est local, non-archim\'edien et de caract\'eristique nulle. On note $p$ sa caract\'eristique r\'esiduelle. Un  entier $n\geq1$ est fix\'e pour tout l'article. On suppose $p\geq 6n+4$.  On  introduit les groupes $G_{iso}$ et $G_{an}$  suivants.   Le groupe $G_{iso}$ est le groupe sp\'ecial orthogonal d'un espace $V_{iso}$ de dimension $2n+1$ sur $F$ muni d'une forme quadratique $Q_{iso}$ et $G_{an}$  est le groupe sp\'ecial orthogonal d'un espace $V_{an}$ de dimension $2n+1$ sur $F$ muni d'une forme quadratique $Q_{an}$. Le groupe $G_{iso}$ est 
d\'eploy\'e et $G_{an}$ en est la forme int\'erieure non d\'eploy\'ee. Pour un indice $\sharp=iso$ ou $an$, on note $Irr_{tunip,\sharp}$ l'ensemble des classes d'isomorphismes de repr\'esentations admissibles irr\'eductibles de $G_{\sharp}(F)$ qui sont temp\'er\'ees et de r\'eduction unipotente, cf. \cite{W3} 1.3 pour la d\'efinition de cette propri\'et\'e. On note $Irr_{tunip}$ la r\'eunion disjointe de $Irr_{tunip,iso}$ et $Irr_{tunip,an}$. Pour une partition symplectique $\lambda$ de $2n$, fixons un homomorphisme alg\'ebrique $\rho_{\lambda}:SL(2;{\mathbb C})\to Sp(2n;{\mathbb C})$ param\'etr\'e par $\lambda$, cf. \cite{W3} 1.3. On note $Z(\lambda)$ le commutant dans $Sp(2n;{\mathbb C})$ de l'image de $\rho_{\lambda}$. Soit $s\in Z(\lambda)$ un \'el\'ement semi-simple dont toutes les valeurs propres sont de module $1$. On note $Z(s,\lambda)$ le commutant de $s$ dans $Z(\lambda)$, ${\bf Z}(\lambda,s)$ son groupe de composantes connexes et ${\bf Z}(\lambda,s)^{\vee}$ le groupe des caract\`eres de ${\bf Z}(\lambda,s)$.  La param\'etrisation   de Langlands prend la forme suivante, cf. \cite{W3} 1.3: $Irr_{tunip}$ est param\'etr\'e par l'ensemble des classes de conjugaison (en un sens facile \`a pr\'eciser) de triplets $(\lambda,s,\epsilon)$, o\`u $\lambda$ et $s$ sont comme ci-dessus et $\epsilon\in {\bf Z}(s,\lambda)^{\vee}$. On note $\mathfrak{Irr}_{tunip}$ cet ensemble de triplets. Ce param\'etrage a \'et\'e obtenu par diff\'erents auteurs: Lusztig, cf.  \cite{L0}; Moeglin, cf. \cite{M2} th\'eor\`eme 5.2; Arthur, dans le cas du groupe $G_{iso}$, cf. \cite{Ar} th\'eor\`eme 2.2.1.
   Dans \cite{MW1} et \cite{W4}, on a montr\'e que les repr\'esentations construites par Lusztig  v\'erifiaient les propri\'et\'es de compatibilit\'e \`a l'endoscopie qui les caract\'erisent. En particulier, dans le cas du groupe $G_{iso}$, ces repr\'esentations  sont les m\^emes que celles d'Arthur. Pour $(\lambda,s,\epsilon)\in \mathfrak{Irr}_{tunip}$, on note $\pi(\lambda,s,\epsilon)$ la repr\'esentation temp\'er\'ee qui lui est associ\'ee par Lusztig. L'involution introduite par Zelevinsky dans le cas du groupe $GL(n)$ a \'et\'e g\'en\'eralis\'ee par Aubert et par Schneider et Stuhler aux groupes r\'eductifs quelconques. On la note $D$ et on pose $\delta(\lambda,s,\epsilon)=D(\pi(\lambda,s,\epsilon))$. 

Soit $\sharp=iso$ ou $an$ et soit $\pi$ une repr\'esentation admissible irr\'eductible de $G_{\sharp}(F)$. Notons $\mathfrak{g}_{\sharp}$ l'alg\`ebre de Lie de $G_{\sharp}$. Harish-Chandra a prouv\'e que, dans un voisinage de l'origine, le caract\`ere de $\pi$, descendu  par l'exponentielle \`a $\mathfrak{g}_{\sharp}(F)$,  \'etait combinaison lin\'eaire de transform\'ees de Fourier d'int\'egrales orbitales nilpotentes. Fixons une cl\^oture alg\'ebrique $\bar{F}$ de $F$ et notons $\bar{{\cal N}}(\pi)$  l'ensemble des orbites nilpotentes $\boldsymbol{{\cal O}}$ dans $\mathfrak{g}_{\sharp}(\bar{F})$ v\'erifiant la condition suivante: il existe une orbite nilpotente ${\cal O}$ dans $\mathfrak{g}_{\sharp}(F)$, qui est incluse dans $\boldsymbol{{\cal O}}$ et qui intervient avec un coefficient non nul dans le d\'eveloppement ci-dessus du caract\`ere de $\pi$. On dit que $\pi$ admet un front d'onde si $\bar{{\cal N}}(\pi)$ admet un unique \'el\'ement maximal. Dans ce cas, on dit que cet \'el\'ement maximal est le front d'onde de $\pi$. Les orbites nilpotentes dans $\mathfrak{g}_{\sharp}(\bar{F})$ sont param\'etr\'ees par les partitions orthogonales de $2n+1$ et nous identifierons ces deux ensembles. On conjecture (ce qui est peut-\^etre hasardeux) que toute repr\'esentation admissible irr\'eductible admet un front d'onde. Signalons que, dans le cas o\`u le corps de base est non pas $p$-adique, mais r\'eel, la notion de front d'onde est \'egalement d\'efinie et se r\'ev\`ele importante, cf. par exemple \cite{BV}.

En modifiant quelque peu une construction de  Spaltenstein, on d\'efinit une "dualit\'e" qui envoie une partition symplectique $\lambda$ de $2n$ sur une partition orthogonale $d(\lambda)$ de $2n+1$. La partition $d(\lambda)$ est toujours sp\'eciale et la dualit\'e $d$  n'est pas  bijective (par contre, sa restriction au sous-ensemble des partitions symplectiques sp\'eciales de $2n$ est une bijection entre cet ensemble et celui des partitions orthogonales sp\'eciales de $2n+1$). On d\'emontre dans cet article le r\'esultat suivant.

\ass{Th\'eor\`eme}{Soit $(\lambda,s,\epsilon)\in \mathfrak{Irr}_{tunip}$. Alors la repr\'esentation $\delta(\lambda,s,\epsilon)$ admet un front d'onde et celui-ci est la partition $d(\lambda)$.}

Remarquons que l'on retrouve dans notre cas particulier le  th\'eor\`eme 1.4 de \cite{M}: ce front d'onde est une partition sp\'eciale. Notre th\'eor\`eme n'est pas tr\`es nouveau.  Moeglin a d\'emontr\'e  un r\'esultat similaire en \cite{M2} th\'eor\`eme 3.3.5. Ses hypoth\`eses \'etaient plus g\'en\'erales que les n\^otres. D'une part,  elle consid\'erait tous les groupes classiques et pas seulement les groupes sp\'eciaux orthogonaux. Surtout, elle consid\'erait les repr\'esentations dont le param\`etre de Langlands, sous sa forme habituelle,  se restreint au groupe de Weil  en une somme de caract\`eres d'ordre au plus $2$, \'eventuellement ramifi\'es. Nous nous limitons au cas de r\'eduction unipotente,  ce qui exclut les caract\`eres ramifi\'es. Toutefois, notre r\'esultat n'est pas inclus dans celui de \cite{M2}: avec nos notations, celui-ci suppose que les termes de $\lambda$ sont  tous distincts. La d\'emonstration est aussi enti\`erement diff\'erente.

 Soit $(\lambda,s,\epsilon)\in \mathfrak{Irr}_{tunip}$, posons $\delta=\delta(\lambda,s,\epsilon)$. Notons $\sharp$ l'indice $iso$ ou $an$ tel que $\delta$ soit une repr\'esentation de $G_{\sharp}(F)$. Dans \cite{W1}, on a donn\'e une formule qui calcule la restriction du caract\`ere de $\delta$ aux \'el\'ements compacts de $G_{\sharp}(F)$ (ceux qui sont contenus dans un sous-groupe compact). A fortiori cette formule calcule la restriction du caract\`ere \`a un voisinage de l'origine. Cette restriction est somme de distributions que l'on peut calculer si l'on conna\^{\i}t   les restrictions de $\delta$ aux diff\'erents sous-groupes  compacts maximaux de $G_{\sharp}(F)$, ou plus exactement les repr\'esentations des groupes "r\'esiduels" qui s'en d\'eduisent. La construction de Lusztig donne les renseignements voulus. A partir de l\`a, en utilisant de nombreux travaux de Lusztig (faisceaux-caract\`eres, correspondance de Springer g\'en\'eralis\'ee etc...), on traduit l'assertion \`a d\'emontrer en termes de repr\'esentations de groupes de Weyl. Il s'agit en gros de savoir quelles sont les repr\'esentations qui peuvent intervenir dans certaines restrictions d'une repr\'esentation d'un produit de groupes de Weyl d\'etermin\'ee par $\delta$. C'est un probl\`eme  combinatoire que nous avons longuement \'etudi\'e dans \cite{W2} et les r\'esultats de cette r\'ef\'erence permettent de conclure.  
 
 {\bf Remarque.} Dans \cite{W2}, le groupe \'etait suppos\'e non ramifi\'e, ce qui est le cas de $G_{iso}$ mais pas de $G_{an}$. En fait, cette hypoth\`ese ne servait qu'\`a utiliser des r\'esultats d'homog\'en\'eit\'e qui n'\'etaient alors connus que sous cette hypoth\`ese restrictive. Ils sont maintenant connus sans cette hypoth\`ese, cf. \cite{D}, et la plupart des r\'esultats de \cite{W2}, en particulier ceux que l'on utilisera, s'\'etendent au cas g\'en\'eral. 
 
 \bigskip

Evidemment, il serait tentant d'appliquer la m\^eme m\'ethode non pas \`a la b\^ete repr\'esentation $\delta(\lambda,s,\epsilon)$, mais \`a la repr\'esentation temp\'er\'ee $\pi(\lambda,s,\epsilon)$. Indiquons o\`u est le probl\`eme. Les repr\'esentations des groupes "r\'esiduels" associ\'es \`a $\delta(\lambda,s,\epsilon)$ sont bien calcul\'ees par Lusztig, mais en termes de repr\'esentations  de groupes de Weyl  peu explicites.  Plus pr\'ecis\'ement, il appara\^{\i}t des repr\'esentations non irr\'eductibles   dont la d\'ecomposition en composantes irr\'eductibles est dict\'ee par des variantes de polyn\^omes de Kazhdan-Lusztig. Ces repr\'esentations sont not\'ees $\boldsymbol{\rho}_{\lambda,\epsilon}$ dans notre article, mais il ne s'agit plus du m\^eme couple $\lambda,\epsilon$, notons-les ici $\boldsymbol{\rho}_{\nu,\tau}$. 
Il y a un ordre (partiel) naturel  sur l'ensemble des repr\'esentations irr\'eductibles  et
  on contr\^ole   tr\`es bien le terme minimal  du d\'eveloppement de $\boldsymbol{\rho}_{\nu,\tau}$ en composantes irr\'eductibles. Il s'av\`ere que cela nous suffit pour conclure. Si l'on remplace $\delta(\lambda,s,\epsilon)$ par $\pi(\lambda,s,\epsilon)$, les repr\'esentations $\boldsymbol{\rho}_{\nu,\tau}$  sont remplac\'ees par leur produit tensoriel avec le caract\`ere signe  $sgn$ du groupe de Weyl sous-jacent. Comme on  peut s'y attendre, cela inverse l'ordre: on conna\^{\i}t le terme maximal du d\'eveloppement de $sgn\otimes \boldsymbol{\rho}_{\nu,\tau}$. Mais maintenant, l'ordre va dans le mauvais sens et conna\^{\i}tre le terme maximal ne permet plus de conclure.

  \bigskip

\section{Combinatoire}

\bigskip

\subsection{Partitions et  repr\'esentations des groupes de Weyl}

On appelle partition une classe d'\'equivalence de suites d\'ecroissantes  finies  de nombres entiers positifs ou nuls, deux suites \'etant \'equivalentes si elles ne diff\`erent que par des termes nuls. Pour une telle partition $\lambda=(\lambda_{1}\geq \lambda_{2}\geq...\geq \lambda_{r})$, on pose $S(\lambda)=\sum_{j=1,...,r}\lambda_{j}$ et on note $l(\lambda)$ le plus grand entier $j$ tel que $\lambda_{j}\not=0$. Cas particulier: on note $\emptyset$ la partition $(0,...,0)$ et on pose $l(\emptyset)=0$. On note $mult_{\lambda}$ la fonction sur ${\mathbb N}-\{0\}$ telle que, pour tout $i$ dans cet ensemble, $mult_{\lambda}(i)$ est le nombre d'entiers $j$ tels que $\lambda_{j}=i$. On pose aussi $mult_{\lambda}(\geq i)=\sum_{i'\geq i}mult_{\lambda}(i')$. On note $Jord(\lambda)$ l'ensemble des $i\geq1$ tels que $mult_{\lambda}(i)\geq1$. 
Soit $k\in {\mathbb N}$. A \'equivalence pr\`es, on peut supposer $r\geq k$ et on pose $ S_{k}(\lambda)=\sum_{j=1,...,k}\lambda_{j}$.  
Pour $N\in {\mathbb N}$, on note ${\cal P}(N)$ l'ensemble des partitions $\lambda$ telles que $S(\lambda)=N$. Plus g\'en\'eralement, pour un entier $k\geq1$, on note ${\cal P}_{k}(N)$ l'ensemble des $k$-uples de partitions $(\lambda_{1},...,\lambda_{k})$ tels que $S(\lambda_{1})+...+S(\lambda_{k})=N$. On utilisera plus loin des variantes de cette notation, par exemple ${\cal P}^{symp}_{k}(2N)$ etc... On d\'efinit de la fa\c{c}on usuelle la transposition $\lambda\mapsto {^t\lambda}$ dans ${\cal P}(N)$ et les applications $(\lambda_{1},\lambda_{2})\mapsto \lambda_{1}+\lambda_{2}$ et $(\lambda_{1},\lambda_{2})\mapsto \lambda_{1}\cup \lambda_{2}$ qui envoient ${\cal P}_{2}(N)$ dans ${\cal P}(N)$. On d\'efinit un ordre partiel sur ${\cal P}(N)$: pour deux partitions $\lambda_{1},\lambda_{2}\in {\cal P}(N)$, $\lambda_{1}\leq \lambda_{2}$ si et seulement si $S_{k}(\lambda_{1})\leq S_{k}(\lambda_{2})$ pour tout $k\in {\mathbb N}$. 

Plusieurs notations  ci-dessus se g\'en\'eralisent aux suites finies $\alpha=(\alpha_{1},...,\alpha_{r})$ de nombres r\'eels pas forc\'ement d\'ecroissantes. Par exemple, si $k$ est un entier tel que $0\leq k\leq r$, on pose $S_{k}(\alpha)=\sum_{j=1,...,k}\alpha_{j}$. On utilisera aussi la notation $\alpha_{\leq k}=S_{k}(\alpha)$. Si $\alpha$ et $\beta$ sont deux suites de m\^eme longueur, on note 
  $\alpha+\beta$ la suite $(\alpha_{1}+\beta_{1},...,\alpha_{r}+\beta_{r})$.  

Pour tout ensemble $X$, on note ${\mathbb C}[X]$ l'espace vectoriel complexe de base $X$. Pour tout groupe fini $W$, on note $\hat{W}$ l'ensemble des classes de repr\'esentations irr\'eductibles de $W$. En identifiant une telle repr\'esentation \`a son caract\`ere, l'espace ${\mathbb C}[X]$ s'identifie \`a celui des fonctions de $W$ dans ${\mathbb C}$ qui sont invariantes par conjugaison. 

Soit $N\in {\mathbb N}$. On note $\mathfrak{S}_{N}$ le groupe des permutations de l'ensemble $\{1,...,N\}$. On sait param\'etrer $\hat{\mathfrak{S}}_{N}$  par ${\cal P}(N)$, on note $\rho(\lambda)$ la repr\'esentation irr\'eductible correspondant \`a une partition $\lambda$ (en particulier  la repr\'esentation triviale de $\mathfrak{S}_{N}$ est param\'etr\'ee par la partition $\lambda=(n)$). On note $sgn$ le caract\`ere signe usuel de $\mathfrak{S}_{N}$. Si une repr\'esentation irr\'eductible $\rho$ est param\'etr\'ee par la partition $\lambda$, $\rho\otimes sgn$ est param\'etr\'ee par $^t\lambda$. 

On note $W_{N}$ le groupe de Weyl d'un syst\`eme de racines de type $B_{N}$ ou $C_{N}$ (avec la convention $W_{0}=\{1\}$). On sait param\'etrer $\hat{W}_{N}$ par ${\cal P}_{2}(N)$, on note $\rho(\alpha,\beta)$ la repr\'esentation irr\'eductible correspondant \`a  un couple de partitions $(\alpha,\beta)$ (en particulier, la repr\'esentation triviale est param\'etr\'ee par $((N),\emptyset)$).On note $sgn$ le caract\`ere signe  usuel de $W_{N}$ et $sgn_{CD}$ le caract\`ere dont le noyau est le sous-groupe $W_{N}^D$ d'un syst\`eme de racines de type $D_{N}$. Si une repr\'esentation irr\'eductible $\rho$ est param\'etr\'ee par le couple de partitions $(\alpha,\beta)$, $\rho\otimes sgn$ est param\'etr\'ee par $(^t\beta,{^t\alpha})$ et $\rho\otimes sgn_{CD}$ est param\'etr\'ee par $(\beta,\alpha)$.

Supposons $N\geq1$. Pour $(\alpha,\beta)\in {\cal P}_{2}(N)$, les restrictions \`a $W_{N}^D$ de $\rho(\alpha,\beta)$ et $\rho(\beta,\alpha)$ sont \'equivalentes. Si $\alpha\not=\beta$, ces restrictions sont irr\'eductibles, on les note $\rho^D(\alpha,\beta)$ ou $\rho^D(\beta,\alpha)$. Si $\alpha=\beta$, la restriction de $\rho(\alpha,\alpha)$ \`a $W_{N}^D$ se d\'ecompose en deux repr\'esentations irr\'eductibles, que l'on  note $\rho^D(\alpha,\alpha,+)$ et $\rho^D(\alpha,\alpha,-)$. Elles sont conjugu\'ees par un \'el\'ement de $W_{N}-W_{N}^D$ et on n'aura pas besoin de les distinguer. Toutes les repr\'esentations irr\'eductibles de $W_{N}^D$ sont ainsi obtenues. 

\bigskip

\subsection{Symboles}

Pour tout ensemble fini $X$, on note $\vert X\vert $ le nombre d'\'el\'ements de $X$. Si $X$ est un ensemble de nombres, on note $S(X)$ la somme des \'el\'ements de $X$. Pour tout nombre r\'eel $x$, on note $[x]$ sa partie enti\`ere.

Soit $N\in {\mathbb N}$. Un symbole de rang $N$  est une classe d'\'equivalence de couples $(X,Y)$ de sous-ensembles finis de ${\mathbb N}$, v\'erifiant la condition
 $$S(X)+S(Y)-[(\frac{\vert X\vert +\vert Y\vert -1}{2})^2]=N.$$
La relation d'\'equivalence est engendr\'ee par les deux relations (qui pr\'eservent l'\'egalit\'e pr\'ec\'edente):

 $(X,Y)\sim (X',Y')$ o\`u 
$$X'=\{x+1;x\in X\}\cup \{0\},\,\, Y'=\{y+1;y\in Y\}\cup \{0\};$$

$(X,Y)\sim (Y,X)$. 

{\bf Remarque.} Par abus de terminologie, on  appellera plut\^ot symbole un couple $(X,Y)$ repr\'esentant une classe d'\'equivalence.

\bigskip 

Le d\'efaut d'un symbole $(X,Y)$ est la valeur absolue de $\vert X\vert -\vert Y\vert $ (il ne d\'epend que de la classe d'\'equivalence de $(X,Y)$). Pour $D\in {\mathbb N}$, on note ${\cal S}_{N,D}$ l'ensemble des symboles de rang $N$ et de d\'efaut $D$. 

On regroupe les symboles en familles: deux symboles sont dans la m\^eme famille si et seulement si on peut les repr\'esenter par des couples $(X,Y)$ et $(X',Y')$ tels que $X\cup Y=X'\cup Y'$ et $X\cap Y=X'\cap Y'$. La parit\'e du d\'efaut est constante sur chaque famille. Toute famille de symboles de d\'efaut impair contient un unique symbole sp\'ecial, c'est-\`a-dire repr\'esent\'e par un couple $(X,Y)$ de la forme $X=(x_{1}\geq...\geq x_{r+1})$, $Y=(y_{1}\geq...\geq y_{r})$ et tel que
$$x_{1}\geq y_{1}\geq x_{2}\geq y_{2}\geq... \geq y_{r}\geq x_{r+1}.$$
Toute famille de symboles de d\'efaut pair contient un unique symbole sp\'ecial, c'est-\`a-dire repr\'esent\'e par un couple $(X,Y)$ de la forme $X=(x_{1}\geq...\geq x_{r})$, $Y=(y_{1}\geq...\geq y_{r})$ et tel que
$$x_{1}\geq y_{1}\geq x_{2}\geq y_{2}\geq... \geq y_{r}.$$

Soit $(X,Y)$ un symbole de rang $N$. Fixons un entier $d$ majorant les \'el\'ements de $X\cup Y$. Posons 
$$X'=\{d,...,0\}-\{d-y;y\in Y\},\,\, Y'=\{d,...,0\}-\{d-x;x\in X\}.$$
On v\'erifie que $(X',Y')$ est un symbole de rang $N$. A \'equivalence pr\`es, il ne d\'epend pas du choix de $d$ et ne d\'epend que de la classe d'\'equivalence de $(X,Y)$.   Cette construction d\'efinit une "dualit\'e" $(X,Y)\mapsto d(X,Y)=(X',Y')$ dans l'ensemble des symboles de rang $N$. Cette dualit\'e conserve le d\'efaut et est involutive: $d\circ d$ est l'identit\'e. Elle se restreint en une involution du sous-ensemble des symboles sp\'eciaux. Enfin, deux symboles sont dans une m\^eme famille si et seulement si leurs images par dualit\'e le sont. Autrement dit, si $(X,Y)$ est dans la famille du symbole sp\'ecial $(X^{sp},Y^{sp})$, alors $d(X,Y)$ est dans la famille du symbole sp\'ecial $d(X^{sp},Y^{sp})$. 

  Soit $\rho\in \hat{W}_{N}$. Comme en 1.1, on lui associe un couple de partitions $(\alpha,\beta)\in {\cal P}_{2}(N)$. On choisit un entier $r\geq l(\alpha),l(\beta)$. On pose $X=\alpha+\{r,...,0\}$, $Y=\beta+\{r-1,...,0\}$. Alors $(X,Y)$ est un symbole de rang $N$ et de d\'efaut $1$ dont la classe ne d\'epend pas de $r$.
On pose $symb(\rho)=(X,Y)$. L'application $symb:\hat{W}_{N}\to {\cal S}_{N,1}$ ainsi d\'efinie est bijective. On a $symb(\rho\otimes sgn)=d\circ symb(\rho)$. 

Soit $\rho\in \hat{W}_{N}^D$. Comme en 1.1, on lui associe un couple de partitions $(\alpha,\beta)\in {\cal P}_{2}(N)$. On choisit un entier $r\geq l(\alpha),l(\beta)$. On pose $X=\alpha+\{r-1,...,0\}$, $Y=\beta+\{r-1,...,0\}$. Alors $(X,Y)$ est un symbole de rang $N$ et de d\'efaut $0$ dont la classe ne d\'epend pas de $r$.
On pose $symb(\rho)=(X,Y)$. L'application $symb:\hat{W}_{N}^D\to {\cal S}_{N,0}$ ainsi d\'efinie est surjective. Ses fibres ont un ou deux \'el\'ements, celles \`a deux \'el\'ements \'etant form\'ees des couples de la forme $\rho(\alpha,\alpha,+),\rho(\alpha,\alpha,-)$. On a $symb(\rho\otimes sgn)=d\circ symb(\rho)$. 

\bigskip

\subsection{Correspondance de Springer, cas symplectique}

Soit $n\in {\mathbb N}$. On note ${\cal P}^{symp}(2n)$ l'ensemble des partitions symplectiques de $2n$, c'est-\`a-dire les $\lambda\in {\cal P}(2n)$ telles que $mult_{\lambda}(i)$ est pair pour tout entier $i$ impair. Pour une telle partition, on note $Jord_{bp}(\lambda)$ l'ensemble des entiers $i\geq2$ pairs tels que $mult_{\lambda}(i)\geq1$. Plus pr\'ecis\'ement, pour un entier $k\geq1$, on note $Jord_{bp}^k(\lambda)$ l'ensemble des $i\geq2$ pairs tels que $mult_{\lambda}(i)=k$. 

On note $\boldsymbol{{\cal P}^{symp}}(2n)$ l'ensemble des couples $(\lambda,\epsilon)$ o\`u $\lambda\in {\cal P}^{symp}(2n)$ et $\epsilon\in \{\pm 1\}^{Jord_{bp}(\lambda)}$. La correspondance de Springer g\'en\'eralis\'ee   \'etablit une bijection entre $\boldsymbol{{\cal P}^{symp}}(2n)$ et l'ensemble des couples $(\rho,k)$ tels que

$k\in {\mathbb N}$ et $k(k+1)\leq 2n$;

$\rho\in \hat{W}_{n-k(k+1)/2}$.

On note $(\rho_{\lambda,\epsilon},k_{\lambda,\epsilon})$ le couple associ\'e \`a un \'el\'ement $(\lambda,\epsilon)\in \boldsymbol{{\cal P}^{symp}}(2n)$. L'entier $k_{\lambda,\epsilon}$ se calcule de la fa\c{c}on suivante. Notons $i_{1}>...>i_{m}>0$ les entiers pairs $i$ tels que $mult_{\lambda}(i)$ soit impair. Posons 
$$h= \sum_{j=1,...,m}(-1)^j(1-\epsilon(i_{j})).$$
Alors $k_{\rho,\epsilon}=sup(h,-h-1)$. En particulier, si $\epsilon=1$, c'est-\`a-dire $\epsilon(i)=1$ pour tout $i\in Jord_{bp}(\lambda)$, on a $k_{\lambda,1}=0$ et $\rho_{\lambda,1}\in \hat{W}_{n}$.

Une partition $\lambda=(\lambda_{1}\geq \lambda_{2}\geq...)\in {\cal P}^{symp}(2n)$ est sp\'eciale si et seulement si $\lambda_{2j-1}$ et $\lambda_{2j}$ sont de m\^eme parit\'e pour tout $j\geq1$.  Cela \'equivaut \`a ce que $^t\lambda$ soit symplectique. En notant ${\cal P}^{symp,sp}(2n)$ l'ensemble des partitions symplectiques sp\'eciales de $2n$, l'application $\lambda\mapsto {^t\lambda}$  est une involution de ${\cal P}^{symp,sp}(2n)$. Consid\'erons une  partition $\lambda\in {\cal P}^{symp,sp}(2n)$ et d\'efinissons $i_{1}>...>i_{m}$ comme ci-dessus.  Si $m$ est pair, on pose $m'=m$; si $m$ est impair, on pose $m'=m+1$ et $i_{m'}=0$. On appelle intervalle de $\lambda$ un ensemble $\Delta$ de l'une des formes suivantes:

pour un entier $h=1,...,m'/2$, $\Delta$ est l'ensemble des $i$ tels que $i=0$ ou $i\geq1$ et $mult_{\lambda}(i)\geq1$ et tels que $i_{2h-1}\geq i \geq i_{2h}$;

$\Delta=\{i\}$ o\`u $i$ est un entier pair tel que $i=0$ ou  $i\geq2$ et $mult_{\lambda}(i)\geq1$ et tel qu'il n'existe pas d'entier $h=1,...,m'/2$ de sorte que $i_{2h-1}\geq i\geq i_{2h}$.

Parce que $\lambda$ est sp\'eciale, on v\'erifie que les intervalles sont form\'es d'entiers pairs (c'est \'evident dans le deuxi\`eme cas ci-dessus, un peu moins dans le premier). Ils forment une partition de l'ensemble $Jord_{bp}(\lambda)\cup\{0\}$. On ordonne les intervalles: $\Delta> \Delta'$ si $i> i'$ pour tous $i\in \Delta$, $i'\in \Delta'$. On note $\Delta_{min}$ le plus petit intervalle (c'est celui qui contient $0$). On note $Int(\lambda)$ l'ensemble des intervalles de $\lambda$. 

 L'application $\lambda\mapsto symb(\rho_{\lambda,1})$ est une bijection de ${\cal P}^{symp,sp}(2n)$ sur l'ensemble des symboles sp\'eciaux de rang $n$ et de d\'efaut $1$. Pour $\lambda\in {\cal P}^{symp,sp}(2n)$, on a d\'efini en \cite{W2} VIII.17 une bijection $fam$ entre la famille du symbole $symb(\rho_{\lambda,1})$ et l'ensemble des $(\tau,\delta)\in({\mathbb Z}/2{\mathbb Z})^{Int(\lambda)}\times ({\mathbb Z}/2{\mathbb Z})^{Int(\lambda)}$ tels que $\tau(\Delta_{min})=\delta(\Delta_{min})=0$.

Soit $\lambda\in {\cal P}^{symp}(2n)$. Il existe une unique partition sp\'eciale $sp(\lambda)\in {\cal P}^{symp,sp}(2n)$ telle que $symb(\rho_{\lambda,1})$ et $symb(\rho_{sp(\lambda),1})$ soient dans la m\^eme famille. Il est connu que $\lambda\leq sp(\lambda)$ et que $sp(\lambda)$ est la plus petite partition symplectique sp\'eciale $\lambda'$ telle que $\lambda\leq \lambda'$. Plus g\'en\'eralement, on a le lemme suivant.

\ass{Lemme}{(i) Soit $(\lambda,\epsilon)\in \boldsymbol{{\cal P}^{symp}}(2n)$, supposons $k_{\lambda,\epsilon}=0$. Notons $sp(\lambda,\epsilon)$ l'unique partition sp\'eciale telle que $symb(\rho_{\lambda,1})$ et $symb(\rho_{sp(\lambda,\epsilon),1})$ soient dans la m\^eme famille. Alors $\lambda\leq sp(\lambda,\epsilon)$.

(ii) Soit $\lambda\in {\cal P}^{symp,sp}(2n)$. Pour  $\epsilon\in \{\pm 1\}^{Jord_{bp}(\lambda)}$, les conditions suivantes sont \'equivalentes:

(a) $k_{\lambda,\epsilon}=0$ et $sp(\lambda,\epsilon)=\lambda$;

(b) $\epsilon$ est constant sur tout $\Delta\in Int(\lambda)$ et, dans le cas o\`u $\Delta_{min}\not=\{0\}$, $\epsilon(i)=1$ pour tout $i\in \Delta_{min}-\{0\}$.

(iii) Soit $\lambda\in {\cal P}^{symp,sp}(2n)$. L'application $\epsilon\mapsto fam\circ symb(\rho_{\lambda,\epsilon})$ est une bijection entre l'ensemble des $\epsilon$ d\'ecrits au (ii) et le sous-ensemble  des $(\tau,\delta)\in({\mathbb Z}/2{\mathbb Z})^{Int(\lambda)}\times ({\mathbb Z}/2{\mathbb Z})^{Int(\lambda)}$ tels que $\delta=0$ et $\tau(\Delta_{min})=0$.}

La preuve est similaire \`a celle du lemme 1.4 ci-dessous.

\bigskip

\subsection{Correspondance de Springer, cas orthogonal impair}

Soit $n\in {\mathbb N}$. On note ${\cal P}^{orth}(2n+1)$ l'ensemble des partitions orthogonales de $2n+1$, c'est-\`a-dire les $\lambda\in {\cal P}(2n+1)$ telles que $mult_{\lambda}(i)$ est pair pour tout entier $i$ pair. Pour une telle partition, on note $Jord_{bp}(\lambda)$ l'ensemble des entiers $i\geq1$ impairs tels que $mult_{\lambda}(i)\geq1$. Plus pr\'ecis\'ement, pour un entier $k\geq1$, on note $Jord_{bp}^k(\lambda)$ l'ensemble des $i\geq1$ impairs tels que $mult_{\lambda}(i)=k$. 

On note $\boldsymbol{{\cal P}^{orth}}(2n+1)$ l'ensemble des couples $(\lambda,\epsilon)$ o\`u $\lambda\in {\cal P}^{orth}(2n+1)$ et $\epsilon\in (\{\pm 1\}^{Jord_{bp}(\lambda)})/\{\pm 1\}$, le groupe $\{\pm 1\}$ s'envoyant diagonalement dans  $\{\pm 1\}^{Jord_{bp}(\lambda)}$. En pratique, on rel\`evera $\epsilon$ en un \'el\'ement de $\{\pm 1\}^{Jord_{bp}(\lambda)}$. Sauf indication contraire, les formules que nous \'ecrirons ne d\'ependront pas du choix de ce rel\`evement.  La correspondance de Springer g\'en\'eralis\'ee   \'etablit une bijection entre $\boldsymbol{{\cal P}^{orth}}(2n+1)$ et l'ensemble des couples $(\rho,k)$ tels que

$k\in {\mathbb N}$, $k$ est impair et $k^2\leq 2n+1$;

$\rho\in \hat{W}_{n- (1-k^2)/2}$.

On note $(\rho_{\lambda,\epsilon},k_{\lambda,\epsilon})$ le couple associ\'e \`a un \'el\'ement $(\lambda,\epsilon)\in \boldsymbol{{\cal P}^{orth}}(2n+1)$. L'entier $k_{\lambda,\epsilon}$ se calcule de la fa\c{c}on suivante. Notons $i_{1}>...>i_{m}$ les entiers impairs $i$ tels que $mult_{\lambda}(i)$ soit impair.  Posons 
$$h= \sum_{j=1,...,m}(-1)^j(1-\epsilon(i_{j})).$$
Alors $k_{\rho,\epsilon}=\vert h+1\vert $. En particulier, si $\epsilon=1$, c'est-\`a-dire $\epsilon(i)=1$ pour tout $i\in Jord_{bp}(\lambda)$, on a $k_{\lambda,1}=1$ et $\rho_{\lambda,1}\in \hat{W}_{n}$.

Une partition $\lambda=(\lambda_{1}\geq \lambda_{2}\geq...)\in {\cal P}^{orth}(2n+1)$ est sp\'eciale si et seulement si $\lambda_{1}$ est impair et
$\lambda_{2j}$ et $\lambda_{2j+1}$ sont de m\^eme parit\'e pour tout $j\geq1$.  Cela \'equivaut \`a ce que  $^t\lambda$ soit orthogonale. En notant ${\cal P}^{orth,sp}(2n+1)$ l'ensemble des partitions orthogonales sp\'eciales de $2n+1$, l'application $\lambda\mapsto {^t\lambda}$  est une involution de ${\cal P}^{orth,sp}(2n+1)$. Consid\'erons une partition $\lambda\in {\cal P}^{orth,sp}(2n+1)$ et d\'efinissons $i_{1}>...>i_{m}$ comme ci-dessus.   L'entier $m$ est forc\'ement impair. On appelle intervalle de $\lambda$ un ensemble $\Delta$ de l'une des formes suivantes:

pour un entier $h=0,...,(m-1)/2$, $\Delta$ est l'ensemble des $i\geq1$ tels que    $mult_{\lambda}(i)\geq1$ et tels que $i_{2h}\geq i \geq i_{2h+1}$, avec la convention $i_{0}=\infty$;

$\Delta=\{i\}$ o\`u $i$ est un entier impair tel que     $mult_{\lambda}(i)\geq1$ et tel qu'il n'existe pas d'entier $h=0,...,(m-1)/2$ de sorte que $i_{2h}\geq i\geq i_{2h+1}$.

Parce que $\lambda$ est sp\'eciale, on v\'erifie que les intervalles sont form\'es d'entiers impairs.  Ils forment une partition de l'ensemble $Jord_{bp}(\lambda)$. On ordonne les intervalles comme dans le cas symplectique. On note $\Delta_{min}$, resp. $\Delta_{max}$, le plus petit, resp. grand, intervalle. On note $Int(\lambda)$ l'ensemble des intervalles.

 L'application $\lambda\mapsto symb(\rho_{\lambda,1})$ est une bijection de ${\cal P}^{orth,sp}(2n+1)$ sur l'ensemble des symboles sp\'eciaux de rang $n$ et de d\'efaut $1$. Pour $\lambda\in {\cal P}^{orth,sp}(2n+1)$, on a d\'efini en \cite{W2} VIII.19 une bijection $fam$ entre la famille du symbole $symb(\rho_{\lambda,1})$ et l'ensemble  des $(\tau,\delta)\in({\mathbb Z}/2{\mathbb Z})^{Int(\lambda)}\times ({\mathbb Z}/2{\mathbb Z})^{Int(\lambda)}$ tels que $\tau(\Delta_{max})=\delta(\Delta_{min})=0$.

Soit $\lambda\in {\cal P}^{orth}(2n+1)$. Il existe une unique partition sp\'eciale $sp(\lambda)\in {\cal P}^{orth,sp}(2n+1)$ telle que $symb(\rho_{\lambda,1})$ et $symb(\rho_{sp(\lambda),1})$ soient dans la m\^eme famille. Il est connu que $\lambda\leq sp(\lambda)$ et que $sp(\lambda)$ est la plus petite partition  orthogonale sp\'eciale $\lambda'$ telle que $\lambda\leq \lambda'$. Plus g\'en\'eralement, on a le lemme suivant.

\ass{Lemme}{(i) Soit $(\lambda,\epsilon)\in \boldsymbol{{\cal P}^{orth}}(2n+1)$, supposons $k_{\lambda,\epsilon}=1$. Notons $sp(\lambda,\epsilon)$ l'unique partition sp\'eciale telle que $symb(\rho_{\lambda,1})$ et $symb(\rho_{sp(\lambda,\epsilon),1})$ soient dans la m\^eme famille. Alors $\lambda\leq sp(\lambda,\epsilon)$.

(ii) Soit $\lambda\in {\cal P}^{orth,sp}(2n+1)$. Pour  $\epsilon\in \{\pm 1\}^{Jord_{bp}(\lambda)}$, les conditions suivantes sont \'equivalentes:

(a) $k_{\lambda,\epsilon}=1$ et $sp(\lambda,\epsilon)=\lambda$;

(b) $\epsilon$ est constant sur tout $\Delta\in Int(\lambda)$.

(iii) Soit $\lambda\in {\cal P}^{orth,sp}(2n+1)$. L'application $\epsilon\mapsto fam\circ symb(\rho_{\lambda,\epsilon})$ est une bijection entre l'ensemble des $\epsilon$ d\'ecrits au (ii), modulo le groupe diagonal $\{\pm 1\}$,  et le sous-ensemble des $(\tau,\delta)\in({\mathbb Z}/2{\mathbb Z})^{Int(\lambda)}\times ({\mathbb Z}/2{\mathbb Z})^{Int(\lambda)}$  tels que $\delta=0$ et $\tau(\Delta_{max})=0$.} 

Preuve. Soit $(\lambda,\epsilon)$ comme en (i). On peut supposer $\lambda=(\lambda_{1},...,\lambda_{2n+1})$. Dans la suite $\lambda+\{2n,...,0\}$, il y  a $n$ nombres pairs not\'es $2z_{1}>...>2z_{n}$ et $n+1$ nombres impairs not\'es $2z'_{1}+1>...>2z'_{n+1}+1$.  On note $z=(z_{1},...,z_{n})$ et $z'=(z'_{1},...,z'_{n+1})$ puis $A^{\sharp}=z'+\{n,...,0\}=(a_{1}^{\sharp},...,a_{n+1}^{\sharp})$, $B^{\sharp}=z+\{n-1,...,0\}=(b_{1}^{\sharp},...,b_{n}^{\sharp})$. On v\'erifie que
$$(1) \qquad a_{1}^{\sharp}\geq b_{1}^{\sharp}\geq a_{2}^{\sharp}\geq  ... \geq b_{n}^{\sharp}\geq a_{n+1}^{\sharp}.$$
On voit aussi qu'il y a une unique bijection croissante $i\mapsto \Sigma_{i}$ entre l'ensemble $Jord_{bp}(\lambda)$ et celui des sous-ensembles non vides de $(A^{\sharp}\cup B^{\sharp})-(A^{\sharp}\cap B^{\sharp})$ form\'es d'entiers cons\'ecutifs, et maximaux pour cette propri\'et\'e. Posons
$$A=(A^{\sharp}-\bigcup_{i\in Jord_{bp}(\lambda),\epsilon_{i}=-1}(\Sigma_{i}\cap A^{\sharp}))\cup(\bigcup_{i\in Jord_{bp}(\lambda),\epsilon_{i}=-1}(\Sigma_{i}\cap B^{\sharp})),$$
$$B=(B^{\sharp}-\bigcup_{i\in Jord_{bp}(\lambda),\epsilon_{i}=-1}(\Sigma_{i}\cap B^{\sharp}))\cup(\bigcup_{i\in Jord_{bp}(\lambda),\epsilon_{i}=-1}(\Sigma_{i}\cap A^{\sharp})).$$
Si on multiplie $\epsilon$ par l'\'el\'ement diagonal $-1$, on \'echange $A$ et $B$. On peut donc choisir le rel\`evement $\epsilon$ de sorte que $\vert A\vert \geq \vert B\vert $. L'hypoth\`ese $k_{\lambda,\epsilon}=1$ entra\^{\i}ne alors que $\vert A\vert =n+1$ et $\vert B\vert =n$. On d\'efinit les suites $X=(x_{1},...,x_{n+1})$ et $Y=(y_{1},...,y_{n})$ par $X+\{n,...,0\}=A$ et $Y+\{n-1,...,0\}=B$. Alors $symb(\rho_{\lambda,\epsilon})=(X,Y)$. 

Soit $k\in \{1,...,2n+1\}$. On va majorer $S_{k}(\lambda)$ en fonction de $(X,Y)$.  Tout d'abord
$$(2) \qquad S_{k}(\lambda)=S_{k}(\lambda+\{2n,...,0\})-\frac{k(4n+1-k)}{2}.$$
Notons $j_{1}<...< j_{s}$ les indices $j$ tels que $\lambda_{j}$ soit impair et $h_{1}<...<h_{t}$ les indices $h$ pour lesquels $\lambda_{h}$ est pair. Puisque la somme des $\lambda_{j}$ vaut $2n+1$ qui est impair, $s$ est impair et $t=2n+1-s$ est pair. Puisque  tout nombre pair non nul  intervient avec multiplicit\'e paire, la parit\'e de $t$ entra\^{\i}ne que $0$ intervient aussi avec multiplicit\'e paire. Il en r\'esulte aussi que $h_{2r}=h_{2r-1}+1$ pour tout $r=1,...,t/2$. Pour $r=1,...s$,  il y a $j_{r}-r$ termes $\lambda_{j}$ qui sont pairs et strictement sup\'erieurs \`a $\lambda_{j_{r}}$. Puisque ces termes interviennent avec multiplicit\'e paire,   $j_{r}$ est de m\^eme parit\'e que $r$. Soient $s_{k}\in \{1,...,s\}$ et $t_{k}\in \{1,...,t\}$ les plus grands entiers tels que $j_{s_{k}}\leq k$ et $h_{t_{k}}\leq k$. On a $s_{k}+t_{k}=k$.   Les $k$ premiers termes de $\lambda+\{2n,...,0\}$ sont les
 
(3) $\lambda_{j_{u}}+2n+1-j_{u}$ pour $u=1,...,s_{k}$

\noindent et les

(4) $\lambda_{h_{v}}+2n+1-h_{v}$ pour $v=1,...,t_{k}$. 

D'apr\`es les propri\'et\'es de nos suites, il y a $[(s_{k}+1)/2]$ \'el\'ements impairs et $[s_{k}/2]$ \'el\'ements pairs parmi les \'el\'ements  (3). Si $t_{k}$ est pair, il y a $t_{k}/2$ \'el\'ements impairs et $t_{k}/2$ \'el\'ements pairs parmi les \'el\'ements  (4). Si $t_{k}$ est impair et $h_{t_{k}}$ est pair, il y a $(t_{k}+1)/2$ \'el\'ements impairs et $(t_{k}-1)/2$ \'el\'ements pairs parmi les \'el\'ements (4). Si $t_{k}$ est impair et $h_{t_{k}}$ est impair, il y a $(t_{k}-1)/2$ \'el\'ements impairs et $(t_{k}+1)/2$ \'el\'ements pairs parmi les \'el\'ements  (4).  En r\'eunissant les deux types d'\'el\'ements, on voit que, parmi les $k$ premiers termes de $\lambda+\{2n,...,0\}$, il y a $[(k+1)/2]+\eta$ termes impairs et $[k/2]-\eta$ termes pairs, o\`u

$\eta=1$ si $t_{k}$ et $s_{k}$ sont impairs et $h_{t_{k}}$ est pair,

$\eta=-1$ si $t_{k}$ et $h_{t_{k}}$ sont impairs et $s_{k}$ est pair,

$\eta=0$ dans les autres cas.

Les $k$ premiers termes de $\lambda+\{2n,...,0\}$ sont donc $2z'_{1}+1,...,2z'_{[(k+1)/2]+\eta}+1$ et $2z_{1},...,2z_{[k/2]-\eta}$. 

Supposons $\lambda_{k}$ impair. Alors $j_{s_{k}}=k$. On a dit que $s_{k}$ est de la m\^eme parit\'e que $j_{s_{k}}$, donc que $k$, donc $t_{k}=k-s_{k}$ est pair. Donc $\eta=0$ et il r\'esulte de la description ci-dessus que
$$(5) \qquad S_{k}(\lambda+\{2n,...,0\})=2S_{[(k+1)/2]}(z')+[(k+1)/2]+2S_{[k/2]}(z).$$

Supposons $\lambda_{k}$ pair. Alors $h_{t_{k}}=k$. Si $\eta=0$, le calcul est le m\^eme que ci-dessus et on a (5). Supposons $\eta=1$. Alors $k=h_{t_{k}}$ est pair et les $k$ premiers termes de $\lambda+\{2n,...,0\}$ sont $2z'_{1}+1,...,2z'_{k/2+1}+1$ et $2z_{1},...,2z_{k/2-1}$. Le dernier de ces termes  est  $\lambda_{k}+2n+1-k$ qui est impair,   donc c'est $2z'_{k/2+1}+1$. Mais, $t_{k}$ \'etant impair, on a $h_{t_{k}+1}=h_{t_{k}}+1=k+1$ et $\lambda_{k}=\lambda_{k+1}$. Le $k+1$-i\`eme terme  de $\lambda+\{2n,...,0\}$ est $\lambda_{k}+2n-k$ qui est pair, c'est donc le premier terme pair strictement inf\'erieur \`a $2z_{k/2-1}$, autrement dit, c'est $2z_{k/2}$. Les \'egalit\'es $\lambda_{k}+2n+1-k =2z'_{k/2+1}+1$ et  $\lambda_{k}+2n-k=2z_{k/2}$ entra\^{\i}nent $z'_{k/2+1}=z_{k/2}$. Les $k$ premiers termes de $\lambda+\{2n,...,0\}$ sont donc aussi bien  $2z'_{1}+1,...,2z'_{ k/2}+1$ et $2z_{1},...,2z_{k/2-1},2z_{k/2}+1$. On obtient alors
$$(6) \qquad S_{k}(\lambda+\{2n,...,0\})=2S_{[(k+1)/2]}(z')+ [(k+1)/2]+1+2S_{[k/2]}(z).$$
Supposons maintenant $\eta=-1$. Alors $k=h_{t_{k}}$ est impair et les $k$ premiers termes de $\lambda+\{2n,...,0\}$ sont $2z'_{1}+1,...,2z'_{ (k-1)/2}+1$ et $2z_{1},...,2z_{( k+1)/2}$. Le dernier de ces termes  est  $\lambda_{k}+2n+1-k$ qui est pair,   donc c'est $2z_{(k+1)/2}$. Comme ci-dessus, on a $\lambda_{k}=\lambda_{k+1}$. Le $k+1$-i\`eme terme  de $\lambda+\{2n,...,0\}$ est $\lambda_{k}+2n-k$ qui est impair, c'est donc le premier terme impair strictement inf\'erieur \`a $2z'_{(k-1)/2}+1$, autrement dit, c'est $2z'_{(k+1)/2}+1$. Les \'egalit\'es $\lambda_{k}+2n+1-k =2z_{(k+1)/2}$ et  $\lambda_{k}+2n-k=2z'_{(k+1)/2}+1$ entra\^{\i}nent $z'_{(k+1)/2}+1=z_{(k+1)/2}$. Les $k$ premiers termes de $\lambda+\{2n,...,0\}$ sont donc aussi bien  $2z'_{1}+1,...,2z'_{ (k-1)/2}+1, 2z'_{(k+1)/2}+2$ et $2z_{1},...,2z_{ (k-1)/2}$. On obtient encore (6). 

Supposons encore $\lambda_{k}$ pair.  Puisque $h_{t_{k}}=k=s_{k}+t_{k}$, les conditions de parit\'e sur $t_{k}$ sont redondantes dans la d\'efinition de $\eta$. On voit que $\eta=\pm 1$ si et seulement si $k+s_{k}$ est impair. On voit aussi que $s_{k}$ est de la m\^eme parit\'e que $S_{k}(\lambda)$. On a obtenu que, si $S_{k}(\lambda)+k$ est pair, on a la formule (5) tandis que, si $S_{k}(\lambda)+k$ est impair, on a la formule (6). Posons alors, pour $k=1,...,2n+1$,

$\nu_{k}(\lambda)=1$ si $\lambda_{k}$ est pair et $S_{k}(\lambda)+k$ est impair, $\nu_{k}(\lambda)=0$ sinon. 

On a la formule g\'en\'erale
$$(7) \qquad S_{k}(\lambda+\{2n,...,0\})=2S_{[(k+1)/2]}(z')+[(k+1)/2]+\nu_{k}(\lambda)+2S_{[k/2]}(z).$$
Remarquons que, d'apr\`es le calcul ci-dessus, on a

(8) si $\nu_{k}(\lambda)=1$, alors $k\leq 2n$ et $\lambda_{k+1}=\lambda_{k}$.

D'apr\`es la d\'efinition des termes $A^{\sharp}$ et $B^{\sharp}$, on a les \'egalit\'es
$$S_{[(k+1)/2]}(z')=S_{[(k+1)/2]}(A^{\sharp})-[(k+1)/2](2n+1-[(k+1)/2])/2,$$
$$S_{[(k-1)/2]}(z)=S_{[(k-1)/2]}(B^{\sharp})-[(k-1)/2](2n-1-[(k-1)/2])/2.$$
Avec les formules (2) et (7), on obtient
$$S_{k}(\lambda)=2S_{[(k+1)/2]}(A^{\sharp})+2S_{[(k-1)/2]}(B^{\sharp})+\nu_{k}(\lambda)+c_{k},$$
o\`u $c_{k}$ est un nombre qui ne d\'epend pas de $\lambda$. Notons $A^{\sharp}\sqcup B^{\sharp}$ la r\'eunion des suites $A^{\sharp}$ et $B^{\sharp}$, les termes \'etant rang\'es en ordre d\'ecroissant mais compt\'es avec leur multiplicit\'e (c'est-\`a-dire qu'un terme intervenant dans les deux suites intervient avec multiplicit\'e $2$). La propri\'et\'e (1) entra\^{\i}ne que la r\'eunion (en ce sens) des $[(k+1)/2]$ plus grands termes de $A^{\sharp}$ et des $[(k-1)/2]$ plus grands termes de $B^{\sharp}$ n'est autre que la suite des $k$ plus grands termes de $A^{\sharp}\sqcup B^{\sharp}$. La formule pr\'ec\'edente se r\'ecrit
$$S_{k}(\lambda)=2S_{k}(A^{\sharp}\sqcup B^{\sharp})+\nu_{k}(\lambda)+c_{k}.$$
Avec une d\'efinition similaire, on a $A^{\sharp}\sqcup B^{\sharp}=A\sqcup B$, donc aussi $S_{k}(A^{\sharp}\sqcup B^{\sharp})=S_{k}(A\sqcup B)$. Il existe deux entiers $e,f\in {\mathbb N}$ tels que $e+f=k$ et que la famille des $k$ plus grands \'el\'ements de $A\sqcup B$ soit la r\'eunion des familles des $e$ plus grands \'el\'ements de $A$ et des $f$ plus grands \'el\'ements de $B$. Alors
$$(9) \qquad S_{k}(A\sqcup B)=S_{e}(A)+S_{f}(B).$$
Par d\'efinition de $X$ et $Y$, on a
$$S_{e}(A)=S_{e}(X)+e(2n+2-e)/2,\,\,S_{f}(B)=S_{f}(Y)+f(2n-f)/2.$$
D'o\`u
$$S_{k}(\lambda)=2S_{e}(X)+2S_{f}(Y)+\nu_{k}(\lambda)-e(2-e)-f^2+c'_{k},$$
o\`u $c'_{k}=c_{k}+2nk$ est ind\'ependant de $\lambda$. Le terme $S_{e}(X)+S_{f}(Y)$ est la somme de $k$ termes de la famille $X\sqcup Y$, donc il est major\'e par la somme des $k$ plus grands termes de cette famille:
$$(10) \qquad S_{e}(X)+S_{f}(Y)\leq S_{k}(X\sqcup Y).$$
D'o\`u
$$(11) \qquad S_{k}(\lambda)\leq 2S_{k}(X\sqcup Y)+\nu_{k}(\lambda)-e(e-2)-f^2+c'_{k}.$$

Posons $\underline{\lambda}=sp(\lambda,\epsilon)$. Reprenons le calcul en rempla\c{c}ant $(\lambda,\epsilon)$ par $(\underline{\lambda},1)$. On souligne les objets associ\'es \`a cette paire. Parce que le caract\`ere $\epsilon$ est remplac\'e par $1$, on a les \'egalit\'es $\underline{A}=\underline{A}^{\sharp}$, $\underline{B}=\underline{B}^{\sharp}$ et l'on voit que $\underline{e}=[(k+1)/2]$ et $\underline{f}=[k/2]$. Parce que le symbole $(\underline{X},\underline{Y})$ est sp\'ecial, la r\'eunion des $[(k+1)/2]$ plus grands termes de $\underline{X}$ et des $[k/2]$ plus grands termes de $\underline{Y}$ n'est autre que la famille des $k$ plus grands termes de $\underline{X}\sqcup \underline{Y}$. L'analogue de l'in\'egalit\'e (10) est donc une \'egalit\'e et on obtient
$$S_{k}(\underline{\lambda})=2S_{k}(\underline{X}\sqcup \underline{Y})+\nu_{k}(\underline{\lambda})-[(k+1)/2]([(k+1)/2]-2)-[k/2]^2+c'_{k}.$$
Par d\'efinition de $sp(\lambda,\epsilon)$, les symboles $(X,Y)$ et $(\underline{X},\underline{Y})$ sont dans la m\^eme famille, d'o\`u $X\sqcup Y=\underline{X}\sqcup \underline{Y}$. En comparant (11) avec l'\'egalit\'e ci-dessus, on obtient
$$(12) \qquad S_{k}(\lambda)\leq S_{k}(\underline{\lambda})+\nu_{k}(\lambda)-e(e-2)-f^2-\nu_{k}(\underline{\lambda}) +[(k+1)/2]([(k+1)/2]-2)+[k/2]^2.$$
On v\'erifie que, pour deux entiers $e,f$ tels que $e+f=k$, on a $e(e-2)+f^2\geq [(k+1)/2]([(k+1)/2]-2)+[k/2]^2$, l'\'egalit\'e n'\'etant v\'erifi\'ee que pour les couples $(e,f)=([(k+1)/2],[k/2])$ ou, si $k$ est pair, $(e,f)=(k/2+1,k/2-1)$. On obtient 
$$(13) \qquad S_{k}(\lambda)\leq S_{k}(\underline{\lambda})+\nu_{k}(\lambda) -\nu_{k}(\underline{\lambda}).$$
Supposons $S_{k}(\lambda)>S_{k}(\underline{\lambda})$. L'in\'egalit\'e pr\'ec\'edente force $\nu_{k}(\lambda)=1$. Donc $\lambda_{k}$ est pair, $k+S_{k}(\lambda)$ est impair et,  d'apr\`es (8), $\lambda_{k+1}=\lambda_{k}$. Les entiers $k-1+S_{k-1}(\lambda)$ et $k+1+S_{k+1}(\lambda)$ sont pairs, donc $\nu_{k-1}(\lambda)=\nu_{k+1}(\lambda)=0$. L'in\'egalit\'e (13) entra\^{\i}ne donc $S_{k-1}(\lambda)\leq S_{k-1}(\underline{\lambda})$ et $S_{k+1}(\lambda)\leq S_{k+1}(\underline{\lambda})$ (pour \^etre pr\'ecis, si $k=1$, notre calcul ne s'applique pas \`a $k-1$ mais, dans ce cas, l'in\'egalit\'e $S_{0}(\lambda)\leq S_{0}(\underline{\lambda})$ est triviale). Les deux in\'egalit\'es $S_{k-1}(\lambda)\leq S_{k-1}(\underline{\lambda})$ et $S_{k}(\lambda)>S_{k}(\underline{\lambda})$ entra\^{\i}nent $\lambda_{k}>\underline{\lambda}_{k}$. Donc $\lambda_{k+1}=\lambda_{k}>\underline{\lambda}_{k}\geq \underline{\lambda}_{k+1}$. Alors l'in\'egalit\'e $S_{k}(\lambda)>S_{k}(\underline{\lambda})$ entra\^{\i}ne $S_{k}(\lambda)>S_{k}(\underline{\lambda})$, contrairement \`a ce que l'on a vu ci-dessus. Cette contradiction prouve l'in\'egalit\'e $S_{k}(\lambda)\leq S_{k}(\underline{\lambda})$. Cela \'etant vrai pour tout $k$, on conclut $\lambda\leq \underline{\lambda}$, ce qui d\'emontre  le (i) de l'\'enonc\'e. 

Supposons maintenant $\lambda$ sp\'eciale et $\lambda=\underline{\lambda}$. Consid\'erons l'in\'egalit\'e (12) pour $k$ impair. Les termes relatifs \`a $\lambda$ et $\underline{\lambda}$ s'annulent et il reste
$$e(e-2)+f^2\leq ((k+1)/2)((k+1)/2-2)+(k/2)^2.$$
Comme on l'a dit, cela entra\^{\i}ne que $e=(k+1)/2$ et $f=k/2$.  
Parce que $\lambda$ est sp\'eciale, on calcule facilement les termes $A^{\sharp}$ et $B^{\sharp}$. Puisque $\lambda_{1}$ est impair, le terme $\lambda_{1}+2n$ l'est aussi, donc c'est $2z'_{1}+1$. Pour $h=1,...,n$, les termes $\lambda_{2h}$ et $\lambda_{2h+1}$ sont de m\^eme parit\'e donc les termes $\lambda_{2h}+2n+1-2h$ et $\lambda_{2h+1}+2n-2h $ sont de parit\'e oppos\'ee. Par r\'ecurrence, ce sont les termes $2z_{h}$ et $2z'_{h+1}+1$. Cela permet le calcul des termes $z_{h}$ et $z'_{h+1}$, puis des termes $a^{\sharp}_{h+1}$ et $b^{\sharp}_{h}$. On obtient

 $a_{1}^{\sharp}=(\lambda_{1}-1)/2+2n,$
 
 pour $h=1,...,n$, $b^{\sharp}_{h}=a^{\sharp}_{h+1}=\lambda_{2h}/2+2n-2h$ si $\lambda_{2h}=\lambda_{2h+1}$ est pair, $b^{\sharp}_{h}=(\lambda_{2h}+1)/2+2n-2h$ et $a^{\sharp}_{h+1}=(\lambda_{2h+1}-1)/2+2n-2h$ si $\lambda_{2h}$ et $\lambda_{2h+1}$ sont impairs. 
 
 Consid\'erons l'intervalle maximal $\Delta_{max}$ de $\lambda$. Notons ses \'el\'ements $i_{1}>...>i_{t}$. Ils sont impairs. Les multiplicit\'es de $i_{1},...,i_{t-1}$ sont paires et celle de $i_{t}$ est impaire. On note ces multiplicit\'es $2m_{1},...,2m_{t-1},2m_{t}+1$. L'intervalle correspond aux \'el\'ements suivants de $A^{\sharp}\sqcup B^{\sharp}$:
 
 $$a_{1}^{\sharp}>b_{1}^{\sharp}>...>a_{m_{1}}^{\sharp}>b_{m_{1}}^{\sharp}>a_{m_{1}+1}^{\sharp}>b_{m_{1}+1}^{\sharp}>...>a_{m_{\leq 2}}^{\sharp}> b_{m_{\leq 2}}^{\sharp}>...$$
 $$>a^{\sharp}_{m_{\leq t-1}+1}>b^{\sharp}_{m_{\leq t-1}+1}>...>b^{\sharp}_{m_{\leq t}}>a^{\sharp}_{m_{\leq t}+1}$$
 (on rappelle que $m_{\leq i}=m_{1}+...+m_{i}$).
 Supposons que $\epsilon$ ne vaut pas $1$ sur $\Delta_{max}$. Soit $s$ le plus petit \'el\'ement de $\{1,...,t\}$ tel que $\epsilon(i_{s})=-1$. Appliquons l'\'egalit\'e (9)  \`a $k=2m_{\leq s-1}+1$. Comme on l'a dit plus haut, on a  $e=m_{\leq s-1}+1$, $f=m_{\leq s-1}$. On obtient
$$ a^{\sharp}_{1}+...+a^{\sharp}_{m_{\leq s-1}+1}+b^{\sharp}_{1}+...+b^{\sharp}_{m_{\leq s-1}}=a_{1}+...+a_{m_{s-1}+1}+b_{1}+...+b_{m_{\leq s-1}}.$$
Par construction de $A$ et $B$ et par d\'efinition de $s$, les termes de ces ensembles sont \'egaux \`a ceux de $A^{\sharp}$ et $B^{\sharp}$ jusqu'\`a l'indice $m_{\leq s-1}$ et l'\'egalit\'e pr\'ec\'edente devient
$$a^{\sharp}_{m_{\leq s-1}+1}=a_{m_{\leq s-1}+1}.$$
 Par contre, passer de $(A^{\sharp},B^{\sharp})$ \`a $(A,B)$ \'echange les termes correspondant \`a l'entier $i_{s}$. Hormis le cas $s=t$ et $m_{t}=0$, on a donc $a_{m_{\leq s-1}+1}=b^{\sharp}_{m_{\leq s-1}+1}$. Quand $s=t$ et $m_{t}=0$, $a_{m_{\leq s-1}+1}$ est un terme de la famille $A^{\sharp}\sqcup B^{\sharp}$ qui n'est pas dans l'ensemble \'ecrit ci-dessus. Dans tous les   cas, on obtient $a_{m_{\leq s-1}+1}<a^{\sharp}_{m_{\leq s-1}+1}$ ce qui contredit l'\'egalit\'e de ces termes prouv\'ee ci-dessus. Cette contradiction conclut:  $\epsilon$ vaut $1$ sur $\Delta_{max}$. 
 
 Consid\'erons maintenant un intervalle $\Delta\not=\Delta_{max}$. On note ses \'el\'ements $i_{1}>...>i_{t}$. Le premier indice $j$ tel que $\lambda_{j}=i_{1}$ est forc\'ement pair. Notons le $2u$. Si  $t=1$, la multiplicit\'e de $i_{1}$ est paire. On la note $2m$.    Alors l'intervalle correspond aux \'el\'ements suivants de $A^{\sharp}\sqcup B^{\sharp}$:
 
 $$b^{\sharp}_{u}< a^{\sharp}_{u+1}<...b^{\sharp}_{u+m-1}<a^{\sharp}_{u+m}.$$

Si $t>1$, les multiplicit\'es de $i_{1}$ et $i_{t}$ sont impaires et, pour $1<s<t$,  celle de $i_{s}$ est paire. On les note respectivement $2m_{1}+1$, $2m_{t}+1$, $2m_{s}$. 
Alors l'intervalle correspond aux \'el\'ements suivants de $A^{\sharp}\sqcup B^{\sharp}$:
 
 $$b^{\sharp}_{u}< a^{\sharp}_{u+1}<...b^{\sharp}_{u+m_{1}}<a^{\sharp}_{u+m_{1}+1}<...<b^{\sharp}_{u+m_{\leq2}}$$
 $$<a^{\sharp}_{u+m_{\leq2}+1}<... <b^{\sharp}_{u+m_{\leq t-1}}<a^{\sharp}_{u+m_{\leq t-1}+1}<...< a^{\sharp}_{u+m_{\leq t}+1}.$$
 
Supposons par r\'ecurrence que $\epsilon$ est constant sur tout intervalle strictement sup\'erieur \`a $\Delta$. Pour ces intervalles, ou bien on ne change pas les termes  de $A^{\sharp}$ et  $B^{\sharp}$ leur correspondant, ou bien on les \'echange. Mais on voit ci-dessus que chaque intervalle contribue autant \`a $A_{\sharp}$ qu'\`a $B_{\sharp}$.  Cela ne perturbe pas les num\'erotations des termes post\'erieurs,  on veut dire par l\`a que la contribution \`a $A$, resp. $B$, de l'intervalle $\Delta$ commence par $a_{u+1}$, resp. $b_{u}$. Si $\Delta$ est r\'eduit \`a $i_{1}$, on n'a rien \`a d\'emontrer: $\epsilon$ est forc\'ement constant sur $\Delta$. Supposons $t>1$ et que $\epsilon$ ne soit pas constant sur $\Delta$. Notons $s$ le plus petit \'el\'ement de $\{2,...,t\}$ tel que $\epsilon_{i_{s-1}}\not=\epsilon_{i_{s}}$. Appliquons l'\'egalit\'e (9) \`a $k=2u-1$ (on sait qu'alors $e=u$ et $f=u-1$) et \`a $k=2u+2m_{\leq s-1}+1$ (on sait qu'alors $e=u+m_{\leq s-1}+1$ et $f=u+m_{\leq s-1}$). Par diff\'erence, on obtient
$$(14) \qquad b^{\sharp}_{u}+ a^{\sharp}_{u+1}+...+b_{u+m_{\leq s-1}}^{\sharp}+a^{\sharp}_{u+m_{\leq s-1}+1}=b_{u}+ a_{u+1}+...+b_{u+m_{\leq s-1}}+a_{u+m_{\leq s-1}+1}.$$

Supposons  d'abord $\epsilon(i_{s-1})=1$ et $\epsilon(i_{s})=-1$ . Les termes de l'ensemble $A$ sont \'egaux \`a ceux de $A^{\sharp}$ entre les indices $u+1$ et $u+m_{\leq s-1}$. Les termes de l'ensemble  $B$ sont \'egaux \`a ceux de $B^{\sharp}$ entre les indices $u$ et $u+m_{\leq s-1}$. L'\'egalit\'e (14) devient
$$a^{\sharp}_{u+m_{\leq s-1}+1}=a_{u+m_{\leq s-1}+1}.$$
Parce que $\epsilon_{i_{s}}=-1$, on \'echange les termes correspondant \`a l'entier $i_{s}$. Hormis le cas $s=t$ et $m_{t}=0$, on a donc
$$a_{u+m_{\leq s-1}+1}=b^{\sharp}_{u+m_{\leq s-1}+1}.$$
Sii $s=t$ et $m_{t}=0$, $a_{u+ m_{\leq s-1}+1}$ est un terme de la famille $A^{\sharp}\sqcup B^{\sharp}$ qui est au-del\`a de ceux   \'ecrits ci-dessus. Dans tous les   cas, on obtient $a_{u+m_{\leq s-1}+1}<a^{\sharp}_{u+m_{\leq s-1}+1}$ ce qui contredit l'\'egalit\'e de ces termes prouv\'ee ci-dessus. 

Supposons maintenant $\epsilon(i_{s-1})=-1$ et $\epsilon(i_{s})=1$. Les entiers $i_{1},...,i_{s-1}$ contribuent \`a $A$ et $B$ en \'echangeant leur contribution \`a $A^{\sharp}$ et $B^{\sharp}$. D'o\`u
$$a_{u+1}=b^{\sharp}_{u},... ,a_{u+m_{\leq s-1}+1}=b^{\sharp}_{u+m_{\leq s-1}},$$
$$b_{u}=a_{u+1}^{\sharp},...,b_{u+m_{\leq s-1}-1}=a_{u+m_{\leq s-1}}^{\sharp}.$$
L'\'egalit\'e (14) devient
$$b_{u+m_{\leq s-1}}=a^{\sharp}_{u+m_{\leq s-1}+1}.$$
Par contre, l'entier $i_{s}$ contribue par les m\^emes termes \`a $A$ et $A^{\sharp}$ comme \`a $B$ et $B^{\sharp}$. Mais les indices sont d\'ecal\'es et on a
$$a_{u+m_{\leq s-1}+2}=a^{\sharp}_{u+m_{\leq s-1}+1}$$
et, hormis le cas $s=t$ et $m_{t}=0$, $b_{u+m_{\leq s-1}}=b^{\sharp}_{u+m_{\leq s-1}+1}$. Si $s=t$ et $m_{t}=0$, $b_{u+ m_{\leq s-1}}$ est un terme de la famille $A^{\sharp}\sqcup B^{\sharp}$ qui est au-del\`a de ceux   \'ecrits ci-dessus. Dans tous les   cas, on obtient $b_{u+m_{\leq s-1}}<a^{\sharp}_{u+m_{\leq s-1}+1}$ ce qui contredit l'\'egalit\'e de ces termes prouv\'ee ci-dessus.

Ces contradictions prouvent que $\epsilon$ est constant sur $\Delta$. Cela prouve que, sous les hypoth\`eses du (ii) de l'\'enonc\'e, la condition (a) implique (b).

Soit maintenant $(\lambda,\epsilon)\in \boldsymbol{{\cal P}}^{orth}(2n+1)$, supposons $\lambda$ sp\'eciale et $\epsilon$ constant sur les intervalles de $\lambda$. En notant $i_{1}>...>i_{m}$ les \'el\'ements de $Jord_{bp}(\lambda)$ intervenant avec multiplicit\'e impaire, on a $\epsilon(i_{2h})=\epsilon(i_{2h+1})$ pour tout $h=1,...,(m-1)/2$. La recette indiqu\'ee plus haut pour calculer $k_{\lambda,\epsilon}$ montre que cet entier vaut $1$. Relevons $\epsilon$ en l'\'el\'ement de $\{\pm 1\}^{Jord_{bp}(\lambda)}$ qui vaut $1$ sur le plus grand intervalle. On a calcul\'e ci-dessus les termes $A^{\sharp}$, $B^{\sharp}$, $A$, $B$. Remarquons que les deux premiers sont aussi les termes $\underline{A}$ et $\underline{B}$ associ\'es \`a $(\lambda,1)$. On en d\'eduit facilement les termes $\underline{X}$, $\underline{Y}$, $X$ et $Y$. On voit que les ensembles suivants contribuent  de la m\^eme fa\c{c}on \`a $\underline{X}$ et $X$ comme \`a $\underline{Y}$ et $Y$:

l'intervalle maximal;  sa contribution est de la forme
$$x_{1}= y_{1}>x_{2}=y_{2}>... >x_{u}= y_{u}>x_{u+1};$$

tout couple $\lambda_{2h},\lambda_{2h+1}$ d'\'el\'ements pairs donc \'egaux; sa contribution    est de la forme $y_{h}=x_{h+1}$;

tout  intervalle non maximal sur lequel $\epsilon$ vaut $1$; sa contribution est de la forme
$$y_{u}>x_{u+1}=y_{u+1}>...>x_{v}=y_{v}>x_{v+1}.$$

Par contre, la contribution \`a $\underline{X}$ et $\underline{Y}$ d'un intervalle sur lequel $\epsilon$ vaut $-1$ est de la forme
$$\underline{y}_{u}>\underline{x}_{u+1}=\underline{y}_{u+1}>...>\underline{x}_{v}=\underline{y}_{v}>\underline{x}_{v+1},$$
tandis que sa contribution \`a $X$ et $Y$ est
$$x_{u+1}=\underline{y}_{u}>y_{u}=\underline{x}_{u+1}=x_{u+2}=\underline{y}_{u+1}>...>y_{v-1}=\underline{x}_{v}=x_{v+1}=\underline{y}_{v}> y_{v}=\underline{x}_{v+1}.$$

Il est imm\'ediat que $X\sqcup Y=\underline{X}\sqcup \underline{Y}$, autrement dit les symboles $(X,Y)$ et $(\underline{X},\underline{Y})$ sont dans la m\^eme famille. Cela prouve que $\underline{\lambda}=sp(\lambda,\epsilon)$, donc que la relation (b) du (ii) de l'\'enonc\'e entra\^{\i}ne la relation (a). 

Conservons les hypoth\`eses sur $(\lambda,\epsilon)$. On rel\`eve $\epsilon$ comme ci-dessus. Posons
  $(\tau,\delta)=fam\circ symb(\rho_{\lambda,\epsilon})$. En utilisant la description du symbole $(X,Y)$ faite ci-dessus et la d\'efinition de l'application $fam$ de \cite{W2} VIII.19, on calcule, pour tout intervalle $\Delta$:
  
  $\delta(\Delta)=0$;
  
  $\tau(\Delta)=0$ si $\epsilon$ vaut $1$ sur $\Delta$ et $\tau(\Delta)=1$ si $\epsilon$ vaut $-1$.
  
  Le (iii) de l'\'enonc\'e s'en d\'eduit. $\square$
  
  \bigskip

\subsection{Correspondance de Springer, cas orthogonal pair}

Soit $n\in {\mathbb N}$. On note ${\cal P}^{orth}(2n)$ l'ensemble des partitions orthogonales de $2n$, c'est-\`a-dire les $\lambda\in {\cal P}(2n)$ telles que $mult_{\lambda}(i)$ est pair pour tout entier $i$ pair. Pour une telle partition, on note $Jord_{bp}(\lambda)$ l'ensemble des entiers $i\geq1$ impairs tels que $mult_{\lambda}(i)\geq1$. Plus pr\'ecis\'ement, pour un entier $k\geq1$, on note $Jord_{bp}^k(\lambda)$ l'ensemble des $i\geq1$ impairs tels que $mult_{\lambda}(i)=k$. 

Pour $\lambda\in {\cal P}^{orth}(2n)$, disons que $\lambda$ est exceptionnelle si $Jord_{bp}(\lambda)=\emptyset$. Si $n>0$, on introduit l'ensemble $\underline{{\cal P}}^{orth}(2n)$ form\'e des partitions $\lambda\in {\cal P}^{orth}(2n)$ non exceptionnelles et des paires $(\lambda,+)$ et $(\lambda,-)$ pour les partitions $\lambda\in {\cal P}^{orth}(2n)$ exceptionnelles.

Justifions  cette d\'efinition. Notons $\bar{{\mathbb F}}_{q}$ une cl\^oture alg\'ebrique de ${\mathbb F}_{q}$ et ${\bf O}(2n)$ le groupe orthogonal  \'evident sur $\bar{{\mathbb F}}_{q}$.  L'ensemble ${\cal P}^{orth}(2n)$ param\`etre les classes de conjugaison unipotentes  par ${\bf O}(2n)(\bar{{\mathbb F}}_{q})$ dans ${\bf SO}(2n)(\bar{{\mathbb F}}_{q})$. Mais il arrive que de telles classes se coupent en deux classes de conjugaison par ${\bf SO}(2n)(\bar{{\mathbb F}}_{q})$. Cela arrive pr\'ecis\'ement quand la classe est param\'etr\'ee par une partition $\lambda$  exceptionnelle. Alors l'ensemble $\underline{{\cal P}}^{orth}(2n)$ param\`etre les classes de conjugaison unipotentes  par ${\bf SO}(2n)(\bar{{\mathbb F}}_{q})$ dans ${\bf SO}(2n)(\bar{{\mathbb F}}_{q})$. Si $n=0$, on pose $\underline{{\cal P}}^{orth}(0)={\cal P}^{orth}(0)=\{\emptyset\}$. Il y a en tout cas une application \'evidente de $\underline{{\cal P}}^{orth}(2n)$ dans ${\cal P}^{orth}(2n)$. Si $\underline{\lambda}$ est un \'el\'ement de $\underline{{\cal P}}^{orth}(2n)$, on note sans plus de commentaire $\lambda\in{\cal P}^{orth}(2n)$ son image. 

On note $\boldsymbol{{\cal P}^{orth}}(2n)$ l'ensemble des couples $(\lambda,\epsilon)$ o\`u $\lambda\in {\cal P}^{orth}(2n)$ et $\epsilon\in (\{\pm 1\}^{Jord_{bp}(\lambda)})/\{\pm 1\}$, le groupe $\{\pm1\}$ s'envoyant diagonalement dans  $\{\pm 1\}^{Jord_{bp}(\lambda)}$.  On note $\boldsymbol{\underline{{\cal P}}^{orth}}(2n)$ l'ensemble des couples $(\underline{\lambda},\epsilon)$ o\`u $\underline{\lambda}\in \underline{{\cal P}}^{orth}(2n)$ et $\epsilon\in (\{\pm 1\}^{Jord_{bp}(\lambda)})/\{\pm 1\}$.
La correspondance de Springer g\'en\'eralis\'ee   \'etablit une bijection entre $\boldsymbol{\underline{{\cal P}}^{orth}}(2n)$ et l'ensemble des couples $(\rho,k)$ tels que

$k\in {\mathbb N}$, $k$ est pair et $k^2\leq 2n$;

si $k>0$, $\rho\in \hat{W}_{n- k^2/2}$; si $k=0$, $\rho\in \hat{W}^D_{n}$.

On note $(\rho_{\underline{\lambda},\epsilon},k_{\lambda,\epsilon})$ le couple associ\'e \`a un \'el\'ement $(\underline{\lambda},\epsilon)\in \boldsymbol{\underline{{\cal P}}^{orth}}(2n)$. L'entier $k_{\lambda,\epsilon}$ ne d\'epend que de l'image $(\lambda,\epsilon)$ de $(\underline{\lambda},\epsilon)$ dans $\boldsymbol{{\cal P}^{orth}}(2n)$. Il se calcule de la fa\c{c}on suivante. Notons $i_{1}>...>i_{m}$ les entiers impairs $i$ tels que $mult_{\lambda}(i)$ soit impair.  Posons 
$$h= \sum_{j=1,...,m}(-1)^j(1-\epsilon(i_{j})).$$
Alors $k_{\lambda,\epsilon}=\vert h\vert $. En particulier, si $\epsilon=1$, c'est-\`a-dire $\epsilon(i)=1$ pour tout $i\in Jord_{bp}(\lambda)$, on a $k_{\lambda,1}=0$ et $\rho_{\lambda,1}\in \hat{W}^D_{n}$.

Une partition $\lambda=(\lambda_{1}\geq \lambda_{2}\geq...)\in {\cal P}^{orth}(2n)$ est sp\'eciale si et seulement si   
$\lambda_{2j-1}$ et $\lambda_{2j}$ sont de m\^eme parit\'e pour tout $j\geq1$.  Cela \'equivaut \`a ce que  $^t\lambda$ soit symplectique. Notons ${\cal P}^{orth,sp}(2n)$ l'ensemble des partitions  orthogonales sp\'eciales de $2n$. Consid\'erons une partition $\lambda\in {\cal P}^{orth,sp}(2n)$ et d\'efinissons $i_{1}>...>i_{m}$ comme ci-dessus.   L'entier $m$ est forc\'ement pair. On appelle intervalle de $\lambda$ un ensemble $\Delta$ de l'une des formes suivantes:

pour un entier $h=1,...,m/2$, $\Delta$ est l'ensemble des $i\geq1$ tels que    $mult_{\lambda}(i)\geq1$ et tels que $i_{2h-1}\geq i \geq i_{2h}$;

$\Delta=\{i\}$ o\`u $i$ est un entier impair tel que     $mult_{\lambda}(i)\geq1$ et tel qu'il n'existe pas d'entier $h=1,...,m/2$ de sorte que $i_{2h-1}\geq i\geq i_{2h}$.

Parce que $\lambda$ est sp\'eciale, on v\'erifie que les intervalles sont form\'es d'entiers impairs. Ils forment une partition de $Jord_{bp}(\lambda)$.  On ordonne les intervalles comme dans le cas symplectique. On note $\Delta_{min}$, resp. $\Delta_{max}$, le plus petit, resp. grand, intervalle. On note $Int(\lambda)$ l'ensemble des intervalles.

 Pour $(\underline{\lambda},\epsilon)\in \boldsymbol{\underline{{\cal P}}^{orth}}(2n)$, le symbole  $symb(\rho_{\underline{\lambda},\epsilon})$ ne d\'epend que de $\lambda$, on le note abusivement $symb(\rho_{\lambda,\epsilon})$. L'application $\lambda\mapsto symb(\rho_{\lambda,1})$ est une bijection de ${\cal P}^{orth,sp}(2n)$ sur l'ensemble des symboles sp\'eciaux de rang $n$ et de d\'efaut $0$. Pour $\lambda\in {\cal P}^{orth,sp}(2n)$, on a d\'efini en \cite{W2} VIII.19 une bijection $fam$ entre la famille du symbole $symb(\rho_{\lambda,1})$ et un  certain sous-ensemble  de $ ({\mathbb Z}/2{\mathbb Z})^{Int(\lambda)}\times ({\mathbb Z}/2{\mathbb Z})^{Int(\lambda)}$. 

Soit $\lambda\in {\cal P}^{orth}(2n)$. Il existe une unique partition sp\'eciale $sp(\lambda)\in {\cal P}^{orth,sp}(2n)$ telle que $symb(\rho_{\lambda,1})$ et $symb(\rho_{sp(\lambda),1})$ soient dans la m\^eme famille. Il est connu que $\lambda\leq sp(\lambda)$ et que $sp(\lambda)$ est la plus petite partition  orthogonale sp\'eciale $\lambda'$ telle que $\lambda\leq \lambda'$. Plus g\'en\'eralement, on a le lemme suivant.

\ass{Lemme}{(i) Soit $(\lambda,\epsilon)\in \boldsymbol{{\cal P}^{orth}}(2n)$, supposons $k_{\lambda,\epsilon}=0$. Notons $sp(\lambda,\epsilon)$ l'unique partition sp\'eciale telle que $symb(\rho_{\lambda,1})$ et $symb(\rho_{sp(\lambda,\epsilon),1})$ soient dans la m\^eme famille. Alors $\lambda\leq sp(\lambda,\epsilon)$.

(ii) Soit $\lambda\in {\cal P}^{orth,sp}(2n)$. Pour  $\epsilon\in \{\pm 1\}^{Jord_{bp}(\lambda)}$, les conditions suivantes sont \'equivalentes:

(a) $k_{\lambda,\epsilon}=0$ et $sp(\lambda,\epsilon)=\lambda$;

(b) $\epsilon$ est constant sur tout $\Delta\in Int(\lambda)$.

(iii) Soit $\lambda\in {\cal P}^{orth,sp}(2n)$. L'application $\epsilon\mapsto fam\circ symb(\rho_{\lambda,\epsilon})$ est une bijection entre l'ensemble des $\epsilon$ d\'ecrits au (ii), modulo le groupe diagonal $\{\pm 1\}$,  et le sous-ensemble des $(\tau,\delta)\in({\mathbb Z}/2{\mathbb Z})^{Int(\lambda)}\times ({\mathbb Z}/2{\mathbb Z})^{Int(\lambda)}$  tels que $\delta=0$ et $\tau(\Delta_{max})=0$.} 

La preuve est similaire \`a celle du lemme pr\'ec\'edent. 

\bigskip

\subsection{Dualit\'e, cas symplectique-orthogonal impair}
Soit $n\in {\mathbb N}$. Notons ${\cal S}_{n,1}^{sp}$ l'ensemble des symboles sp\'eciaux de rang $n$ et de d\'efaut impair (ce d\'efaut est alors $1$).  On dispose de bijections 
$$\begin{array}{ccc}{\cal P}^{symp,sp}(2n)&\to& {\cal S}_{n,1}^{sp}\\ \lambda&\mapsto & symb(\rho_{\lambda,1})\\ \end{array}$$
$$\begin{array}{ccc}{\cal P}^{orth,sp}(2n+1)&\to& {\cal S}_{n,1}^{sp}\\ \lambda&\mapsto & symb(\rho_{\lambda,1})\\ \end{array}$$
et d'une involution $d$ de ${\cal S}_{n,1}^{sp}$. On en d\'eduit des bijections $d:{\cal P}^{symp,sp}(2n)\to {\cal P}^{orth,sp}(2n+1)$ et $d:{\cal P}^{orth,sp}(2n+1)\to {\cal P}^{symp,sp}(2n)$ inverses l'une de l'autre d\'efinies par la formule commune $symb(\rho_{d(\lambda),1})=d\circ symb(\rho_{\lambda,1})$. 

Soit $\lambda\in {\cal P}^{symp,sp}(2n)$. On v\'erifie qu'il y a une unique bijection d\'ecroissante de $Int(\lambda)$ sur $Int(d(\lambda))$. Notons $\Delta_{1}>...>\Delta_{r}$ les intervalles de $\lambda$ et $\Delta'_{1}>...> \Delta'_{r}$ ceux de $d(\lambda)$. On a dit que l'involution $d$ des symboles \'echangeait les familles de $symb(\rho_{\lambda,1})$ et de $symb(\rho_{d(\lambda),1})$. D'autre part, ces familles sont param\'etr\'ees par des sous-ensembles de $({\mathbb Z}/2{\mathbb Z})^{Int(\lambda)}\times ({\mathbb Z}/2{\mathbb Z})^{Int(\lambda)}$, resp. $({\mathbb Z}/2{\mathbb Z})^{Int(d(\lambda))}\times ({\mathbb Z}/2{\mathbb Z})^{Int(d(\lambda))}$.  Soit $(X,Y)$ un symbole dans la famille de $symb(\rho_{\lambda,1})$, notons $(\tau,\delta)=fam(X,Y)$ et $(\tau',\delta')=fam\circ d(X,Y)$. On v\'erifie les \'egalit\'es:
$$\tau'(\Delta_{h})=\tau(\Delta_{r+1-h}),\,\, \delta'(\Delta_{h})=\delta(\Delta_{r-h})$$
pour tout $h=1,...,r$, avec la convention $\delta(\Delta_{0})=0$. En particulier cette application \'echange l'ensemble des $(\tau,\delta)$ tels que $\delta=0$ et celui des $(\tau',\delta')$ tels que $\delta'=0$.

On a d\'efini des applications $sp:{\cal P}^{symp}(2n)\to {\cal P}^{symp,sp}(2n)$ et $sp:{\cal P}^{orth}(2n+1)\to {\cal P}^{orth,sp}(2n+1)$. On  \'etend les bijections  $d$ en des applications encore not\'ees $d:{\cal P}^{symp}(2n)\to {\cal P}^{orth,sp}(2n+1)$ et $d:{\cal P}^{orth}(2n+1)\to {\cal P}^{symp,sp}(2n)$ par la formule commune $d(\lambda)=d\circ sp(\lambda)$. Il est connu que ces applications sont d\'ecroissantes: pour $\lambda,\lambda'\in {\cal P}^{symp}(2n)$, resp. $\lambda,\lambda'\in {\cal P}^{orth}(2n+1)$, $\lambda\leq \lambda'$ entra\^{\i}ne $d(\lambda')\leq d(\lambda)$. 

Soit $\lambda\in {\cal P}^{symp}(2n)$. On v\'erifie que $d(\lambda)$ est la plus grande partition $\mu\in {\cal P}^{orth}(2n+1)$ telle que $\mu\leq {^t(\lambda\cup\{1\})}$, cf. \cite{MR} paragraphe 7.

Soit $\mu\in {\cal P}^{orth}(2n+1)$. Ecrivons $\mu=(\mu_{1}=...=\mu_{s}>\mu_{s+1}\geq...)$. Posons $\mu'=(\mu_{1}=...=\mu_{s-1}\geq \mu_{s}-1\geq \mu_{s+1}\geq...)$. 
  On v\'erifie que $d(\mu)$ est la plus grande partition $\lambda\in {\cal P}^{symp}(2n)$ telle que $\lambda\leq{^t\mu'}$, cf. \cite{MR} paragraphe 7. 

Soit $\lambda\in {\cal P}^{symp,sp}(2n)$, resp. $\lambda\in {\cal P}^{orth,sp}(2n+1)$. On a d\'efini l'ensemble $Int(\lambda)$. Soit $\Delta\in Int(\lambda)$. On note $J(\Delta)$ l'ensemble des indices $j\geq1$ tels que $\lambda_{j}\in \Delta$. Hormis les cas particuliers ci-dessous, on note $j_{min}(\Delta)$, resp. $j_{max}(\Delta)$, le plus petit, resp. grand,  \'el\'ement de $J(\Delta)$. Les cas particuliers sont: $\lambda$ symplectique et $\Delta=\Delta_{min}$, auquel cas on pose $j_{max}(\Delta)=\infty$; $\lambda$ orthogonal et $\Delta=\Delta_{max}$, auquel cas $j_{min}(\Delta))$ n'est pas d\'efini (plus exactement,  on peut le d\'efinir en appliquant la d\'efinition ci-dessus, on obtient $j_{min}(\Delta_{max})=1$, mais cette valeur perturberait nos calculs et on consid\`ere que $j_{min}(\Delta_{max})$ n'est pas d\'efini). Remarquons que, si $\lambda\in {\cal P}^{symp,sp}(2n)$, les $j_{min}(\Delta)$ sont impairs et les $j_{max}(\Delta)$ sont pairs (ou $\infty$); si $\lambda\in {\cal P}^{orth,sp}(2n+1)$,  les $j_{min}(\Delta)$ sont pairs et les $j_{max}(\Delta)$ sont impairs.

On d\'efinit une suite de nombres $\zeta(\lambda)=(\zeta(\lambda)_{1},\zeta(\lambda)_{2},...)$ par 
$$\zeta(\lambda)_{j}=\left\lbrace\begin{array}{cc}1,& si\,\,il\,\,existe\,\,\Delta\in Int(\lambda)\,\,tel \,\,que\,\,j=j_{min}(\Delta),\\-1,& si\,\,il\,\,existe\,\,\Delta\in Int(\lambda)\,\,tel \,\,que\,\,j=j_{max}(\Delta),\\ 0,&dans\,\,les\,\,autres\,\,cas.\\ \end{array}\right.$$

\ass{Lemme}{Soit $\lambda\in {\cal P}^{symp,sp}(2n)$, resp. $\lambda\in {\cal P}^{orth,sp}(2n+1)$. On a l'\'egalit\'e $^td(\lambda)=\lambda+\zeta(\lambda)$.}

Preuve. On suppose $\lambda\in {\cal P}^{symp,sp}(2n)$, la preuve \'etant similaire dans le cas orthogonal.  Montrons d'abord

(1) $^td(\lambda)$ est la plus petite partition orthogonale sp\'eciale $\nu$ de $2n+1$ telle que $\nu\geq \lambda\cup \{1\}$.

Puisque $d(\lambda)$ est orthogonale et sp\'eciale, sa transpos\'ee l'est \'egalement. L'in\'egalit\'e $d(\lambda)\leq {^t(\lambda\cup\{1\})}$ entra\^{\i}ne $^td(\lambda)\geq \lambda\cup\{1\}$. Inversement, soit $\nu$ une partition orthogonale sp\'eciale de $2n+1$ telle que $\nu\geq \lambda\cup\{1\}$. Alors $^t\nu$ est encore orthogonale et v\'erifie $^t\nu\leq {^t(\lambda\cup\{1\})}$. Donc $^t\nu\leq d(\lambda)$ puis $\nu\geq {^td}(\lambda)$. Cela d\'emontre (1). 

On v\'erifie facilement que $\lambda+\zeta(\lambda)$ est une partition, c'est-\`a-dire $\lambda_{j}+\zeta(\lambda)_{j}\geq \lambda_{j+1}+\zeta(\lambda)_{j+1}$ pour tout $j\geq1$. Puisque tout intervalle $\Delta\not=\Delta_{min}$ cr\'ee un terme $j_{min}(\Delta)$ pour lequel $\zeta(\lambda)_{j_{min}(\Delta)}=1$ et un terme $j_{max}(\Delta)$ pour lequel $\zeta(\lambda)_{j_{max}(\Delta)}=-1$ et puisque le dernier intervalle $\Delta_{min}$ cr\'ee seulement un $j_{min}(\Delta_{min})$ la somme totale des $\zeta(\lambda)_{j}$ vaut $1$ et $\lambda+\zeta(\lambda)\in {\cal P}(2n+1)$.  Si $\lambda_{1}$ est impair, $\lambda_{1}$ n'est pas dans un intervalle et $\zeta(\lambda)_{1}=0$ donc $\lambda_{1}+\zeta(\lambda)_{1}$ est impair. Si $\lambda_{1}$ est pair, il appartient \`a un intervalle $\Delta$ (le plus grand intervalle). On a $j_{min}(\Delta)=1$, d'o\`u $\zeta(\lambda)_{1}=1$ et $\lambda_{1}+\zeta(\lambda)_{1}$ est  encore impair. Consid\'erons un entier $h\geq1$ et distinguons les cas:

$\lambda_{2h}$ et $\lambda_{2h+1}$ sont impairs; comme ci-dessus, on a alors $\zeta(\lambda)_{2h}=\zeta(\lambda)_{2h+1}=0$ et les termes $\lambda_{2h}+\zeta(\lambda)_{2h}$ et $\lambda_{2h+1}+\zeta(\lambda)_{2h+1}$ sont impairs;

$\lambda_{2h}$ est impair et $\lambda_{2h+1}$ est pair; dans ce cas $\zeta(\lambda)_{2h}=0$ mais $\lambda_{2h+1}$ appartient \`a un intervalle $\Delta$ tel que $j_{min}(\Delta)=2h+1$, donc $\zeta(\lambda)_{2h+1}=1$;  les termes $\lambda_{2h}+\zeta(\lambda)_{2h}$ et $\lambda_{2h+1}+\zeta(\lambda)_{2h+1}$ sont impairs;

$\lambda_{2h}$ est pair et $\lambda_{2h+1}$ est impair; dans ce cas $\zeta(\lambda)_{2h+1}=0$ mais $\lambda_{2h}$ appartient \`a un intervalle $\Delta$ tel que $j_{max}(\Delta)=2h$, donc $\zeta(\lambda)_{2h}=-1$;  les termes $\lambda_{2h}+\zeta(\lambda)_{2h}$ et $\lambda_{2h+1}+\zeta(\lambda)_{2h+1}$ sont impairs;

$\lambda_{2h}$ et $\lambda_{2h+1}$ sont pairs et distincts; dans ce cas, $\lambda_{2h}$ appartient \`a un intervalle $\Delta$ tel que $j_{max}(\Delta)=2h$ et $\lambda_{2h+1}$ appartient \`a l'intervalle suivant $\Delta'$ tel que $j_{min}(\Delta')=2h+1$; on a $\zeta(\lambda)_{2h}=-1$ et $\zeta(\lambda)_{2h+1}=1$; les termes $\lambda_{2h}+\zeta(\lambda)_{2h}$ et $\lambda_{2h+1}+\zeta(\lambda)_{2h+1}$ sont impairs;
  
$\lambda_{2h}$ et $\lambda_{2h+1}$ sont pairs et \'egaux; dans ce cas, $\lambda_{2h} =\lambda_{2h+1}$ appartient \`a un intervalle $\Delta$ tel que $j_{min}(\Delta)<2h<2h+1<j_{max}(\Delta)$ et $\zeta(\lambda)_{2h}=\zeta(\lambda)_{2h+1}=0$;  les termes $\lambda_{2h}+\zeta(\lambda)_{2h}$ et $\lambda_{2h+1}+\zeta(\lambda)_{2h+1}$ sont pairs et \'egaux.

Cela montre d'abord que les termes pairs de la partition $\lambda+\zeta(\lambda)$ interviennent par paires, donc sont de multiplicit\'e paire, c'est-\`a-dire que $\lambda+\zeta(\lambda)$ est orthogonale. Cela montre ensuite que deux termes $\lambda_{2h}+\zeta(\lambda)_{2h}$ et $\lambda_{2h+1}+\zeta(\lambda)_{2h+1}$ sont de la m\^eme parit\'e. Donc $\lambda+\zeta(\lambda)$ est sp\'eciale.

Pour $k\geq1$, on voit que $S_{k}(\zeta(\lambda))$ vaut $1$ s'il existe $\Delta\in Int(\lambda)$ tel que $j_{min}(\Delta)\leq k<j_{max}(\Delta)$    et vaut $0$ sinon. Ecrivons $\lambda=(\lambda_{1}\geq...\geq \lambda_{l}>0)$. Alors $\lambda\cup \{1\}=(\lambda_{1},...,\lambda_{l},1,0)$. Si $k\leq l$, on a  
 $$S_{k}(\lambda+\zeta(\lambda))=S_{k}(\lambda)+S_{k}(\zeta(\lambda))\geq S_{k}(\lambda)=S_{k}(\lambda\cup \{1\}).$$
Si $k\geq l+1$, on a $k\geq j_{min}(\Delta_{min})$ et $S_{k}(\zeta(\lambda))=1$. Le m\^eme calcul conduit  \`a  l'\'egalit\'e
 $$S_{k}(\lambda+\zeta(\lambda))=S_{k}(\lambda\cup \{1\}).$$
 Donc $\lambda+\zeta(\lambda)\geq \lambda\cup\{1\}$.
 
 Soit maintenant $\nu$ une partition orthogonale sp\'eciale de $2n+1$ telle que $\nu\geq \lambda\cup\{1\}$. Soit $k\geq1$. On a $S_{k}(\nu)\geq S_{k}(\lambda\cup\{1\})\geq S_{k}(\lambda)$.  Supposons  que $S_{k}(\nu)< S_{k}(\lambda+\zeta(\lambda))$. Alors    $S_{k}(\nu)=S_{k}(\lambda)$ et $S_{k}(\zeta(\lambda))=1$. Donc il existe un intervalle $\Delta$ tel que $j_{min}(\Delta)\leq k < j_{max}(\Delta)$. On v\'erifie alors que $S_{k}(\lambda)$ est pair. Supposons de plus $k$ impair. Alors $S_{k}(\nu)$ est impair parce que $\nu$ est sp\'eciale. L'\'egalit\'e $S_{k}(\nu)=S_{k}(\lambda)$ est contradictoire. Cela d\'emontre que, pour $k$ impair, $S_{k}(\nu)\geq S_{k}(\lambda+\zeta(\lambda))$. Supposons maintenant que $k$ est  pair.  Les in\'egalit\'es $j_{min}(\Delta)\leq k< j_{max}(\Delta)$ et le fait que $j_{min}(\Delta)$ est impair tandis que $j_{max}(\Delta)$ est pair ou infini entra\^{\i}nent  que $j_{min}(\Delta)\leq k-1< j_{max}(\Delta)$ et $j_{min}(\Delta)< k+1< j_{max}(\Delta)$.   D'apr\`es ce que l'on vient de d\'emontrer, on a $S_{k-1}(\nu)\geq S_{k-1}(\lambda+\zeta(\lambda))=S_{k-1}(\lambda)+1$. Avec l'\'egalit\'e $S_{k}(\nu)=S_{k}(\lambda)$, cela entra\^{\i}ne $\nu_{k}< \lambda_{k}$. D'autre part, $\lambda_{k}$ et $\lambda_{k+1}$ sont dans un m\^eme intervalle. L'entier $k$ \'etant pair, cela entra\^{\i}ne $\lambda_{k+1}=\lambda_{k}$, donc $\nu_{k+1}\leq \nu_{k}< \lambda_{k}= \lambda_{k+1}$. Avec l'\'egalit\'e $S_{k}(\nu)=S_{k}(\lambda)$, cela entra\^{\i}ne $S_{k+1}(\nu)< S_{k+1}(\lambda)$, ce qui contredit l'hypoth\`ese $\nu\geq \lambda\cup\{1\}$. Cette contradiction d\'emontre encore l'in\'egalit\'e $S_{k}(\nu)\geq S_{k}(\lambda+\zeta(\lambda))$. Celle-ci est donc vraie pour tout $k$, d'o\`u $\nu\geq \lambda+\zeta(\lambda)$.

On a donc prouv\'e que $\lambda+\zeta(\lambda)$  \'etait la plus petite partition orthogonale sp\'eciale $\nu$ de $2n+1$ telle que $\nu\geq \lambda\cup \{1\}$. Le lemme r\'esulte alors de (1). $\square$

 \bigskip

\subsection{Dualit\'e, cas  orthogonal pair}
Soit $n\in {\mathbb N}$. Notons ${\cal S}_{n,0}^{sp}$ l'ensemble des symboles sp\'eciaux de rang $n$ et de d\'efaut pair (ce d\'efaut  est alors $0$).  On dispose d'une bijection  
$$\begin{array}{ccc}{\cal P}^{orth,sp}(2n)&\to& {\cal S}_{n,0}^{sp}\\ \lambda&\mapsto & symb(\rho_{\lambda,1})\\ \end{array}$$
 et d'une involution $d$ de ${\cal S}_{n,0}^{sp}$. On en d\'eduit une  involution $d:{\cal P}^{orth,sp}(2n)\to {\cal P}^{orth,sp}(2n)$  d\'efinie par la formule $symb(\rho_{d(\lambda),1})=d\circ symb(\rho_{\lambda,1})$. 

Soit $\lambda\in {\cal P}^{orth,sp}(2n)$. On v\'erifie qu'il y a une unique bijection d\'ecroissante de $Int(\lambda)$ sur $Int(d(\lambda))$. Notons $\Delta_{1}>...>\Delta_{r}$ les intervalles de $\lambda$ et $\Delta'_{1}>...> \Delta'_{r}$ ceux de $d(\lambda)$. On a dit que l'involution $d$ des symboles \'echangeait les familles de $symb(\rho_{\lambda,1})$ et de $symb(\rho_{d(\lambda),1})$. D'autre part, ces familles sont param\'etr\'ees par des sous-ensembles de $({\mathbb Z}/2{\mathbb Z})^{Int(\lambda)}\times ({\mathbb Z}/2{\mathbb Z})^{Int(\lambda)}$, resp. $({\mathbb Z}/2{\mathbb Z})^{Int(d(\lambda))}\times ({\mathbb Z}/2{\mathbb Z})^{Int(d(\lambda))}$.  Soit $(X,Y)$ un symbole dans la famille de $symb(\rho_{\lambda,1})$, notons $(\tau,\delta)=fam(X,Y)$ et $(\tau',\delta')=fam\circ d(X,Y)$. On v\'erifie    les \'egalit\'es suivantes, pour tout $h=1,...,r$:

$$ \delta'(\Delta_{h})=\delta(\Delta_{r-h}),$$
avec la convention $\delta(\Delta_{0})=0$;

si le d\'efaut de $(X,Y)$ est strictement positif,
$$\tau'(\Delta_{h})=\tau(\Delta_{r+1-h}) ;$$

si ce d\'efaut est nul,
$$\tau'(\Delta_{h})=\tau(\Delta_{r+1-h})-\tau(\Delta_{1}) .$$
 
 En particulier cette application \'echange l'ensemble des $(\tau,\delta)$ tels que $\delta=0$ et celui des $(\tau',\delta')$ tels que $\delta'=0$.

On a d\'efini  l'application $sp:{\cal P}^{orth}(2n)\to {\cal P}^{orth,sp}(2n)$. On  \'etend l'involution   $d$ en une application encore not\'ee $d:{\cal P}^{orth}(2n)\to {\cal P}^{orth,sp}(2n)$ par la formule  $d(\lambda)=d\circ sp(\lambda)$. Il est connu que cette application est d\'ecroissante: pour $\lambda,\lambda'\in {\cal P}^{orth}(2n)$, $\lambda\leq \lambda'$ entra\^{\i}ne $d(\lambda')\leq d(\lambda)$. 

Soit $\lambda\in {\cal P}^{orth}(2n)$. On v\'erifie que $d(\lambda)$ est la plus grande partition $\mu\in {\cal P}^{orth}(2n)$ telle que $\mu\leq {^t\lambda }$, cf. \cite{MR} paragraphe 7.

Soit $\lambda\in {\cal P}^{orth,sp}(2n)$. On a d\'efini l'ensemble $Int(\lambda)$. Soit $\Delta\in Int(\lambda)$.  On note  $J(\Delta)$ l'ensemble des indices $j\geq1$ tels que $\lambda_{j}\in \Delta$. On note $j_{min}(\Delta)$, resp. $j_{max}(\Delta)$, le plus petit, resp. grand, \'el\'ement de $J(\Delta)$.   Le nombre $j_{min}(\Delta)$, resp. $j_{max}(\Delta)$, est impair, resp. pair. On d\'efinit une suite de nombres $\zeta(\lambda)=(\zeta(\lambda)_{1},\zeta(\lambda)_{2},...)$ par 
$$\zeta(\lambda)_{j}=\left\lbrace\begin{array}{cc}1,& si\,\,il\,\,existe\,\,\Delta\in Int(\lambda)\,\,tel \,\,que\,\,j=j_{min}(\Delta),\\-1,& si\,\,il\,\,existe\,\,\Delta\in Int(\lambda)\,\,tel \,\,que\,\,j=j_{max}(\Delta),\\ 0,&dans\,\,les\,\,autres\,\,cas.\\ \end{array}\right.$$

\ass{Lemme}{Soit $\lambda\in {\cal P}^{orth,sp}(2n)$. On a l'\'egalit\'e $^td(\lambda)=\lambda+\zeta(\lambda)$.}

La preuve est similaire \`a celle du lemme pr\'ec\'edent.

\bigskip

\subsection{Dualit\'e et induction}
 Consid\'erons une famille ${\bf n}=(n_{1},...,n_{t},n_{0})$ d'entiers positifs ou nuls. Posons
 $n=\sum_{j=0,...,t}n_{j}$. Posons
 $${\cal P}^{orth}({\bf n})={\cal P}(n_{1})\times...\times{\cal P}(n_{t})\times {\cal P}^{orth}(2n_{0}+1)$$
 et
 $${\cal P}^{symp}({\bf n})={\cal P}(n_{1})\times...\times{\cal P}(n_{t})\times {\cal P}^{symp}(2n_{0}).$$
 On d\'efinit une op\'eration d'induction
 $$\begin{array}{ccc}{\cal P}^{orth}({\bf n})&\to&{\cal P}^{orth}(2n+1),\\ \boldsymbol{\lambda}=(\lambda_{1},...,\lambda_{t},\lambda_{0})&\mapsto&ind(\boldsymbol{\lambda})\\ \end{array}$$
 de la fa\c{c}on suivante: $ind(\boldsymbol{\lambda})$ est la plus grande partition orthogonale $\lambda$ telle que 
 $$\lambda\leq (\lambda_{1}+\lambda_{1})+...+(\lambda_{t}+\lambda_{t})+\lambda_{0}.$$
 L'ensemble ${\cal P}^{orth}({\bf n})$ \'etant le produit d'ensembles ordonn\'es, il l'est aussi par l'ordre produit. On  v\'erifie que l'application d'induction est strictement croissante.
 
 On d\'efinit l'application
 $$cup:{\cal P}^{symp}({\bf n})\to {\cal P}^{symp}(2n)$$
 par la formule 
 $$cup(\lambda_{1},...,\lambda_{t},\lambda_{0})=(\lambda_{1}\cup \lambda_{1})\cup...\cup(\lambda_{t}\cup \lambda_{t})\cup\lambda_{0}.$$
 
 On d\'efinit  enfin une dualit\'e $d:{\cal P}^{symp}({\bf n})\to {\cal P}^{orth}({\bf n})$. C'est le produit des applications $\lambda\mapsto {^t\lambda}$ sur chaque facteur ${\cal P}(n_{i})$ pour $i=1,...,t$ et de la dualit\'e $d:{\cal P}^{symp}(2n_{0})\to {\cal P}^{orth,sp}(2n_{0}+1)\subset {\cal P}^{orth}(2n_{0}+1)$. On a alors
 
 \ass{Lemme}{Pour $\boldsymbol{\lambda}=(\lambda_{1},...,\lambda_{t},\lambda_{0})\in {\cal P}^{symp}({\bf n})$, on a l'\'egalit\'e $d\circ cup(\boldsymbol{\lambda})=ind\circ d(\boldsymbol{\lambda})$.}
 
 Cf. \cite{BV} corollaire A.4.

\bigskip

\subsection{Induction endoscopique}
Soient $n_{1},n_{2}\in {\mathbb N}$, posons $n=n_{1}+n_{2}$. Soient $\lambda_{1}\in {\cal P}^{symp,sp}(2n_{1})$ et $\lambda_{2}\in {\cal P}^{orth,sp}(2n_{2})$.  Rappelons la d\'efinition de   l'induite endoscopique $ind(\lambda_{1},\lambda_{2})\in {\cal P}^{symp}(2n)$, cf. \cite{W2} XI.6. On note $J^+$ l'ensemble des entiers $j\geq1$ tels que

  $\lambda_{1,j}$ est pair, $\lambda_{2,j}$ est impair et il  existe   $\Delta\in Int(\lambda_{1})\cup Int(\lambda_{2})$ de sorte que $j=j_{min}(\Delta)$ (cela entra\^{\i}ne que $j$ est impair).

On note $J^-$  l'ensemble des entiers $j\geq1$ tels que

  $\lambda_{1,j}$ est pair, $\lambda_{2,j}$ est impair et il  existe   $\Delta\in Int(\lambda_{1})\cup Int(\lambda_{2})$ de sorte que $j=j_{max}(\Delta)$ (cela entra\^{\i}ne que $j$ est pair).
  
  On v\'erifie que $J^+$ et $J^-$ ont m\^eme nombre d'\'el\'ements et que, si on note leurs \'el\'ements $j_{1}^+<... <j_{r}^+$ et $j_{1}^-<...<j_{r}^-$, on a
  $$j_{1}^+<j_{1}^-<j_{2}^+<j_{2}^-<... <j_{r}^+<j_{r}^-.$$
  On note $\xi=(\xi_{1},\xi_{2},...)$ la famille d\'efinie par $\xi_{j}=1$ si $j\in J^+$, $\xi_{j}=-1$ si $j\in J^-$ et $\xi_{j}=0$ pour $j\geq1$ tel que $j\not\in J^+\cup J^-$. Alors $ind(\lambda_{1},\lambda_{2})=\lambda_{1}+\lambda_{2}+\xi$.
  
 \ass{Proposition}{On a l'in\'egalit\'e $d(\lambda_{1})\cup d(\lambda_{2})\leq d(ind(\lambda_{1},\lambda_{2}))$. } 
 
 Preuve. Les deux membres de l'in\'egalit\'e \`a prouver sont des partitions de $2n+1$.  Les partitions $d(\lambda_{1})$ et $d(\lambda_{2})$ sont orthogonales, leur r\'eunion l'est aussi. D'apr\`es la caract\'erisation de $d(ind(\lambda_{1},\lambda_{2}))$ donn\'ee en 1.6, il suffit de prouver l'in\'egalit\'e
 
  $$d(\lambda_{1})\cup d(\lambda_{2})\leq {^t(ind(\lambda_{1},\lambda_{2})\cup\{1\}}),$$
  ou encore
  $$^td(\lambda_{1})+{^td(\lambda_{2})}\geq ind(\lambda_{1},\lambda_{2})\cup\{1\},$$
  ou encore, d'apr\`es les lemmes 1.6 et 1.7 et la d\'efinition ci-dessus
  $$\lambda_{1}+\zeta(\lambda_{1})+\lambda_{2}+\zeta(\lambda_{2})\geq (\lambda_{1}+\lambda_{2}+\xi)\cup\{1\}.$$
  Cette in\'egalit\'e se traduit par les in\'egalit\'es suivantes, pour tout $k\geq1$:
  
 (1)  $S_{k}(\zeta(\lambda_{1}))+S_{k}(\zeta(\lambda_{2}))\geq S_{k}(\xi)$, si $k\leq l(ind(\lambda_{1},\lambda_{2}))$;
  
 (2)  $S_{k}(\zeta(\lambda_{1}))+S_{k}(\zeta(\lambda_{2}))\geq S_{k}(\xi)+1$, si $k> l(ind(\lambda_{1},\lambda_{2}))$.
  
  Les entiers $S_{k}(\zeta(\lambda_{1}))$, $S_{k}(\zeta(\lambda_{2}))$ et  $S_{k}(\xi)$ valent toujours $0$ ou $1$. L'in\'egalit\'e (1) est donc v\'erifi\'ee si $S_{k}(\xi)=0$. Supposons $S_{k}(\xi)=1$. Avec les notations introduites plus haut, il existe alors $s\in \{1,...,r\}$ tel que $j_{s}^+\leq k< j_{s}^-$. L'entier $\lambda_{1,j_{s}^-}$ est pair et l'entier $\lambda_{2,j_{s}^-}$ est impair. Il existe donc $\Delta_{1}\in Int(\lambda_{1})$ et $\Delta_{2}\in Int(\lambda_{2})$ tels que $\lambda_{1,j_{s}^-}\in \Delta_{1}$ et $\lambda_{2,j_{s}^-}\in \Delta_{2}$.   Posons $u=max(j_{min}(\Delta_{1}),j_{min}(\Delta_{2}))$. Alors $\lambda_{1,u}\in \Delta_{1}$ est pair et $\lambda_{2,u}\in \Delta_{2}$ est impair. En appliquant la d\'efinition de $J^+$, on voit que $u\in J^+$. On a aussi $u\leq j_{s}^-$, donc $u\leq j_{s}^+$. On a alors $j_{min}(\Delta_{2})\leq u\leq j_{s}^+\leq k< j_{s}^-\leq j_{max}(\Delta_{2})$ et $j_{min}(\Delta_{1})\leq u\leq j_{s}^+\leq k< j_{s}^-\leq j_{max}(\Delta_{1})$. Ces in\'egalit\'es entra\^{\i}nent $\zeta(\lambda_{1})_{k} =\zeta(\lambda_{2})_{k}=1$. On a alors l'\'egalit\'e 
  $$S_{k}(\zeta(\lambda)_{1})+S_{k}(\zeta(\lambda_{2}))=2= S_{k}(\xi)+1,$$
  qui est plus forte que (1). Cela prouve cette in\'egalit\'e (1).
  
  Supposons $k>l(ind(\lambda_{1},\lambda_{2}))$. Si $S_{k}(\xi)=1$, on vient de voir que l'in\'egalit\'e (2) est v\'erifi\'ee (et que c'est une \'egalit\'e). On peut donc supposer $S_{k}(\xi)=0$ et il suffit de montrer que $S_{k}(\zeta(\lambda_{1}))=1$. Puisque $k>l(ind(\lambda_{1},\lambda_{2}))$, on a $\lambda_{1,k}+\lambda_{2,k}+\xi_{k}=0$. Si $\xi_{k}\not=-1$, cela force $\lambda_{1,k}=0$. Si $\xi_{k}=-1$, alors d'une part $\lambda_{1,k}\leq 1$, d'autre part $k\in J^-$. Donc $\lambda_{1,k}$ est pair et on a encore $\lambda_{1,k}=0$. Donc $k\in J(\Delta_{1,min})$, o\`u $\Delta_{1,min}$ est le plus petit \'el\'ement de $Int(\lambda_{1})$.  Cela entra\^{\i}ne $S_{k}(\zeta(\lambda_{1}))=1$, ce qui ach\`eve la d\'emonstration. $\square$

  \bigskip
  
  \subsection{Intervalles relatifs, induction endoscopique r\'eguli\`ere}
  On conserve les donn\'ees $n_{1},n_{2},\lambda_{1}$ et $\lambda_{2}$. On pose $\lambda=ind(\lambda_{1},\lambda_{2})$.
  
  On a d\'efini en \cite{W2} XI.11 un ensemble d'intervalles de $\lambda$. La terminologie est mal choisie car il se peut que $\lambda$ soit sp\'eciale et que cet ensemble ne soit pas celui d\'efini en 1.3 ci-dessus. Nous appellerons ici intervalles relatifs (\`a $\lambda_{1}$ et $\lambda_{2}$) ces nouveaux intervalles. Rappelons leur d\'efinition. On pose
  $${\cal J}=\{j_{min}(\Delta); \Delta\in Int(\lambda_{1})\cup Int(\lambda_{2})\}\cup \{j_{max}(\Delta); \Delta\in Int(\lambda_{1})\cup Int(\lambda_{2})\} ;$$
  
  $${\cal J}^+=\{j_{min}(\Delta);\Delta\in Int(\lambda_{1})\}\cap \{j_{min}(\Delta);\Delta\in Int(\lambda_{2})\};$$

$${\cal J}^-=\{j_{max}(\Delta);\Delta\in Int(\lambda_{1})\}\cap \{j_{max}(\Delta);\Delta\in Int(\lambda_{2})\}.$$

Remarquons que  ${\cal J}$ contient $\infty$ qui est $j_{max}(\Delta)$ pour le plus petit $\Delta\in Int(\lambda_{1})$. 
Appelons intervalle relatif d'indices tout intervalle d'entiers $\{j,...,j'\}$ (avec \'eventuellement $j'=\infty$) v\'erifiant l'une des conditions suivantes:
 
 (1) $j=j'\in {\cal J}^+\cup {\cal J}^-$;
 
 (2) $j<j'$, $j$ et $j'$ sont deux termes cons\'ecutifs de ${\cal J}$ et il existe un unique $d\in \{1,2\}$ et un unique $\Delta_{d}\in Int(\lambda_{d})$ de sorte que $j_{min}(\Delta_{d})\leq j<j'\leq j_{max}(\Delta_{d})$.
 
 Pour tout tel intervalle relatif d'indices $J=\{j,...,j'\}$, on pose $D(J)=\{\lambda_{j''}; j\leq j''\leq j'\}$. On appelle intervalle relatif un tel ensemble $D(J)$. Inversement, pour un intervalle relatif $D$, on note $J(D)$ l'intervalle relatif d'indices $J$ dont il provient et on note $j_{min}(D)$, resp. $j_{max}(D)$, le plus petit, resp.  grand, terme de $J(D)$. 
 On note $Int_{\lambda_{1},\lambda_{2}}(\lambda)$ l'ensemble de ces intervalles relatifs. On montre que cet ensemble d'intervalles relatifs forme une partition de $Jord_{bp}(\lambda)\cup\{0\}$.

 On dit que $\lambda_{1}$ et $\lambda_{2}$ induisent r\'eguli\`erement $\lambda$ si et seulement si tout intervalle relatif est r\'eduit \`a un \'el\'ement. Autrement dit, $Int_{\lambda_{1},\lambda_{2}}(\lambda)$ est la partition maximale de $Jord_{bp}(\lambda)\cup\{0\}$. 
 
 Supposons que $\lambda_{1}$ et $\lambda_{2}$ induisent r\'eguli\`erement $\lambda$. On d\'efinit alors une fonction $\tau_{\lambda_{1},\lambda_{2}}:Jord_{bp}(\lambda)\to {\mathbb Z}/2{\mathbb Z}$ de la fa\c{c}on suivante. Soit $i\in Jord_{bp}(\lambda)$. L'ensemble $\{i\}$ est un intervalle relatif.  Remarquons que $mult_{\lambda}(i)=1$ si et seulement $J(\{i\})$ n'a qu'un \'el\'ement, autrement dit $J(\{i\})$ est du type (1). Si $mult_{\lambda}(i)=1$, on pose $\tau_{\lambda_{1},\lambda_{2}}(i)=0$. Si $mult_{\lambda}(i)\geq2$, $J(\{i\})$ est du type (2) et on note $d(i)\in\{1,2\}$ l'indice tel qu'il existe $\Delta_{d(i)}\in Int(\lambda_{i})$ de sorte que $J(\{i\})\subset \{j_{min}(\Delta_{d(i)}),...,j_{max}(\Delta_{d(i)})\}$. On pose $\tau_{\lambda_{1},\lambda_{2}}(i)=d(i)+1\,\,mod\,\,2{\mathbb Z}$. 
 
  \bigskip
 
 \subsection{Une proposition d'existence}
 
 Soient $n\in {\mathbb N}$ et $\lambda\in {\cal P}^{symp}(2n)$. On se limite ici au cas o\`u tous les termes de $\lambda$ sont pairs. En particulier, $\lambda$ est sp\'eciale. Fixons une fonction $\tau:Jord_{bp}(\lambda)\to {\mathbb Z}/2{\mathbb Z}$ telle que $\tau(i)=0$ pour tout $i\in Jord_{bp}(\lambda)$ tel que $mult_{\lambda}(i)=1$. 
 
 \ass{Proposition}{Soient $\lambda$ et $\tau$ comme ci-dessus. Il existe $n_{1},n_{2}\in {\mathbb N}$ tels que $n_{1}+n_{2}=n$ et il existe $\lambda_{1}\in {\cal P}^{symp,sp}(2n_{1})$ et $\lambda_{2}\in {\cal P}^{orth,sp}(2n_{2})$ tels que
 
 (a) $\lambda_{1}$ et $\lambda_{2}$ induisent r\'eguli\`erement $\lambda$;
 
 (b) $d(\lambda_{1})\cup d(\lambda_{2})=d(\lambda)$;
 
 (c) $\tau_{\lambda_{1},\lambda_{2}}=\tau$.}
 
  Preuve. Notons $\mathfrak{J}^+$ l'ensemble des $j\geq1$ tels que $j$ soit impair et $\lambda_{j}> \lambda_{j+1}$. Notons $\mathfrak{J}^-$ l'ensemble des $j\geq2$ tels que $j$ soit pair et $\lambda_{j-1}>\lambda_{j}$. Les ensembles $\mathfrak{J}^+$ et $\mathfrak{J}^-$ sont disjoints et leur r\'eunion est \'egale \`a  la r\'eunion des couples $\{2k-1,2k\}$, pour $k\geq1$, tels que $\lambda_{2k-1}> \lambda_{2k}$. On note $\mathfrak{x}=(\mathfrak{x}_{1},\mathfrak{x}_{2},...)$ la suite telle que $\mathfrak{x}_{j}=1$ si $j\in \mathfrak{J}^+$, $\mathfrak{x}_{j}=-1$ si $j\in \mathfrak{J}^-$ et $\mathfrak{x}_{j}=0$ si $j\not\in \mathfrak{J}^+\cup\mathfrak{J}^-$. 
 
 On prolonge la fonction $\tau$ \`a $Jord_{bp}(\lambda)\cup \{0\}$ en posant $\tau(0)=0$. Soit $d\in \{1,2\}$. Pour $j\geq1$, disons que $j$ et $j+1$ sont $d$-li\'es si et seulement s'ils v\'erifient l'une des conditions suivantes:
 
 (1) $\lambda_{j}=\lambda_{j+1}$ et $\tau(\lambda_{j})= d+1$ (on veut dire par l\`a $\tau(\lambda_{j})\equiv d+1\,\,mod\,\, 2{\mathbb Z}$);
 
 (2) $j$ est impair et $\lambda_{j}>\lambda_{j+1}$.
 
 Pour deux entiers $1\leq j\leq j'$, disons qu'ils sont $d$-li\'es si et seulement si $k$ et $k+1$ sont $d$-li\'es pour tout $k=j,...,j'-1$. C'est une relation d'\'equivalence. On note $\mathfrak{Int}_{d}$ l'ensemble des classes d'\'equivalence dont le nombre d'\'el\'ements est au moins $2$. Pour $\mathfrak{I}\in \mathfrak{Int}_{d}$, on note $j_{min}(\mathfrak{I})$, resp. $j_{max}(\mathfrak{I})$, le plus petit, resp.  plus grand, \'el\'ement de $\mathfrak{I}$ (\'eventuellement $j_{max}(\mathfrak{I})=\infty$). Montrons que
 
 (3) l'ensemble $\mathfrak{Int}_{d}$ est fini; il contient un \'el\'ement infini   si et seulement si $d=1$;
 
 (4) pour $\mathfrak{I}\in \mathfrak{Int}_{d}$, $j_{min}(\mathfrak{I})$ est impair et $j_{max}(\mathfrak{I})$ est pair ou infini;
 
 (5) pour tout $k\geq1$, il existe au moins un $d\in \{1,2\}$ et un  $\mathfrak{I}\in \mathfrak{Int}_{d}$ tel que $\{2k-1,2k\}\subset \mathfrak{I}_{d}$; les deux \'el\'ements de $\{1,2\}$ v\'erifient cette condition si et seulement si $2k-1\in \mathfrak{J}^+$ (ce qui \'equivaut \`a $2k\in \mathfrak{J}^-$).
 
 (6) pour tout $k\geq1$, il existe au plus un $d\in \{1,2\}$ et un  $\mathfrak{I}\in \mathfrak{Int}_{d}$ tel que $\{2k,2k+1\}\subset \mathfrak{I}_{d}$;
 
 (7) pour tout $j\geq1$, il existe au moins un $d\in \{1,2\}$ et un  $\mathfrak{I}\in \mathfrak{Int}_{d}$ tel que $j\in \mathfrak{I}_{d}$; les deux \'el\'ements de $\{1,2\}$ v\'erifient cette condition si et seulement si $j\in \mathfrak{J}^+ \cup \mathfrak{J}^-$. 
 
 Pour $j\geq l(\lambda)+1$, on a $\lambda_{j}=\lambda_{j+1}=0$  et, puisque $\tau(0)=0$, $j$ et $j+1$ sont $1$-li\'es mais pas $2$-li\'es. Donc $\{l(\lambda)+1,...,\infty\}$ est contenu dans une classe infinie $\mathfrak{I}_{1,min}\in \mathfrak{Int}_{1}$ tandis que, pour $j\geq l(\lambda)+2$, $\{j\}$ forme une classe pour la $2$-\'equivalence donc n'est pas contenu dans un \'el\'ement de $\mathfrak{Int}_{2}$. Cela prouve (3).
 
 Soit $\mathfrak{I}\in \mathfrak{Int}_{d}$, posons $j=j_{min}(\mathfrak{I})$. Montrons que $j$ est impair. Puisque $\mathfrak{I}$ a au moins deux \'el\'ements, $j$ et $j+1$ sont $d$-li\'es. Si la condition (2) est v\'erifi\'ee, $j$ est impair et on a termin\'e. Si (1) est v\'erifi\'ee, on  a $\tau(\lambda_{j})=d+1$. Si $j=1$, $j$ est impair et on a termin\'e. Sinon, puisque $j$ est l'\'el\'ement minimal de $\mathfrak{I}$, $j-1$ et $j$ ne sont pas  $d$-li\'es. Alors le couple $(j-1,j)$ ne v\'erifie pas (2). Donc $j-1$ est pair ou $\lambda_{j-1}=\lambda_{j}$. Dans le premier cas, $j$ est impair et on a termin\'e. Dans le deuxi\`eme cas, on a $\tau(\lambda_{j-1})=\tau(\lambda_{j})=d+1$ mais alors $(j-1,j)$ v\'erifie (1) et $j-1$ et $j$ sont $d$-li\'es, ce qui n'est pas le cas. Cela d\'emontre l'assertion. Un raisonnement similaire prouve que $j_{max}(\mathfrak{I})$ est pair s'il n'est pas infini. D'o\`u (4). 
 
 Soit $k\geq1$. Pour $d=1,2$, dire qu'il existe $\mathfrak{I}\in \mathfrak{Int}_{d}$ tel que $\{2k-1,2k\}\subset \mathfrak{I}$ \'equivaut \`a ce que $2k-1$ et $2k$ soient $d$-li\'es.  Si $\lambda_{2k-1}> \lambda_{2k}$, $2k-1$ et $2k$ v\'erifient (2) et sont $d$-li\'es pour les deux \'el\'ements $d=1,2$.   Mais on a aussi $2k-1\in \mathfrak{J}^+$ et l'assertion (5) est v\'erifi\'ee dans ce cas. Si $\lambda_{2k-1}=\lambda_{2k}$, (2) n'est pas v\'erifi\'ee. Alors $2k-1$ et $2k$ sont $d$-li\'es pour l'unique \'el\'ement $d=\tau(\lambda_{2k-1})+1$. On a aussi $2k-1\not\in \mathfrak{J}^+$ et (5) est encore v\'erifi\'ee.
 
 Soient $k\geq1$ et $d=1,2$.  Le couple $(2k,2k+1)$ ne v\'erifie pas la condition (2). Si $2k,2k+1$ sont $d$-li\'es, la condition (1) est satisfaite. Donc $d=\tau(\lambda_{2k})+1$ est uniquement d\'etermin\'e. D'o\`u (6).

  Soit $j\geq1$.  Posons $k=[(j+1)/2]$. On a $j\in \{2k-1,2k\}$.  Soit $d=1,2$ et $\mathfrak{I}\in \mathfrak{Int}_{d}$. D'apr\`es (4), les conditions $j\in \mathfrak{I}$ et $\{2k-1,2k\}\subset \mathfrak{I}$ sont \'equivalentes. Alors (7) r\'esulte de (5). 
  
  \bigskip
  
  Pour $d=1,2$, d\'efinissons une fonction $p_{d}:{\mathbb N}-\{0\}\to {\mathbb Z}/2{\mathbb Z}$: $p_{d}(j)=1$ s'il existe $\mathfrak{I}\in \mathfrak{Int}_{d}$ tel que $j\in \mathfrak{I}$, $p_{d}(j)=0$ sinon. La relation (4) entra\^{\i}ne
  
  $p_{d}(j)=p_{d}(j+1)$ si $j$ est impair.
  
   La d\'efinition de $\mathfrak{x}$ et l'assertion (7) entra\^{\i}nent l'\'egalit\'e 
  $$(8)\qquad  \mathfrak{x}_{j}\equiv p_{1}(j)+p_{2}(j)+1\,\,mod\,\,2{\mathbb Z}.$$
  
  On va montrer qu'il existe des  suites d'entiers positifs ou nuls $\lambda_{1}$ et $\lambda_{2}$ v\'erifiant les conditions suivantes, pour $j\geq1$:
  
  (9) pour $j\geq1$, $\lambda_{1,j}+\lambda_{2,j}+\mathfrak{x}_{j}=\lambda_{j}$;
  
  (10) pour $j\geq1$ et $d=1,2$, $\lambda_{d,j}\equiv d+p_{d}(j)\,\,mod\,\,2{\mathbb Z}$;
  
  (11) pour $j\geq1$ et $d=1,2$, on a
  
  (a) $\lambda_{d,j}=\lambda_{d,j+1}$ si $j$ est pair, $p_{d}(j)=1$ et il n'existe pas de $\mathfrak{I}\in \mathfrak{Int}_{d}$ tel que $j=j_{max}(\mathfrak{I})$ ou si $j$ est impair et $p_{d}(j)=0$;
  
 (b)  $\lambda_{d,j}> \lambda_{d,j+1}$ si $j$ est pair et il existe $\mathfrak{I}\in \mathfrak{Int}_{d}$ tel que $j=j_{max}(\mathfrak{I})$;
  
  (c) $\lambda_{d,j}\geq \lambda_{d,j+1}$ si $j$ est impair et $p_{d}(j)=1$ ou si $j$ est pair et $p_{d}(j)=0$.

  On raisonne par r\'ecurrence descendante sur $j$. Pour $j\geq l(\lambda)+2$, on pose $\lambda_{1,j}=\lambda_{2,j}=0$. On a vu dans la preuve de (3) que $j$ \'etait contenu dans $\mathfrak{I}_{1,min}$ et n'\'etait  contenu dans aucun \'el\'ement de $\mathfrak{Int}_{2}$. On a aussi $j\not\in \mathfrak{J}^+\cup \mathfrak{J}^-$ donc $\mathfrak{x}_{j}=0$. On voit que toutes nos conditions sont v\'erifi\'ees.
  
  On fixe $j$ et on suppose que l'on a fix\'e des termes $\lambda_{1,j'}$ et $\lambda_{2,j'}$ pour $j'>j$ de sorte que les conditions soient v\'erifi\'ees pour ces $j'$. Pour $d=1,2$, soit $e_{d}\in {\mathbb N}$, posons  $\lambda_{d,j}=\lambda_{d,j+1}+e_{d}$. Traduisons les conditions ci-dessus en termes des entiers $e_{1}$ et $e_{2}$.  La condition (9) \'etant v\'erifi\'ee pour $j+1$, on voit que cette condition pour $j$ \'equivaut \`a
  
  (12) $e_{1}+e_{2}=\lambda_{j}-\lambda_{j+1}+\mathfrak{x}_{j+1}-\mathfrak{x}_{j}$.
  
  La condition (10) \'etant v\'erifi\'ee pour $j+1$, cette condition pour $j$ \'equivaut \`a
  
  (13) $e_{d}\equiv p_{d}(j)+p_{d}(j+1)\,\,mod\,\,2{\mathbb Z}$.
  
  Remarquons que, si (12) est v\'erifi\'ee et si (13) l'est pour un $d\in \{1,2\}$, cette condition (13) est aussi v\'erifi\'ee pour l'autre \'el\'ement de $\{1,2\}$: cela r\'esulte de la parit\'e de $\lambda_{j}$ et de $\lambda_{j+1}$ et de la relation (8). La condition (11) se traduit par les conditions $e_{d}=0$ dans le cas (a), $e_{d}>0$ dans le cas (b) et $e_{d}\geq 0$ dans le cas (c). Remarquons que, dans le cas (a), la condition $e_{d}=0$ est compatible avec (13), autrement dit on a $p_{d}(j)+p_{d}(j+1)\equiv 0\,\,mod\,\,2{\mathbb Z}$. En effet, si $j$ est impair, on a toujours $p_{d}(j)=p_{d}(j+1)$. Si $j$ est pair, la condition de (11)(a) est d'une part que $p_{d}(j)=1$ donc il existe $\mathfrak{I}\in \mathfrak{Int}_{d}$ tel que $j\in \mathfrak{I}$, d'autre part que $j$ n'est pas l'\'el\'ement maximal de $\mathfrak{I}$. Donc $j+1\in \mathfrak{I}$ et $p_{d}(j+1)=p_{d}(j)$. 
  
 Supposons la condition (11)(a) v\'erifi\'ee pour au moins un $d=1,2$, disons pour $d=1$ pour fixer la notation. On n'a pas le choix pour $e_{1}$: on pose $e_{1}=0$. Comme on vient de le dire, la condition (13) est v\'erifi\'ee pour $d=1$. La condition (12) ne laisse plus le choix pour $e_{2}$: on pose $e_{2}=\lambda_{j}-\lambda_{j+1}+\mathfrak{x}_{j+1}-\mathfrak{x}_{j}$. Puisque  (12) est v\'erifi\'ee et aussi (13) pour $d=1$, (13) est aussi v\'erifi\'ee pour $d=2$. Il reste \`a v\'erifier que $e_{2}$ v\'erifie les conditions r\'esultant de (11). Supposons d'abord $j$ pair. L'hypoth\`ese que (11)(a) est v\'erifi\'ee pour $d=1$ signifie, comme on l'a vu ci-dessus, qu'il existe $\mathfrak{I}_{1}\in \mathfrak{Int}_{1}$ tel que $\{j,j+1\}\subset \mathfrak{I}_{1}$. D'apr\`es (6), cette condition ne peut pas \^etre r\'ealis\'ee pour $d=2$. Donc (11)(a) n'est pas v\'erifi\'ee pour $d=2$. Si (11)(c) est v\'erifi\'ee pour $d=2$, on doit seulement voir que $e_{2}\geq0$. Or, puisque $j$ est pair, on a $-\mathfrak{x}_{j}\geq0$ et $\mathfrak{x}_{j+1}\geq0$, donc $\lambda_{j}-\lambda_{j+1}+\mathfrak{x}_{j+1}-\mathfrak{x}_{j}\geq0$ comme on le voulait. Si (11)(b) est v\'erifi\'ee pour $d=2$, on doit montrer que $e_{2}>0$. On a $p_{1}(j)=1$ d'apr\`es (11)(a) pour $d=1$ et $p_{2}(j)=1$ d'apr\`es (11)(b) pour $d=2$. Alors $j\in \mathfrak{J}^-$ d'apr\`es (7) et $-\mathfrak{x}_{j}=1$.  Donc  $\lambda_{j}-\lambda_{j+1}+\mathfrak{x}_{j+1}-\mathfrak{x}_{j}>0$ comme on le voulait. Supposons maintenant $j$ impair. L'hypoth\`ese que (11)(a) est v\'erifi\'ee pour $d=1$ signifie que $p_{1}(j)=0$. D'apr\`es (7), on a $p_{2}(j)=1$ et $j\not\in \mathfrak{J}^+$, donc aussi $j+1\not\in \mathfrak{J}^-$. Ces deux derni\`eres relations entra\^{\i}nent $\mathfrak{x}_{j}=\mathfrak{x}_{j+1}=0$ et $e_{2}=\lambda_{j}-\lambda_{j+1}$. La relation $p_{2}(j)=1$ entra\^{\i}ne que (11)(c) est v\'erifi\'ee pour $d=2$ et que l'on doit seulement prouver que $e_{2}\geq0$, ce qui est clair d'apr\`es la formule pr\'ec\'edente.

Supposons maintenant que (11)(a) n'est v\'erifi\'ee ni  pour $d=1$, ni pour $d=2$. Supposons la condition (11)(b)  v\'erifi\'ee pour au moins un $d=1,2$, disons pour $d=1$. Cela entra\^{\i}ne que $j$ est pair.   Choisissons pour $e_{1}$ le plus petit entier strictement positif v\'erifiant la condition (13). On a $e_{1}=1$ ou $2$. Posons $e_{2}=\lambda_{j}-\lambda_{j+1}+\mathfrak{x}_{j+1}-\mathfrak{x}_{j}-e_{1}$. Comme ci-dessus, on doit montrer que $e_{2}$ v\'erifie les conditions r\'esultant de (11). On a suppos\'e que (11)(a) n'\'etait pas v\'erifi\'ee pour $d=2$. Supposons que (11)(c) soit v\'erifi\'ee pour $d=2$. Il faut voir que $e_{2}\geq0$. D'apr\`es (11)(b) pour $d=1$,   il existe $\mathfrak{I}_{1}\in \mathfrak{Int}_{1}$ tel que $j=j_{max}(\mathfrak{I}_{1})$. Donc $j$ et $j+1$ ne sont pas $1$-li\'es. D'apr\`es (11)(c) pour $d=2$ et parce que $j$ est pair, on a $p_{2}(j)=0$ donc $j$ et $j+1$ ne sont pas $2$ li\'es. Si $\lambda_{j}=\lambda_{j+1}$ la condition (1) est v\'erifi\'ee pour un $d$ donc $j$ et $j+1$ sont $d$-li\'es pour ce $d$. Puisque ce n'est pas le cas, on a $\lambda_{j}\not=\lambda_{j+1}$, donc $\lambda_{j}\geq \lambda_{j+1}+2$, puisque les termes de $\lambda$ sont pairs. Le m\^eme calcul que plus haut conduit \`a l'in\'egalit\'e cherch\'ee $e_{2}\geq0$. Supposons maintenant (11)(b) v\'erifi\'ee pour $d=2$. On doit prouver $e_{2}>0$. On vient de montrer que $j$ et $j+1$ n'\'etaient pas $1$-li\'es. Pour la m\^eme raison, ils ne sont pas $2$-li\'es et cela entra\^{\i}ne encore $\lambda_{j}\geq \lambda_{j+1}+2$. Les conditions (11)(b) pour $d=1,2$ entra\^{\i}nent que $p_{1}(j)=p_{2}(j)=1$, donc $j\in \mathfrak{J}^-$ d'apr\`es (7). Alors $-\mathfrak{x}_{j}=1$ et on voit  que $e_{2}>0$.

Il reste le cas o\`u (11)(c) est v\'erifi\'ee pour $d=1,2$. Puisque $p_{d}(j)=1$ pour au moins un $d$, cette hypoth\`ese entra\^{\i}ne que $j$ est impair et $p_{1}(j)=p_{2}(j)=1$. Donc $j\in \mathfrak{J}^+$, puis $j+1\in \mathfrak{J}^-$. Ces relations entra\^{\i}nent que $\mathfrak{x}_{j}=1$ et $\mathfrak{x}_{j+1}=-1$ et aussi que $\lambda_{j}> \lambda_{j+1}$, donc $\lambda_{j}\geq \lambda_{j+1}+2$.  Puisque $j$ est impair, on a $p_{1}(j+1)=p_{1}(j)=1$. La condition (13) pour $d=1$ signifie  que $e_{1}$ doit \^etre pair. Choisissons $e_{1}=0$, qui v\'erifie la condition r\'esultant de (11)(c) pour $d=1$. Posons 
  $e_{2}=\lambda_{j}-\lambda_{j+1}+\mathfrak{x}_{j+1}-\mathfrak{x}_{j}-e_{1}=\lambda_{j}-\lambda_{j+1}-2$. On a $e_{2}\geq0$, ce qui v\'erifie la condition r\'esultant de (11)(c) pour $d=2$.
 Cela d\'emontre l'existence de nos suites $\lambda_{1}$ et $\lambda_{2}$.

  \bigskip
  
  Fixons donc de telles suites $\lambda_{1}$ et $\lambda_{2}$. La condition (11) entra\^{\i}ne que ce sont des partitions, c'est-\`a-dire qu'elles sont d\'ecroissantes. Montrons que
  
  (14) il existe des entiers $n_{1},n_{2}$ tels que $n_{1}+n_{2}=n$, que $\lambda_{1}$ appartienne \`a ${\cal P}^{symp,sp}(2n_{1})$ et que $\lambda_{2}$ appartienne \`a ${\cal P}^{orth,sp}(2n_{2})$. 
  
  On voit qu'il s'agit de prouver que, pour $d=1,2$ et $k\geq1$, les termes $\lambda_{d,2k-1}$ et $\lambda_{d,2k}$ sont de m\^eme parit\'e et que, quand cette parit\'e est celle de $d$, on a $\lambda_{d,2k-1}=\lambda_{d,2k}$. La premi\`ere propri\'et\'e r\'esulte de (10) et de l'\'egalit\'e $p_{d}(2k-1)=p_{d}(2k)$. Si la parit\'e de $\lambda_{d,2k-1}$ est celle de $d$, cette m\^eme relation (10) entra\^{\i}ne $p_{d}(2k-1)=0$. Mais alors (11)(a) est v\'erifi\'ee pour $2k-1$, d'o\`u $\lambda_{d,2k-1}=\lambda_{d,2k}$. D'o\`u (14).
  
  Gr\^ace \`a cette relation, on peut d\'efinir les ensembles d'intervalles $Int(\lambda_{1})$ et $Int(\lambda_{2})$ et, comme en 1.9, les ensembles $J^+$ et $J^-$ et la fonction $\xi$.
  Montrons que
  
  (15)
   $\{J(\Delta); \Delta\in Int(\lambda_{1})\}=\mathfrak{Int}_{1}$, $\{J(\Delta); \Delta\in Int(\lambda_{2})\}=\mathfrak{Int}_{2}$, $J^+=\mathfrak{J}^+$, $J^-=\mathfrak{J}^-$, $\xi=\mathfrak{x}$.
  
 Soit $d=1,2$.  La r\'eunion des $J(\Delta)$ quand $\Delta$ parcourt $Int(\lambda_{d})$ est l'ensemble des $j\geq1$ tels que $\lambda_{d,j}\equiv d+1\,\,mod\,\,2{\mathbb Z}$.    En vertu de (10), c'est l'ensemble des $j\geq1$ tels que $p_{d}(j)=1$, autrement dit c'est la r\'eunion des \'el\'ements de $\mathfrak{Int}_{d}$. On a donc un m\^eme ensemble d'indices d\'ecoup\'e de deux fa\c{c}ons en intervalles disjoints: les $J(\Delta)$ pour $\Delta\in Int(\lambda_{d})$ ou les $\mathfrak{I}\in \mathfrak{Int}_{d}$. Pour prouver que ces d\'ecoupages sont les m\^emes, il suffit de prouver que les \'el\'ements maximaux de ces intervalles sont les m\^emes, c'est-\`a-dire
 $$\{j_{max}(\Delta); \Delta\in Int(\lambda_{d})\}=\{j_{max}(\mathfrak{I});\mathfrak{I}\in \mathfrak{Int}_{d}\}.$$
 Comme on l'a vu en (3), l'infini intervient dans les deux ensembles si $d=1$ et n'y intervient pas si $d=2$.  Soit $j\geq1$. Par d\'efinition de $Int(\lambda_{d})$,  $j$ appartient \`a l'ensemble de gauche ci-dessus si et seulement si $j$ est pair, $\lambda_{d,j}\equiv d+1\,\,mod\,\,2{\mathbb Z}$ et $\lambda_{d,j}> \lambda_{d,j+1}$. On vient de voir que la congruence est \'equivalente \`a $p_{d}(j)=1$.   Les relations (11) entra\^{\i}nent alors que ces conditions \'equivalent \`a ce que $j$ soit de la forme $j_{max}(\mathfrak{I})$ pour un $\mathfrak{I}\in \mathfrak{Int}_{d}$. Cela d\'emontre les deux premi\`eres \'egalit\'es de (15). Soit $j\in J^+$. Alors $j$ est impair $\lambda_{1,j}$ et $\lambda_{2,j}$ sont "de bonne parit\'e", d'o\`u, comme on l'a vu, $p_{d}(j) =1$ pour $d=1,2$. Alors $j\in \mathfrak{J}^+$ d'apr\`es (7) et l'imparit\'e de $j$. Inversement, soit $j\in \mathfrak{J}^+$. Alors $j$ est impair et, en inversant le raisonnement pr\'ec\'edent, $\lambda_{1,j}$ et $\lambda_{2,j}$ sont de bonne parit\'e. Il existe $ \Delta_{1}\in Int(\lambda_{1})$ et $ \Delta_{2}\in Int(\lambda_{2})$ tels que $j\in  J(\Delta_{1})\cap J(\Delta_{2})$. Si $j=1$,  on a \'evidemment $j=j_{min}(\Delta_{1})=j_{min}(\Delta_{2})$  et $j\in J^+$. Si $j>1$, l'assertion (6) implique qu'il existe $d$ tel que $j-1$ n'appartienne pas \`a $\mathfrak{I}_{d}$, o\`u $\mathfrak{I}_{d}=J(\Delta_{d})$. Alors $j=j_{min}(\mathfrak{I}_{d})$ pour ce $d$, ou encore $j=j_{min}(\Delta_{d})$. Par d\'efinition de l'ensemble $J^+$, on a alors $j\in J^+$. Cela prouve l'\'egalit\'e $J^+=\mathfrak{J}^+$ et l'\'egalit\'e $J^-=\mathfrak{J}^-$ se d\'emontre de m\^eme.    Ces \'egalit\'es et les d\'efinitions de $\xi$ et $\mathfrak{x}$ entra\^{\i}nent la derni\`ere \'egalit\'e de (15). 
 
 L'\'egalit\'e $\xi=\mathfrak{x}$ et la relation (9) entra\^{\i}nent l'\'egalit\'e $Ind(\lambda_{1},\lambda_{2})=\lambda$. Montrons que
 
 (16) $\lambda_{1}$ et $\lambda_{2}$ induisent r\'eguli\`erement $\lambda$. 
 
 Cela signifie que tout intervalle relatif est r\'eduit \`a un seul \'el\'ement. Soit $D$ un tel intervalle relatif. Evidemment, si $J(D)$ est r\'eduit \`a un seul \'el\'ement, $D$ aussi. Supposons que $J(D)$ a au moins deux \'el\'ements. Il v\'erifie la relation (2) de 1.10. Pour fixer la notation, supposons que l'entier $d$ qui figure dans cette relation soit $1$. Il existe donc $\mathfrak{I}_{1}\in \mathfrak{Int}_{1}$ tel que $J(D)\subset \mathfrak{I}_{1}$. Consid\'erons deux \'el\'ements cons\'ecutifs $j,j+1\in J(D)$. Supposons qu'il existe $\mathfrak{I}_{2}\in \mathfrak{Int}_{2}$ tel que $\{j,j+1\}\subset \mathfrak{I}_{2}$. On a $j_{min}(\mathfrak{I}_{2})<j+1\leq j_{max}(D)$. Puisque les termes $j_{min}(D)$ et $j_{max}(D)$ sont par d\'efinition des \'el\'ements cons\'ecutifs de ${\cal J}$, cela entra\^{\i}ne $j_{min}(\mathfrak{I}_{2})\leq j_{min}(D)$. De m\^eme $j_{max}(D)\leq j_{max}(\mathfrak{I}_{2})$. Alors $J(D)\subset \mathfrak{I}_{2}$ ce qui est exclu par 1.10(2). Cela d\'emontre que, pour deux elements $j,j+1\in J(D)$, il n'existe pas de $\mathfrak{I}_{2}\in \mathfrak{Int}_{2}$ tel que $\{j,j+1\}\subset \mathfrak{I}_{2}$. Donc $j$ et $j+1$ sont $1$-li\'es mais pas $2$-li\'es. En se reportant aux relations (1) et (2) qui d\'efinissent la liaison, on voit que, si $j$ est impair, le fait que $j$ et $j+1$ ne sont pas $2$-li\'es   entra\^{\i}ne  que $\lambda_{j}=\lambda_{j+1}$, tandis que, si $j$ est pair, le fait que $j$ et $j+1$ sont $1$-li\'es entra\^{\i}ne la m\^eme \'egalit\'e.  Cette \'egalit\'e pour tout couple $j,j+1\in J(D)$ entra\^{\i}ne que $\lambda_{j}$ est constant pour $j\in J(D)$, ce que l'on voulait d\'emontrer.  
 
 Montrons que 
 
 (17) $\tau_{\lambda_{1},\lambda_{2}}=\tau$.
 
 Soit $i\in Jord_{bp}(\lambda)$. Si $mult_{\lambda}(i)=1$, on a $\tau_{\lambda_{1},\lambda_{2}}(i)=\tau(i)=0$ par d\'efinition. Supposons $mult_{\lambda}(i)\geq2$. Comme ci-dessus, il existe un unique $d=1,2$ et un unique $\mathfrak{I}_{d}\in \mathfrak{Int}_{d}$ tel que $J(\{i\})\subset \mathfrak{I}_{d}$. On a alors $\tau_{\lambda_{1},\lambda_{2}}(i)=d+1$. Consid\'erons un couple $j,j+1\in J(\{i\})$. Ils sont $d$-li\'es et on a $\lambda_{j}=\lambda_{j+1}=i$. L'une des relations (1) ou (2) est v\'erifi\'ee pour $d$ et ce ne peut \^etre que (1). Donc $\tau(i)=d+1$, d'o\`u l'\'egalit\'e cherch\'ee 
$\tau_{\lambda_{1},\lambda_{2}}(i)=\tau(i)$. 

Montrons qu'on a l'\'egalit\'e

(18)  $\zeta(\lambda_{1})+\zeta(\lambda_{2})=\zeta(\lambda)+\xi$.

Soit $j\geq1$. Supposons  $j$ impair. Chacune des quatre fonctions vaut $0$ ou $1$ en $j$.  Supposons d'abord $\zeta(\lambda_{1})_{j}=\zeta(\lambda_{2})_{j}=1$. Alors il existe $\mathfrak{I}_{1}\in \mathfrak{Int}_{1}$ et $\mathfrak{I}_{2}\in \mathfrak{Int}_{2}$ de sorte que $j=j_{min}(\mathfrak{I}_{1})=j_{min}(\mathfrak{I}_{2})$. D'apr\`es (7) et (15), on a $j\in J^+$, d'o\`u $\xi_{j}=1$. Si $j=1$,  $j$ est le plus petit indice tel que $\lambda_{j}$ appartienne au plus grand intervalle de $\lambda$ (il s'agit ici des intervalles au sens des partitions sp\'eciales) donc $\zeta(\lambda)_{j}=1$. Si $j>1$, l'hypoth\`ese sur $j$ implique que $j-1$ et $j$ ne sont ni $1$-li\'es, ni $2$-li\'es. Si $\lambda_{j-1}=\lambda_{j}$, $j-1$ et $j$ sont $d$-li\'es pour le $d$ tel que $\tau(\lambda_{j})=d+1$, cf. (1). C'est impossible donc $\lambda_{j-1}> \lambda_{j}$. Puisque $j$ est impair, c'est la condition pour que $j$ soit de la forme $j=j_{min}(\Delta)$ pour un $\Delta\in Int(\lambda)$. Donc $\zeta(\lambda)_{j}=1$. L'\'egalit\'e (18) est v\'erifi\'ee en $j$. Supposons maintenant $\zeta(\lambda_{1})_{j}=1$ et $\zeta(\lambda_{2})_{j}=0$ (un raisonnement analogue vaut si on \'echange les indices $1$ et $2$). Il existe $\mathfrak{I}_{1}\in \mathfrak{Int}_{1}$ tel que $j=j_{min}(\mathfrak{I}_{1})$ mais il n'y a pas de $\mathfrak{I}_{2}$ v\'erifiant la m\^eme \'egalit\'e. Supposons d'abord $p_{2}(j)=1$. De nouveau, $j\in J^+$ et $\xi_{j}=1$. Puisque $p_{2}(j)=1$, le fait  que $j$ ne soit pas le plus petit \'el\'ement d'un \'el\'ement de $\mathfrak{Int}_{2}$ entra\^{\i}ne que $j\geq2$ et que $j-1$ et $j$ sont $2$-li\'es. Puisque $j-1$ est pair, cette condition implique $\lambda_{j-1}=\lambda_{j}$. Donc $\zeta(\lambda)_{j}=0$ et on obtient l'\'egalit\'e cherch\'ee. Supposons au contraire $p_{2}(j)=0$. Alors $j\not\in J^+$ et $\xi_{j}=0$. Si $j=1$, on a $\zeta(\lambda)_{j}=1$ comme ci-dessus. Sinon, $j-1$ et $j$ ne sont pas $1$-li\'es (car $j=j_{min}(\mathfrak{I}_{1})$) et ne sont pas $2$-li\'es (car $p_{2}(j)=0$). Comme ci-dessus, cela entra\^{\i}ne $\lambda_{j-1}>\lambda_{j}$ et $\zeta(\lambda)_{j}=1$. D'o\`u l'\'egalit\'e cherch\'ee. Supposons enfin $\zeta(\lambda_{1})_{j}=\zeta(\lambda_{2})_{j}=0$. D'apr\`es (7), on peut supposer par exemple $p_{1}(j)=1$. Comme ci-dessus, l'hypoth\`ese $\zeta(\lambda_{1})_{j}=0$ implique alors $j\geq2$ et $j-1$ et $j$ sont $1$-li\'es. D'o\`u $\lambda_{j-1}=\lambda_{j}$ et $\tau(\lambda_{j})=0$. La premi\`ere relation entra\^{\i}ne $\zeta(\lambda)_{j}=0$. La seconde entra\^{\i}ne que $j-1$ et $j$ ne sont pas $2$-li\'es. Si $p_{2}(j)=1$,   $j$ est de la forme $j_{min}(\mathfrak{I}_{2})$ et alors $\zeta(\lambda)_{j}=1$ contrairement \`a l'hypoth\`ese. Donc $p_{2}(j)=0$ et $j\not\in J^+$. Donc $\xi_{j}=0$ et on obtient l'\'egalit\'e cherch\'ee. Des calculs similaires valent dans le cas $j$ pair. Cela prouve (18). 

 Cette \'egalit\'e entra\^{\i}ne
 $$\lambda_{1}+\zeta(\lambda_{1})+\lambda_{2}+\zeta(\lambda_{2})=\lambda_{1}+\lambda_{2}+\xi+\zeta(\lambda)=\lambda+\zeta(\lambda).$$
 En utilisant les lemmes 1.6 et 1.7, cette \'egalit\'e se transforme en
 $$^td(\lambda_{1})+{^td(\lambda_{2})}={^td(\lambda)},$$
 qui \'equivaut \`a
 $$d(\lambda_{1})\cup d(\lambda_{2})=d(\lambda).$$
 Cela ach\`eve la d\'emonstration. $\square$

\bigskip

\subsection{Multiplicit\'es}
Soient $n,n',n'',n_{1},n_{2}\in {\mathbb N}$ tels que $n=n'+n''=n_{1}+n_{2}$. 
Soient   $\rho_{1}\in \hat{W}_{n_{1}}$ et $\rho_{2}\in \hat{W}^D_{n_{2}}$. A $\rho_{1}$ est associ\'e un symbole $(X_{1},Y_{1})$ de rang $n_{1}$ et de d\'efaut $1$. On note $\lambda_{1}$ la partition symplectique sp\'eciale de $2n_{1}$ associ\'ee \`a la famille de $(X_{1},Y_{1})$ et on pose $(\tau_{1},\delta_{1})=fam(X_{1},Y_{1})$. A $\rho_{2}$ est associ\'e un symbole $(X_{2},Y_{2})$ de rang $n_{2}$ et de d\'efaut $0$. On note $\lambda_{2}$ la partition orthogonale sp\'eciale de $2n_{2}$ associ\'ee \`a la famille de $(X_{2},Y_{2})$ et on pose $(\tau_{2},\delta_{2})=fam(X_{2},Y_{2})$. On d\'efinit des repr\'esentations $\rho_{2}^+$ et $\rho_{2}^-$ de $W_{n_{2}}$ de la fa\c{c}on suivante.  Introduisons le couple $(\alpha_{2},\beta_{2})\in {\cal P}_{2}(n_{2})$ qui param\`etre $\rho_{2}$. C'est-\`a-dire que, si $\alpha_{2}\not=\beta_{2}$, $\rho_{2}=\rho^D(\alpha_{2},\beta_{2})$; si $\alpha_{2}=\beta_{2}$, il existe un signe $\eta=\pm$ tel que $\rho_{2}=\rho^D(\alpha_{2},\alpha_{2},\eta)$.  Dans ce dernier cas, on pose $\rho_{2}^+=\rho_{2}^-=\rho(\alpha_{2},\alpha_{2})$. Si $\alpha_{2}\not=\beta_{2}$, on sait que l'on peut permuter $\alpha_{2}$ et $\beta_{2}$. Supposons $\alpha_{2}$ plus grand que $\beta_{2}$ pour l'ordre lexicographique (pour le plus petit indice $j$ tel que $\alpha_{2,j}\not=\beta_{2,j}$, on a $\alpha_{2,j}>\beta_{2,j}$). On pose $\rho_{2}^+=\rho(\alpha_{2},\beta_{2})$ et $\rho_{2}^-=\rho(\beta_{2},\alpha_{2})$. 
  
 Soient $(\lambda',\epsilon')\in \boldsymbol{{\cal P}}^{symp}(2n')$, $(\lambda'',\epsilon'')\in \boldsymbol{{\cal P}}^{symp}(2n'')$. Consid\'erons l'hypoth\`ese
 
 (Hyp) $k_{\lambda',\epsilon'}=k_{\lambda'',\epsilon''}=0$. 
 
 Supposons-la v\'erifi\'ee.   Dans \cite{W3} 1.8 et 1.10, on a d\'efini des espaces ${\cal R}=\oplus_{\gamma\in \Gamma}{\cal R}(\gamma)$, ${\cal R}^{glob}\subset {\cal R}$ et une application lin\'eaire $\rho\iota:{\cal R}\to {\cal R}^{glob}$ (ces objets sont relatifs \`a l'entier $n$).   Posons $\gamma=(0,0,n',n'')$. C'est un \'el\'ement de $\Gamma$ et $\rho_{\lambda',\epsilon'}\otimes \rho_{\lambda'',\epsilon''}$ s'identifie \`a un \'el\'ement de ${\cal R}(\gamma)$. On dispose donc de l'\'el\'ement $\rho\iota(\rho_{\lambda',\epsilon'}\otimes \rho_{\lambda'',\epsilon''})\in {\cal R}$. Remarquons en passant que l'\'el\'ement ${\bf a}$ de \cite{W3} 1.10 vaut $(0,0,0,1)$. Posons $\theta=(0,0,n_{1},n_{2})$. C'est aussi un \'el\'ement de $\Gamma$ et, pour $\zeta=\pm$, $\rho_{1}\otimes \rho_{2}^{\zeta}$ s'identifie \`a un \'el\'ement de ${\cal R}(\theta)$. On peut d\'efinir la multiplicit\'e $m(\rho_{1},\rho_{2}^{\zeta};\rho_{\lambda',\epsilon'},\rho_{\lambda'',\epsilon''})$ de $\rho_{1}\otimes \rho_{2}^{\zeta}$ dans $\rho\iota(\rho_{\lambda',\epsilon'}\otimes \rho_{\lambda'',\epsilon''})$ par la formule usuelle
 $$m(\rho_{1},\rho_{2}^{\zeta};\rho_{\lambda',\epsilon'},\rho_{\lambda'',\epsilon''})=\vert W_{n_{1}}\vert ^{-1}\vert W_{n_{2}}\vert ^{-1}$$
 $$\sum_{w_{1}\in W_{n_{1}},w_{2}\in W_{n_{2}}}\rho_{1}(w_{1})\rho_{2}^{\zeta}(w_{2})\rho\iota(\rho_{\lambda',\epsilon'}\otimes \rho_{\lambda'',\epsilon''})(w_{1}\times w_{2}).$$
 On n'a pas besoin d'introduire des conjugaisons complexes dans cette formule puisqu'on sait que les repr\'esentations irr\'eductibles des groupes de type $W_{n}$ ont des caract\`eres r\'eels. En r\'efl\'echissant \`a la d\'efinition de $\rho\iota(\rho_{\lambda',\epsilon'}\otimes \rho_{\lambda'',\epsilon''})$, on voit que sa restriction \`a ${\cal R}(\theta)$ est une "vraie" repr\'esentation, ce qui entra\^{\i}ne que la multiplicit\'e ci-dessus est un entier naturel. 
 
 On a d\'efini en 1.9 l'induite endoscopique $ind(\lambda_{1},\lambda_{2})\in {\cal P}^{symp}(2n)$.
 
\ass{ Proposition}{On suppose v\'erifi\'ee l'hypoth\`ese (Hyp). Soit $\zeta=\pm$.  Si $m(\rho_{1},\rho_{2}^{\zeta};\rho_{\lambda',\epsilon'},\rho_{\lambda'',\epsilon''})\not=0$, alors $\lambda'\cup \lambda''\leq ind(\lambda_{1},\lambda_{2})$.}

Cela r\'esulte de \cite{W2} proposition XI.28. Le lien entre cette proposition et l'\'enonc\'e ci-dessus n'est pas imm\'ediat mais il est expliqu\'e dans la preuve de \cite{W2} proposition XII.7. 

\bigskip

\subsection{Multiplicit\'es, cas  particulier}
On conserve les donn\'ees du paragraphe pr\'ec\'edent. Posons $\lambda=ind(\lambda_{1},\lambda_{2})$. On suppose de plus

$\lambda$ est \`a termes pair; $\lambda_{1}$ et $\lambda_{2}$ induisent r\'eguli\`erement $\lambda$; $\delta_{1}=\delta_{2}=0$.   

On d\'efinit des fonctions $\delta^+,\delta^-, \tau^+,\tau^-:Jord_{bp}(\lambda)\to {\mathbb Z}/2{\mathbb Z}$ de la fa\c{c}on suivante, o\`u on utilise les notations des paragraphes 1.9 et 1.10.  Soit $i\in Jord_{bp}(\lambda)$. On a $\{i\}\in Int_{\lambda_{1},\lambda_{2}}(\lambda)$ puisque  $\lambda_{1}$ et $\lambda_{2}$ induisent r\'eguli\`erement $\lambda$. On pose  $\delta^+(i)=\delta^-(i)=0$ sauf dans le cas o\`u $j_{max}(\{i\})\in J^+$. Dans ce cas, il existe d'uniques $\Delta_{1}\in Int(\lambda_{1})$ et $\Delta_{2}\in Int(\lambda_{2})$ tels que $j_{max}(\{i\})\in J(\Delta_{1})\cap J(\Delta_{2})$ et on pose
$$\delta^+(i)=\tau_{1}(\Delta_{1})+\tau_{2}(\Delta_{2})+1,\,\,\delta^-(i)=\tau_{1}(\Delta_{1})+\tau_{2}(\Delta_{2}).$$
Si $mult_{\lambda}(i)=1$, il existe comme ci-dessus  d'uniques $\Delta_{1}\in Int(\lambda_{1})$ et $\Delta_{2}\in Int(\lambda_{2})$ tels que $j_{max}(\{i\})\in J(\Delta_{1})\cap J(\Delta_{2})$ et on pose $\tau^+(i)=\tau^-(i)=\tau_{1}(\Delta_{1})$. Supposons $mult_{\lambda}(i)\geq2$. Alors il existe un unique $d=1,2$ et un unique $\Delta_{d}\in Int(\lambda_{d})$ tels que $J(\{i\})\subset J(\Delta_{d})$. Si $d=1$, on pose $\tau^+(i)=\tau^-(i)=\tau_{1}(\Delta_{1})$. Si $d=2$, on pose 
$$\tau^+(i)=\tau_{2}(\Delta_{2}),\,\,\tau^-(i)=\tau_{2}(\Delta_{2})+1.$$

Si $i$ n'est pas l'\'el\'ement maximal de $Jord_{bp}(\lambda)$, on note $i^+$ le plus petit \'el\'ement de $Jord_{bp}(\lambda)$ strictement sup\'erieur \`a $i$. Si $i$ est l'\'el\'ement maximal, on pose par convention $\delta^+(i^+)=\delta^-(i^+)=1$. 

{\bf Remarques.} (1) On v\'erifie sur ces formules que $\tau^++\tau^-=\tau_{\lambda_{1},\lambda_{2}}$, cf. 1.10. 

(2) On a montr\'e en \cite{W2} XI.29 remarque, que, pour tout $i\in Jord_{bp}(\lambda)$, on a l'\'egalit\'e $\delta^+(i)+\delta^-(i)=mult_{\lambda}(\geq i)$. Cela \'equivaut \`a $mult_{\lambda}(i)=\delta^+(i)+\delta^-(i)-\delta^+(i^+)-\delta^-(i^+)$. 

\bigskip

Soit $\zeta=\pm$. On introduit les deux conditions suivantes

$$(A)^{\zeta}\,\left\lbrace\begin{array}{c} (i) \quad \lambda'\cup \lambda''=\lambda;\\

(ii) \text{ pour tout }i\in Jord_{bp}(\lambda), \,\,mult_{\lambda''}( i)\equiv \delta^{-\zeta}(i)-\delta^{-\zeta}(i^+)\,\,mod\,\,2{\mathbb Z};\\

(iii) \text{ pour tout }i\in Jord_{bp}(\lambda'), \,\,\epsilon'(i)=(-1)^{\tau^{\zeta}(i)};\\
\text{ pour tout }i\in Jord_{bp}(\lambda''),\,\,\epsilon''(i)=(-1)^{\tau^{-\zeta}(i)}.\\ \end{array}\right.$$

$$(B)^{\zeta}\,\left\lbrace\begin{array}{c}  (i) \quad \lambda'\cup \lambda''=\lambda;\\

(ii)\text{ l'hypoth\`ese (Hyp) de 1.12 est v\'erifi\'ee et }m(\rho_{1},\rho_{2}^{\zeta};\rho_{\lambda',\epsilon'},\rho_{\lambda'',\epsilon''})\not=0. \\ \end{array}\right.$$

\ass{Proposition}{Pour $\zeta=\pm$, les conditions $(A)^{\zeta}$ et $(B)^{\zeta}$ sont \'equivalentes. Si elles sont v\'erifi\'ees, on a $m(\rho_{1},\rho_{2}^{\zeta};\rho_{\lambda',\epsilon'},\rho_{\lambda'',\epsilon''})=1$.}

Cela r\'esulte de \cite{W2} proposition XI.29.  De nouveau, le lien entre cette proposition et l'\'enonc\'e ci-dessus n'est pas imm\'ediat mais il est expliqu\'e dans la preuve de \cite{W2} proposition XII.7. 

\bigskip

\section{Calcul de caract\`eres}

 \subsection{Caract\`eres de repr\'esentations}
 Dans cette deuxi\`eme section, on reprend les donn\'ees et notations de \cite{W3} et \cite{W4}. Rappelons les principales. Le corps de base $F$ est local non-archim\'edien et de caract\'eristique nulle.  On note $p$ sa caract\'eristique r\'esiduelle et $val_{F}$ la valuation usuelle de $F$. Un entier $n\geq1$ est fix\'e. On suppose 
 
 $p>  6n+4$. 
 
  On consid\`ere deux espaces vectoriels sur $F$ de dimension $2n+1$, not\'es $V_{iso}$ et $V_{an}$, munis de formes quadratiques non d\'eg\'en\'er\'ees $Q_{iso}$ et $Q_{an}$. On note $G_{iso}$ et $G_{an}$ les groupes sp\'eciaux orthogonaux de $(V_{iso},Q_{iso})$ et $(V_{an},Q_{an})$. On suppose $G_{iso}$ d\'eploy\'e et $G_{an}$ non quasi-d\'eploy\'e. Pour un indice $\sharp=iso$ ou $an$, on fixe une mesure de Haar sur $G_{\sharp}(F)$ comme en \cite{W4} 1.1. On note $Irr_{unip,\sharp}$ l'ensemble des classes d'isomorphismes de repr\'esentations admissibles irr\'eductibles de $G_{\sharp}(F)$ qui sont de r\'eduction unipotente, cf. \cite{W3} 1.3. 
 
  Soit $\pi\in Irr_{unip,\sharp}$. A $\pi$ est associ\'e son caract\`ere-distribution, c'est-\`a-dire la forme lin\'eaire $\Theta_{\pi}$ sur $C_{c}^{\infty}(G_{\sharp}(F))$ d\'efinie par $\Theta_{\pi}(f)=trace\,\pi(f)$. Restreignons-nous aux fonctions $f$ dont le support est form\'e d'\'el\'ements compacts de $G_{\sharp}(F)$, c'est-\`a-dire d'\'el\'ements dont les valeurs propres dans une cl\^oture alg\'ebrique de $F$ sont de valuation nulle. La repr\'esentation $\pi$ \'etant de niveau $0$, on a donn\'e dans \cite{W1} une formule pour $\Theta_{\pi}(f)$, que nous allons expliciter. 
 
 Dans \cite{W1} paragraphe 10, on a introduit un ensemble $Fac^*_{max}(G_{\sharp})$. A tout $({\cal F},\nu)\in Fac^*_{max}(G_{\sharp})$ sont associ\'es un sous-groupe compact $K_{{\cal F}}^{\dag}$ de $G_{\sharp}(F)$ et un sous-ensemble $K_{{\cal F}}^{\nu}\subset K_{{\cal F}}^{\dag}$. Le groupe $G_{\sharp}(F)$ agit naturellement sur $Fac^*_{max}(G_{\sharp})$. Il r\'esulte facilement des d\'efinitions que l'ensemble des orbites pour cette action est en bijection avec l'ensemble des triplets $(n',n'',\zeta)$, o\`u $(n',n'')\in D(n)$ (c'est-\`a-dire $n',n''\in {\mathbb N}$ et $n'+n''=n$) et $\zeta=\pm$, soumis aux restrictions suivantes
 
  dans le cas o\`u  $\sharp=iso$, on a $\zeta=+ $ si $n''=0$ et $\zeta=-$ si $n''=1$;
 
  dans le cas o\`u $\sharp=an$, on a $n''\geq1$.
 
 On peut choisir un ensemble de repr\'esentants des orbites dans $Fac^*_{max}(G_{\sharp})$ de sorte que, si un \'el\'ement  $({\cal F},\nu)$ de cet ensemble correspond \`a un triplet $(n',n'',\zeta)$, le groupe $K_{{\cal F}}^{\dag}$ soit \'egal au groupe $K_{n',n''}^{\pm}$ de \cite{W3} 1.2 et l'ensemble $K_{{\cal F}}^{\nu}$ soit \'egal \`a $K_{n',n''}^{\zeta}$. 
 
 Consid\'erons un triplet $(n',n'',\zeta)$ comme ci-dessus.  On dispose de la fonction 
 
 \noindent $proj_{cusp}(Res_{n',n''}^{\zeta}(\pi))\in {\cal R}^{par,glob}$, cf. \cite{W3} 1.5. On peut consid\'erer que c'est une fonction sur $K_{n',n''}^{\zeta}$, invariante par $K_{n',n''}^{u}$. On pose
 $$\Theta_{\pi,cusp}(f)=\sum_{(n',n'',z\eta)}mes(K_{n',n''}^{\pm})^{-1}\int_{G_{\sharp}(F)}\int_{G_{\sharp}(F)}f(g^{-1}hg)proj_{cusp}(Res_{n',n''}^{\zeta}(\pi))(h)\,dh\,dg.$$
 Cette int\'egrale est convergente dans cet ordre. Les $(n',n'',\zeta)$ sont soumis aux restrictions ci-dessus. Mais on peut en fait lever celles-ci parce la fonction $proj_{cusp}(Res_{n',n''}^{\zeta}(\pi))$ est nulle si elles ne sont pas v\'erifi\'ees.
 
 Consid\'erons maintenant une partition ${\bf m}=(m_{1}\geq...\geq m_{t}>0)\in {\cal P}(\leq n)$ (c'est-\`a-dire $S({\bf m}):=m_{1}+...+m_{t}\leq n$), posons $n_{0}=n-S({\bf m})$. On suppose
    $n_{0}\geq 1$ si $\sharp=an$.   On associe \`a  ${\bf m}$ un sous-groupe de Levi  $M\subset G_{\sharp}$. Avec les notations de \cite{W3} 1.1, c'est l'ensemble des \'el\'ements qui, pour tout $j=1,...,t$, stabilisent les deux sous-espaces de $V_{\sharp}$ engendr\'es respectivement par $v_{n_{0}+m_{1}+...+m_{j-1}+1},....,v_{n_{0}+m_{1}+...+m_{j}}$ et par 
   
   \noindent $v_{2n+2-n_{0}-m_{1}-...-m_{j}},...,v_{2n+2-n_{0}-m_{1}-...-m_{j-1}-1}$. On a
 $$M\simeq GL(m_{1})\times...\times GL(m_{t})\times G_{n_{0},\sharp},$$
 o\`u $G_{n_{0},\sharp}$ est l'analogue de $G_{\sharp}$ quand $n$ est remplac\'e par $n_{0}$ (ce  groupe est trivial si $\sharp=iso$ et $n_{0}=0$). Pour tout $j=1,...,t$, fixons un sous-groupe compact maximal $K_{m_{j}}\subset GL(m_{j};F)$ et notons $K_{m_{j}}^{u}$ son radical pro-$p$-unipotent. On note $K_{{\bf m}}$, resp. $K_{{\bf m}}^{u}$, le produit de ces groupes.  On note aussi $A_{M}$ le plus grand tore d\'eploy\'e central dans $M$, c'est-\`a-dire le produit des centres des groupes $GL(m_{j})$. On a d\'efini en \cite{W1} paragraphe 11 un ensemble $Fac^*_{max}(M)_{G_{\sharp}-comp}$. A tout \'el\'ement $({\cal F}_{M},\nu)$ de cet ensemble est associ\'e un sous-groupe $K_{{\cal F}_{M}}^{\dag}$ de $M(F)$. Le groupe $M(F)$ agit naturellement sur $Fac^*_{max}(M)_{G_{\sharp}-comp}$. On voit que l'ensemble des orbites est en bijection avec l'ensemble des triplets $(n',n'',\zeta)$ tels que $(n',n'')\in D(n_{0})$ et $\zeta=\pm$, soumis aux restrictions similaires \`a celles ci-dessus. On peut choisir un ensemble de repr\'esentants des orbites de sorte que, si un \'el\'ement $({\cal F}_{M},\nu)$ de cet ensemble correspond \`a un triplet $(n',n'',\zeta)$, le groupe $K_{{\cal F}_{M}}^{\dag}$ soit \'egal \`a $A_{M}(F)K_{{\bf m}}\times K_{n',n''}^{\pm}$.

  L'analogue pour ce groupe $M$ de l'espace ${\cal R}^{par,glob}$ est l'espace
 $${\cal R}^{par,glob}_{{\bf m}}=C^{GL(m_{1})}\otimes...\otimes C^{GL(m_{t})}\otimes {\cal R}^{par,glob}_{n_{0}},$$
 cf. \cite{W3} 1.5. 
  On introduit  l'application lin\'eaire $Res^M:{\mathbb C}[Irr_{unip,M}]\to {\cal R}^{par,glob}_{{\bf m}}$ analogue \`a $Res$. Soit $P$ un sous-groupe parabolique de $G_{\sharp}$ de composante de Levi $M$. Le semi-simplifi\'e du module de Jacquet $\pi_{P}$  s'identifie \`a un \'el\'ement de ${\mathbb C}[Irr_{unip,M}]$. D'apr\`es \cite{W3} 1.5(1) (qui r\'esulte directement de \cite{MP} proposition 6.7) le terme $Res(\pi_{P})$ ne d\'epend pas du choix de $P$ et on a l'\'egalit\'e
   $$Res^M(\pi_{P})=res_{{\bf m}}\circ Res(\pi).$$
   Notons ce terme $Res_{{\bf m}}(\pi)$ et notons ses diverses composantes $Res_{{\bf m},n',n''}^{\zeta}(\pi)$.  On dispose de  la projection cuspidale  $proj_{cusp}(Res_{{\bf m},n',n''}^{\zeta}(\pi))$. On peut consid\'erer que c'est une fonction sur $K_{{\bf m}}\times K_{n',n''}^{\zeta}$, invariante par $K_{{\bf m}}^{u}\times K_{n',n''}^{u}$. 
  Pour une fonction $\phi\in C_{c}^{\infty}(M(F))$, posons
 $$\Theta^M_{\pi,cusp}(\phi)=\sum_{(n',n'',\zeta)}mes((A_{M}(F)K_{{\bf m}}\times K_{n',n''}^{\pm})/A_{M}(F))^{-1}\int_{M(F)/A_{M}(F)}$$
 $$\int_{M(F)}\phi(m^{-1}ym)proj_{cusp}(Res_{{\bf m},n',n''}^{\zeta}(\pi))(y)\,dy\,dm.$$
 Cette int\'egrale converge dans cet ordre. Fixons un groupe $P$ comme ci-dessus, notons $U$ son radical unipotent. Fixons une mesure de Haar sur $U(F)$. De la mesure sur $M(F)$ (fix\'ee comme en \cite{W4} 1.1)  et de celle sur $U(F)$ se d\'eduit une mesure invariante \`a gauche sur $P(F)$, puis une pseudo-mesure sur $P(F)\backslash G(F)$ (pseudo parce qu'elle s'applique \`a des fonctions qui ne sont pas invariantes \`a gauche par $P(F)$ mais qui se transforment selon le module usuel $\delta_{P}$).   D\'efinissons une fonction $f_{U}$ sur $M(F)$ par
 $$f_{U}(m)=\delta_{P}(m)^{1/2}\int_{U(F)}f(mu)\,du.$$
 En vertu de notre hypoth\`ese sur le support de $f$, on peut aussi bien supprimer le facteur $\delta_{P}(m)^{1/2}$, il vaut $1$ si l'int\'egrale est non nulle. 
 D'autre part, pour $g\in G_{\sharp}(F)$, on d\'efinit la fonction $^gf$ sur $G_{\sharp}(F)$ par $^gf(h)=f(g^{-1}hg)$. 
On pose
 $$\Theta_{\pi,{\bf m}, cusp}(f))=\int_{P(F)\backslash G_{\sharp}(F)}\Theta^M_{\pi,cusp}((^gf)_{U})\,dg.$$
 Ce terme ne d\'epend pas du choix de $P$. Remarquons que le terme $\Theta_{\pi,cusp}(f)$ introduit plus haut est \'egal \`a $\Theta_{\pi,\emptyset,cusp}(f)$, o\`u on a not\'e  $\emptyset$  l'unique partition de $0$.
 
 Rappelons que l'on suppose que le support de $f$ est form\'e d'\'el\'ements compacts de $G_{\sharp}(F)$. Le th\'eor\`eme 12 de \cite{W1} affirme l'\'egalit\'e
 $$\Theta_{\pi}(f)=\sum_{{\bf m}}2^{-l({\bf m})}mult!_{{\bf m}} ^{-1}\Theta_{\pi,{\bf m},cusp}(f),$$
  o\`u on a pos\'e
 $$mult!_{{\bf m}}=\prod_{i\geq1}mult_{{\bf m}}(i)!$$
 et not\'e $l({\bf m})$ le nombres de termes non nuls de ${\bf m}$ (qui est not\'e $t$ plus haut). 
  La somme porte sur les partitions   indiqu\'es plus haut, c'est-\`a-dire ${\bf m}\in {\cal P}(\leq n)$ si $\sharp=iso$ et ${\bf m}\in {\cal P}(\leq n-1)$  si $\sharp=an$.

 {\bf Remarque.} Le th\'eor\`eme 12 de \cite{W1} n'est pas tout-\`a-fait \'enonc\'e comme ci-dessus mais on voit facilement que les deux \'enonc\'es sont \'equivalents. 
 
 \bigskip
 
 \subsection{Un lemme \'el\'ementaire}
Soit $\sharp=iso$ ou $an$. Pour $g\in G_{\sharp}(F)$, on dit que $g$ est topologiquement unipotent si et seulement si $lim_{m\to \infty}g^{p^m}=1$.   Pour $X\in \mathfrak{g}_{\sharp}(F)$, on dit que $X$ est topologiquement nilpotent si et seulement si $lim_{m\to \infty}X^m=0$.   Sous certaines  hypoth\`eses sur $p$ (du type $p> A+B\,val_{F}(p)$), l'exponentielle est d\'efinie  sur l'ensemble des \'el\'ements topologiquement nilpotents de $\mathfrak{g}_{\sharp}(F)$ et est une bijection de cet ensemble sur celui des \'el\'ements topologiquement unipotents de $G_{\sharp}(F)$. Pour simplifier les hypoth\`eses sur $p$, on remplace l'exponentielle par l'application $E$ d\'efinie par $E(X)=\frac{1+X/2}{1-X/2}$. Pour $p>2$, c'est une bijection de l'ensemble des \'el\'ements topologiquement nilpotents de $\mathfrak{g}_{\sharp}(F)$ sur celui des \'el\'ements topologiquement unipotents de $G_{\sharp}(F)$. Rappelons que l'on a suppos\'e $p>6n+4$, a fortiori $p>2$.

Soit $(n',n'')\in D(n)$. On suppose $n''\geq1$ si $\sharp=an$. On a d\'efini en \cite{W3} 1.2 le r\'eseau $L_{n',n''}\subset V_{\sharp}$,  le sous-groupe compact $K_{n',n''}^{+}$ de $G_{\sharp}(F)$ et son radical pro-$p$-unipotent $K_{n',n''}^{u}$. On d\'efinit deux r\'eseaux $\mathfrak{k}_{n',n''}$ et $\mathfrak{k}_{n',n''}^{u}$ de $\mathfrak{g}_{\sharp}(F)$: ce sont les sous-ensembles des \'el\'ements $X\in \mathfrak{g}_{\sharp}(F)$ tels que $X(L_{n',n''})\subset L_{n',n''}$ (ce qui entra\^{\i}ne aussi  $X(L^*_{n',n''})\subset L^*_{n',n''}$), resp. $X(L_{n',n''})\subset \varpi L^*_{n',n''}$ et $X(L^*_{n',n''})\subset L_{n',n''}$. On v\'erifie que, pour $X\in \mathfrak{g}_{\sharp}(F)$ topologiquement nilpotent, on a
$$X\in \mathfrak{k}_{n',n''}\iff exp(X)\in K_{n',n''}^+\,\,  X\in \mathfrak{k}^{u}_{n',n''}\iff exp(X)\in K_{n',n''}^{u}.$$
Posons ${\bf G}={\bf SO}(2n'+1)\times {\bf SO}(2n'')_{\sharp}$, avec les notations de \cite{W3} 1.1. On sait que $K^+_{n',n''}/K_{n',n''}^{u}\simeq {\bf G}({\mathbb F}_{q})$. Notons $\boldsymbol{\mathfrak{g}}$ l'alg\`ebre de Lie de ${\bf G}$. On v\'erifie que $\mathfrak{k}_{n',n''}/\mathfrak{k}_{n',n''}^{u}\simeq \boldsymbol{\mathfrak{g}}({\mathbb F}_{q})$. On note encore $E$ l'application d\'efinie par  $E(X)=\frac{1+X/2}{1-X/2}$ sur  
 l'ensemble des \'el\'ements nilpotents de $\boldsymbol{\mathfrak{g}}({\mathbb F}_{q})$. C'est une bijection de cet ensemble sur celui des \'el\'ements unipotents de ${\bf G}({\mathbb F}_{q})$.

Soit $f\in C_{c}^{\infty}(G_{\sharp}(F))$. Supposons que le support de $f$ est form\'e d'\'el\'ements topologiquement unipotents. On d\'eduit de $f$ une fonction $f_{Lie}\in C_{c}^{\infty}(\mathfrak{g}_{\sharp}(F))$. Son support est form\'e d'\'el\'ements topologiquement nilpotents. Pour un tel \'el\'ement $X$, on a $f_{Lie}(X)=f(E(X))$. On d\'eduit aussi de $f$ une fonction $f_{red}$ sur $K_{n',n''}^+$ telle que, pour tout $g\in K_{n',n''}$,
$$f_{red}(g)=\int_{K_{n',n''}^{u}}f(gh)\,dh.$$
Cette fonction est invariante par $K_{n',n''}^{u}$, on peut consid\'erer que c'est une fonction sur ${\bf G}({\mathbb F}_{q})$. Elle est alors \`a support unipotent. On en d\'eduit une fonction $f_{red,Lie}$ sur $\boldsymbol{\mathfrak{g}}({\mathbb F}_{q})$: celle-ci est \`a support nilpotent et, pour un \'el\'ement nilpotent $X\in \boldsymbol{\mathfrak{g}}({\mathbb F}_{q})$, on a l'\'egalit\'e $f_{red,Lie}(X)=f_{red}(E(X))$. Enfin, on d\'eduit de $f_{Lie}$ une fonction $f_{Lie,red}$ sur $\mathfrak{k}_{n',n''}$: pour $X$ dans cet ensemble, 
$$f_{Lie,red}(X)=\int_{\mathfrak{k}_{n',n''}^{u}}f_{Lie}(X+Y)\,dY.$$
Cette fonction est invariante par translations par $\mathfrak{k}_{n',n''}^{u}$. On peut la consid\'erer comme une fonction sur $\boldsymbol{\mathfrak{g}}({\mathbb F}_{q})$. Elle est alors \`a support nilpotent.

\ass{Lemme}{On a l'\'egalit\'e $f_{red,Lie}=f_{Lie,red}$. }

Preuve. En d\'etaillant les d\'efinitions, on voit qu'il s'agit de d\'emontrer l'assertion suivante:

(1) soit $X\in \mathfrak{k}_{n',n''}$ un \'el\'ement topologiquement nilpotent; alors l'application $Y\mapsto E(X)^{-1}E(X+Y)$ envoie bijectivement $\mathfrak{k}_{n',n''}^{u}$ sur $K_{n',n''}^{u}$ et pr\'eserve les mesures.

La d\'emonstration est  \'el\'ementaire, on la laisse au lecteur. $\square$

 On a effectu\'e les constructions ci-dessus pour le groupe $G_{\sharp}$ afin de ne pas introduire de notations suppl\'ementaires. Mais il est clair que les m\^emes constructions et le m\^eme lemme valent pour les groupes de Levi de $G_{\sharp}$ et nous les utiliserons pour ceux-ci.

\subsection{Calcul du caract\`ere sur les \'el\'ements topologiquement unipotents}
Pour $\sharp=iso$ ou $an$, soit $\pi\in Irr_{unip,\sharp}$. On a d\'efini l'\'el\'ement $Res(\pi)\in {\cal R}^{par,glob}$ et l'isomorphisme $k:{\cal R}^{glob}\to {\cal R}^{par,glob}$ en \cite{W3} 1.5 et 1.9. On note $\kappa_{\pi}$ l'\'el\'ement de ${\cal R}^{glob}$ tel que $Res(\pi)=k(\kappa_{\pi})$. 
Soit $f\in C_{c}^{\infty}(G_{\sharp}(F))$. On suppose que tout \'el\'ement $g$ du support de $f$ est topologiquement unipotent. Un tel \'el\'ement est compact, donc $\Theta_{\pi}(f)$ est donn\'e par la formule de 2.1. Nous allons expliciter cette formule \`a l'aide de l'\'el\'ement $\kappa_{\pi}\in {\cal R}^{glob}$.

Consid\'erons  un entier $n_{0}\in \{1,...,n\}$,  une d\'ecomposition $n_{0}=n'+n''$ et une partition ${\bf m}=(m_{1},...,m_{t}>0)\in {\cal P}(n-n_{0})$. Ces donn\'ees sont soumises aux m\^emes restrictions qu'en 2.1:  si $\sharp=an$, on a $n_{0}\geq1$ et $n''\geq1$.
  On a associ\'e \`a ces donn\'ees un groupe de Levi $M$ de $G_{\sharp}$. Soit $\phi\in C_{c}^{\infty}(M(F))$, supposons que le support de $\phi$ est form\'e d'\'el\'ements topologiquement unipotents.  On va d'abord calculer
 $$I=\int_{M(F)}\phi(y)\sum_{\zeta}proj_{cusp}(Res_{{\bf m},n',n''}^{\zeta}(\pi))(y)\,dy.$$
 La somme porte sur les signes $\zeta=\pm$, soumis aux conditions:
 si $\sharp=iso$, on a $\zeta=+$ si $n''=0$ et $\zeta=-$ si $n''=1$. La deuxi\`eme fonction dans l'int\'egrale est \`a support dans le groupe compact $K_{{\bf m}}\times K_{n',n''}^{\pm}$ et est invariante par $K_{{\bf m}}^{u}\times K_{n',n''}^{u}$. On peut l'identifier \`a une fonction sur le groupe ${\bf M}^{\pm}({\mathbb F}_{q})$, o\`u
 $${\bf M}^{\pm}={\bf GL}(m_{1})\times...\times {\bf GL}(m_{t})\times {\bf SO}(2n'+1)\times {\bf O}(2n'')_{\sharp}.$$
 D\'efinissons une fonction $\phi_{res}$ sur $K_{{\bf m}}\times K_{n',n''}^{\pm}$ par
 $$\phi_{res}(y)=\int_{K_{{\bf m}}^{u}\times K_{n',n''}^{u}}\phi(yh)\,dh.$$
 On peut la consid\'erer elle-aussi comme une fonction sur ${\bf M}^{\pm}({\mathbb F}_{q})$.  On a l'\'egalit\'e
 $$I=\sum_{y\in  {\bf M}^{\pm}({\mathbb F}_{q})}\phi_{res}(y)\sum_{\zeta}proj_{cusp}(Res_{{\bf m},n',n''}^{\zeta}(\pi))(y).$$

 On dispose de l'application 
 $$res_{{\bf m}}:{\cal R}^{glob}\to {\mathbb C}[\hat{\mathfrak{S}}_{m_{1}}]\times...\times {\mathbb C}[\hat{\mathfrak{S}}_{m_{t}}]\otimes {\cal R}_{n_{0}}^{glob}$$
 obtenue en it\'erant la construction de \cite{W3} 1.8. 
Notons $\Gamma_{n',n''}$ l'ensemble des $\gamma=(r',r'',N',N'')\in \Gamma_{n_{0}}$ tels que $r^{_{'}2}+r'+N'=n'$, $r^{_{''}2}+N''=n''$. Pour un tel $\gamma$, notons $res_{{\bf m}}(\kappa_{\pi})_{\gamma}$ la composante dans
 $${\mathbb C}[\hat{\mathfrak{S}}_{m_{1}}]\times...\times {\mathbb C}[\hat{\mathfrak{S}}_{m_{t}}]\otimes {\cal R}_{\gamma}$$
 de   $res_{{\bf m}}(\kappa_{\pi})$.   Excluons d'abord le cas o\`u $\sharp=iso$ et $n''=1$. On  voit que
 $$\sum_{\zeta}proj_{cusp}(Res_{{\bf m},n',n''}^{\zeta}(\pi))=\sum_{\gamma\in \Gamma_{n',n''}}proj_{cusp}\circ k^M(res_{{\bf m}}(\kappa_{\pi})_{\gamma}).$$
 Fixons $\gamma=(r',r'',N',N'')\in \Gamma_{n',n''}$. Dans \cite{MW1} 2.12 et 2.13, on a introduit des fonctions $k(r',w')$ sur ${\bf SO}(2n'+1)({\mathbb F}_{q})$ pour $w'\in W_{N'}$ et $k(r'',w)$ sur ${\bf O}(2n'')_{\sharp}({\mathbb F}_{q})$ pour $w''\in W_{N''}$. Une construction analogue vaut pour les groupes ${\bf GL}(m_{j})$: pour $w\in \mathfrak{S}_{m_{j}}$, on d\'efinit une fonction $k(w)$ sur ${\bf GL}(m_{j};{\mathbb F}_{q})$. Posons
 $$W({\bf m},N',N'')=\mathfrak{S}_{m_{1}}\times...\times\mathfrak{S}_{m_{t}}\times W_{N'}\times W_{N''}.$$ 
 La fonction $res_{{\bf m}}(\kappa_{\pi})_{\gamma}$ est d\'efinie sur ce groupe. Pour $w=(w_{1},...,w_{t},w',w'')\in W({\bf m},N',N'')$, on pose $k(r',r'';w)=k(w_{1})\otimes...\otimes k(w_{t})\otimes k(r',w')\otimes k(r'',w'')$. 
 Il r\'esulte des d\'efinitions que
 $$ k^M(res_{{\bf m}}(\kappa_{\pi})_{\gamma})=\vert W({\bf m},N',N'')\vert ^{-1}\sum_{w\in W({\bf m},N',N'')}res_{\mu}(\kappa_{\pi})_{\gamma}(w)k(r',r'';w).$$
 Dans chacun des groupes $\mathfrak{S}_{m_{j}}$, $W_{N'}$ et $W_{N''}$, on d\'efinit usuellement la notion d'\'el\'ement elliptique. L'application $k$ entrelace la projection $proj_{cusp}$ et la projection sur les \'el\'ements elliptiques. Donc
  $$proj_{cusp}\circ k^M(res_{{\bf m}}(\kappa_{\pi})_{\gamma})=\vert W({\bf m},N',N'')\vert ^{-1}\sum_{w\in W({\bf m},N',N'')_{ell}}res_{{\bf m}}(\kappa_{\pi})_{\gamma}(w)k(r',r'';w),$$
  o\`u $W({\bf m},N',N'')_{ell}$ est le sous-ensemble des \'el\'ements elliptiques de $W({\bf m},N',N'')$. On obtient
  $$I=\sum_{\gamma=(r',r'',N',N'')\in \Gamma_{n',n''}}\vert W({\bf m},N',N'')\vert ^{-1}\sum_{w\in W({\bf m},N',N'')_{ell}}res_{{\bf m}}(\kappa_{\pi})_{\gamma}(w)$$
  $$\sum_{y\in  {\bf M}^{\pm}({\mathbb F}_{q})}\phi_{res}(y)k(r',r'';w)(y).$$
  Les hypoth\`eses sur le support de $\phi$ entra\^{\i}nent que $\phi_{res}$ est \`a support unipotent. D'apr\`es la proposition \cite{MW1} 2.16, $k(r',r'';w)$ est nulle sur les unipotents sauf si $r'=r''=0$. Il ne reste  qu'un seul $\gamma$ qui contribue, \`a savoir l'\'el\'ement $\gamma_{n',n''}=(0,0,n',n'')$. D'o\`u
 $$I=\vert W({\bf m},n',n'')\vert ^{-1}\sum_{w\in W({\bf m},n',n'')_{ell}}res_{{\bf m}}(\kappa_{\pi})_{\gamma_{n',n''}}(w)\sum_{y\in  {\bf M}^{\pm}({\mathbb F}_{q})}\phi_{res}(y)k(w)(y),$$
 o\`u on a pos\'e simplement $k(w)=k(0,0;w)$.  
   A ce point, on peut supprimer l'hypoth\`ese restrictive faite plus haut. Si $\sharp=iso$ et $n''=1$, on a $I=0$ car on se limite \`a $\zeta=-$ et $K_{n'',iso}^{-}$ ne contient pas d'\'el\'ement topologiquement unipotent. Mais la formule ci-dessus donne le m\^eme r\'esultat, car pour l'unique \'el\'ement elliptique $w''\in W_{1,ell}$, on a $k(0,w'')_{iso}=0$, cf. \cite{MW1} 2.13.  Notons ${\bf M}$ la composante neutre de ${\bf M}^{\pm}$ et $\boldsymbol{\mathfrak{m}}$ son alg\`ebre de Lie.  On dispose
   de l'application $E$ de 2.2, qui est une bijection de l'ensemble $\boldsymbol{\mathfrak{m}}_{nil}({\mathbb F}_{q})$ des \'el\'ements nilpotents de $\boldsymbol{\mathfrak{m}}({\mathbb F}_{q})$ sur l'ensemble ${\bf M}_{unip}({\mathbb F}_{q})$ des \'el\'ements nilpotents de ${\bf M}({\mathbb F}_{q})$. Notons $\phi_{red,Lie}$ la fonction sur $\boldsymbol{\mathfrak{m}}({\mathbb F}_{q})$ qui est nulle hors des \'el\'ements nilpotents et qui v\'erifie
   $\phi_{red,Lie}(X)=\phi_{red}(E(X))$ pour tout $X$ nilpotent. Pour $w\in W({\bf m},n',n'')_{ell}$, d\'efinissons de m\^eme une fonction $k(w)_{Lie}$. On obtient
    $$(2) \qquad I=\vert W({\bf m},n',n'')\vert ^{-1}\sum_{w\in W({\bf m},n',n'')_{ell}}res_{{\bf m}}(\kappa_{\pi})_{\gamma_{n',n''}}(w)I_{w},$$
    o\`u
    $$I_{w}=\sum_{X\in \boldsymbol{\mathfrak{m}}({\mathbb F}_{q})}\phi_{red,Lie}(X)k(w)_{Lie}(X).$$

   Fixons $w=(w_{1},...,w_{t},w',w'')\in W({\bf m},n',n'')_{ell}$. On peut supposer $sgn_{CD}(w'')=1$ si $\sharp=iso$, $sgn_{CD}(w'')=-1$ si $\sharp=an$, sinon la fonction  $k(0,w'')$ est nulle sur ${\bf SO}_{\sharp}({\mathbb F}_{q})$, cf. \cite{MW1} 2.13. A tout $w_{j}$ est associ\'ee une classe de conjugaison de sous-tore maximal elliptique dans ${\bf GL}(m_{j})$ (qui est d'ailleurs l'unique telle classe). A $w'$, resp. $w''$, est associ\'ee une classe de conjugaison de sous-tore maximal elliptique dans ${\bf SO}(2n'+1)$, resp. ${\bf SO}(2n'')_{\sharp}$. On fixe des tores dans ces classes de conjugaison et on note ${\bf T}_{w}$ leur produit qui est donc un sous-tore maximal elliptique dans ${\bf M}$. On dispose de l'induction de Deligne-Lusztig de ${\bf T}_{w}$ \`a ${\bf M}$. Ce foncteur vaut aussi pour les alg\`ebres de Lie. Notons $\boldsymbol{\mathfrak{t}}_{w}$ l'alg\`ebre de Lie de ${\bf T}_{w}$ et consid\'erons la fonction caract\'eristique de $\{0\}$ dans $\boldsymbol{\mathfrak{t}}_{w}({\mathbb F}_{q})$. On note $Q_{w}$ son image par induction de Deligne-Lusztig, qui est une fonction sur $\boldsymbol{\mathfrak{m}}({\mathbb F}_{q})$, \`a support nilpotent.  On a l'\'egalit\'e
   
   (3) $k(w)_{Lie}=2^{\beta}Q_{w}$, o\`u $\beta=0$ si $n''=0$, $\beta=1$ si $n''>0$. 
   
   En effet, d'apr\`es nos d\'efinitions de \cite{MW1} 2.12 et 2.13, $k(w)$ est \'egal \`a $(-1)^n2^{\beta}$ fois la trace d'un Frobenius sur un faisceau-caract\`ere. D'apr\`es \cite{L1} th\'eor\`eme 1.14, cette trace est \'egale, sur les unipotents, $(-1)^n$ fois l'image par induction de Deligne-Lusztig de la fonction caract\'eristique de $1$ dans ${\bf T}_{w}({\mathbb F}_{q})$. En descendant par l'application $E$ \`a l'alg\`ebre de Lie, on obtient (3). 
   
 En \cite{W4} 1.1, on a fix\'e un  caract\`ere $\psi_{F}$ de $F$ de conducteur $\varpi\mathfrak{o}$.   Il lui   est associ\'e un caract\`ere  de ${\mathbb F}_{q}$ gr\^ace auquel on d\'efinit comme  en \cite{W4} 1.1  une transformation  de Fourier $\varphi\mapsto \hat{\varphi}$  dans $C_{c}^{\infty}(\boldsymbol{\mathfrak{m}}({\mathbb F}_{q}))$. On la normalise de sorte que $\hat{\hat{\varphi}}(X)=\varphi(-X)$. D'apr\`es \cite{L2} proposition 7.2 et \'egalit\'e 6.15(a), on a l'\'egalit\'e
   
   (4) $\hat{Q}_{w}(X)=sgn(w)q^{-n/2}Q_{w}(X)$ pour tout \'el\'ement nilpotent $X\in \boldsymbol{\mathfrak{m}}_{nil}({\mathbb F}_{q})$.
   
   Fixons un point ${\bf X}_{w}\in \boldsymbol{\mathfrak{t}}_{w}({\mathbb F}_{q})$ en position g\'en\'erale. Notons $\varphi[{\bf X}_{w}]$ la fonction caract\'eristique de la classe de conjugaison par ${\bf M}({\mathbb F}_{q})$ de ${\bf X}_{w}$. D'apr\`es \cite{W2} proposition II.8, on a l'\'egalit\'e
   
   (5) $\hat{\varphi}[{\bf X}_{w}](X)=\hat{Q}_{w}(X)$ pour tout \'el\'ement nilpotent $X\in \boldsymbol{\mathfrak{m}}_{nil}({\mathbb F}_{q})$.
   
   En rassemblant (3), (4) et (5), on obtient l'\'egalit\'e $k(w)_{Lie}(X)=sgn(w)q^{n/2}2^{\beta}\hat{\varphi}[{\bf X}_{w}](X)$, pour $X\in \boldsymbol{\mathfrak{m}}_{nil}({\mathbb F}_{q})$.  D'o\`u
  $$ I_{w}=sgn(w)q^{n/2}2^{\beta}\sum_{X\in \boldsymbol{\mathfrak{m}}({\mathbb F}_{q})}\phi_{red,Lie}(X)\hat{\varphi}[{\bf X}_{w}](X),$$
  puis, par la formule de Parseval,
  $$  I_{w}=sgn(w)q^{n/2}2^{\beta}\sum_{X\in \boldsymbol{\mathfrak{m}}({\mathbb F}_{q})}\hat{\phi}_{red,Lie}(X)\varphi[{\bf X}_{w}](X).$$
  Ou encore, en explicitant la fonction $\varphi[{\bf X}_{w}]$, 
   $$  I_{w}=sgn(w)q^{n/2}2^{\beta}\vert {\bf T}_{w}({\mathbb F}_{q})\vert ^{-1}\sum_{x\in {\bf M}({\mathbb F}_{q})}\hat{\phi}_{red,Lie}(x^{-1}{\bf X}_{w}x).$$
 La conjugaison se fait ici par le groupe ${\bf M}({\mathbb F}_{q})$ et on rappelle que ${\bf M}$ est la composante neutre de ${\bf M}^{\pm}$. Mais, dans la formule ci-dessus, on peut remplacer ${\bf X}_{w}$ par un conjugu\'e quelconque par un \'el\'ement de ${\bf M}^{\pm}({\mathbb F}_{q})$. Un tel conjugu\'e v\'erifie les m\^emes propri\'et\'es que ${\bf X}_{w}$. On peut donc remplacer la conjugaison par ${\bf M}({\mathbb F}_{q})$ par la conjugaison par ${\bf M}^{\pm}({\mathbb F}_{q})$ tout entier, \`a condition de diviser par $[  {\bf M}^{\pm}({\mathbb F}_{q}):{\bf M}({\mathbb F}_{q})]$, qui vaut pr\'ecis\'ement $2^{\beta}$. D'o\`u
  $$(6) \qquad    I_{w}=sgn(w)q^{n/2} \vert {\bf T}_{w}({\mathbb F}_{q})\vert ^{-1}\sum_{x\in {\bf M}^{\pm}({\mathbb F}_{q})}\hat{\phi}_{red,Lie}(x^{-1}{\bf X}_{w}x).$$  
   
  Notons $\mathfrak{m}$ l'alg\`ebre de Lie de $M$.  Comme  en 2.2, on d\'efinit une fonction $\phi_{Lie}$ sur $\mathfrak{m}(F)$: elle est \`a support topologiquement nilpotent; pour $X\in \mathfrak{m}(F)$ topologiquement nilpotent, on a $\phi_{Lie}(X)=\phi(E(X))$. On en d\'eduit une fonction
    $\phi_{Lie,red}$ sur $\mathfrak{k}_{{\bf m}}\oplus \mathfrak{k}_{n',n''}$ (avec une d\'efinition \'evidente de $\mathfrak{k}_{{\bf m}}$ et, ci-dessous, de $\mathfrak{k}_{{\bf m}}^{u}$) par 
    $$\phi_{Lie,red}(X)=\int_{\mathfrak{k}^{u}_{{\bf m}}\oplus \mathfrak{k}_{n',n''}^{u}}\phi_{Lie}(X+Y)\,dY$$
    pour tout $X\in \mathfrak{k}_{{\bf m}}\oplus \mathfrak{k}_{n',n''}$. On peut consid\'erer que c'est une fonction sur 
    $\boldsymbol{\mathfrak{m}}({\mathbb F}_{q})$. Le lemme 2.2 dit que $\phi_{red,Lie}=\phi_{Lie,red}$.  On dispose de la fonction $\hat{\phi}_{Lie}$ (la transform\'ee de Fourier de $\phi_{Lie}$) dont on d\'eduit comme ci-dessus une fonction $(\hat{\phi}_{Lie})_{red}$ sur $\mathfrak{k}_{{\bf m}}\oplus \mathfrak{k}_{n',n''}$, que l'on peut consid\'erer comme une fonction sur $\boldsymbol{\mathfrak{m}}({\mathbb F}_{q})$. On v\'erifie l'\'egalit\'e
    $$(\hat{\phi}_{Lie})_{red}=\hat{\phi}_{Lie,red}.$$
    Dans la formule (6), rempla\c{c}ons $\hat{\phi}_{red,Lie}$ par $(\hat{\phi}_{Lie})_{red}$. Les termes de la formule vivent dans $\boldsymbol{\mathfrak{m}}({\mathbb F}_{q})$ mais on peut les relever dans $ \mathfrak{k}_{{\bf m}}\oplus \mathfrak{k}_{n',n''}$. On rel\`eve ainsi ${\bf X}_{w}$ en un \'el\'ement de ce r\'eseau que l'on note  $X_{w}$. La somme en $x\in {\bf M}^{\pm}({\mathbb F}_{q})$ devient une int\'egrale sur $K_{{\bf m}}\times K_{n',n''}^{\pm}$, divis\'ee par la mesure de $K_{{\bf m}}^{u}\times K_{n',n''}^{u}$. On obtient
    $$(7) \qquad I_{w}=sgn(w)q^{n/2}\vert {\bf T}_{w}({\mathbb F}_{q})\vert ^{-1}mes(K_{{\bf m}}^{u}\times K_{n',n''}^{u})^{-1}\int_{K_{{\bf m}}\times K_{n',n''}^{\pm}}(\hat{\phi}_{Lie})_{red}(x^{-1}X_{w}x)\,dx.$$
    Notons $T_{w}$ le centralisateur de $X_{w}$ et $\mathfrak{t}_{w}$ son alg\`ebre de Lie. Le tore $T_{w}$ est non ramifi\'e sur $F$ et poss\`ede une structure naturelle sur $\mathfrak{o}$. On a $\mathfrak{t}_{w}(\mathfrak{o})=\mathfrak{t}_{w}(F)\cap (\mathfrak{k}_{{\bf m}}\oplus \mathfrak{k}_{n',n''})$ et $\varpi\mathfrak{t}_{w}(\mathfrak{o})=\mathfrak{t}_{w}(F)\cap (\mathfrak{k}^{u}_{{\bf m}}\oplus \mathfrak{k}_{n',n''}^{u})$.        Posons ${\cal X}_{w}=X_{w}+\varpi\mathfrak{t}_{w}(\mathfrak{o})$.   Montrons que
    
    (8) pour tout $x\in K_{{\bf m}}\times K_{n',n''}^{\pm}$, on a l'\'egalit\'e
    $$(\hat{\phi}_{Lie})_{red}(x^{-1}X_{w}x)=mes({\cal X}_{w})^{-1} \int_{ K_{{\bf m}}^{u}\times K_{n',n''}^{u}}\int_{ {\cal X}_{w}}\hat{\phi}_{Lie}(x^{-1}y^{-1}Y_{w}yx)\,dY_{w}\,dy.$$
    
On se ram\`ene imm\'ediatement au cas $x=1$ en  conjuguant par $x$ la fonction $\hat{\phi}_{Lie}$. Supposons donc     $x=1$.   
 Posons $T_{w}(F)^{u}=T_{w}(F)\cap   (K_{{\bf m}}^{u}\times K_{n',n''}^{u})$. C'est l'image par $E$ de $\varpi\mathfrak{t}_{w}(\mathfrak{o})$, on a donc $mes(T_{w}(F)^{u})=mes({\cal X}_{w})$. Les \'el\'ements $Y_{w}$ appartiennent \`a $\mathfrak{t}_{w}(F)$ donc $T_{w}(F)$ commute \`a ces \'el\'ements. On peut remplacer l'int\'egrale en $y\in K_{{\bf m}}^{u}\times K_{n',n''}^{u}$  du membre de droite ci-dessus par une int\'egrale en $y\in T_{w}(F)^{u}\backslash (K_{{\bf m}}^{u}\times K_{n',n''}^{u})$, multipli\'ee par $mes({\cal X}_{w})$. Ce facteur fait dispara\^{\i}tre son inverse qui figure dans ce  membre de droite. Consid\'erons l'application
 $$\begin{array}{cccc}\iota:&T_{w}(F)^{u}\backslash (K_{{\bf m}}^{u}\times K_{n',n''}^{u})\times \varpi\mathfrak{t}_{w}(\mathfrak{o})&\to&\mathfrak{m}(F)\\ &(y,Z)&\mapsto& y^{-1}(X_{w}+Z)y-X_{w}\\ \end{array}$$
 Il est clair que son image est contenue dans $\mathfrak{k}_{{\bf m}}^{u}\oplus \mathfrak{k}_{n',n''}^{u}$. Montrons qu'elle est injective. Si $(y,Z)$ et $(y',Z')$ ont m\^eme image, on a $y^{-1}(X_{w}+Z)y=y^{_{'}-1}(X_{w}+Z')y'$. Le point ${\bf X}_{w}$ est en position g\'en\'erale et ses valeurs propres (dans une cl\^oture alg\'ebrique $\bar{{\mathbb F}}_{q}$ de ${\mathbb F}_{q}$) sont distinctes. Les valeurs propres de $X_{w}+Z$ et $X_{w}+Z'$ sont enti\`eres (dans une cl\^oture alg\'ebrique de $F$) et leurs r\'eductions dans $\bar{{\mathbb F}}_{q}$ sont les m\^emes que celles de ${\bf X}_{w}$. On en d\'eduit ais\'ement que les points $X_{w}+Z$ et $X_{w}+Z'$ ne peuvent \^etre conjugu\'es que s'ils sont \'egaux. Donc $Z=Z'$. Alors $y'y^{-1}$ commute \`a $X_{w}+Z'$ et appartient donc \`a $T_{w}(F)$. Cela prouve l'injectivit\'e de $\iota$. L'application $\iota$ est diff\'erentiable. Sa d\'eriv\'ee en un point $(y,Z)$ est l'application
$$\begin{array}{ccc}\mathfrak{t}_{w}(F)\backslash\mathfrak{m}(F)\times \mathfrak{t}_{w}(F)&\to&\mathfrak{m}(F)\\ (\mathfrak{Y},\mathfrak{Z})&\mapsto&y^{-1}([X_{w}+Z,\mathfrak{Y}]+\mathfrak{Z})y\\ \end{array}$$
Celle-ci est bijective    et, parce que les valeurs propres de $X_{w}+Z$ sont enti\`eres et de r\'eductions toutes distinctes, on v\'erifie qu'elle pr\'eserve les mesures. Donc $\iota$ est un isomorphisme local, de jacobien constant de valeur $1$. On en d\'eduit que l'image de $\iota$ est ouverte dans  $\mathfrak{k}_{{\bf m}}^{u}\oplus \mathfrak{k}_{n',n''}^{u}$ et que cette image a m\^eme mesure que l'espace de d\'epart. D'autre part, l'image de $\iota$ est clairement compacte et l'espace de d\'epart a m\^eme mesure que $\mathfrak{k}_{{\bf m}}^{u}\oplus \mathfrak{k}_{n',n''}^{u}$. Cela entra\^{\i}ne que $\iota$ est un isomorphisme pr\'eservant les mesures de $T_{w}(F)^{u}\backslash (K_{{\bf m}}^{u}\times K_{n',n''}^{u})\times \varpi\mathfrak{t}_{w}(\mathfrak{o})$ sur $\mathfrak{k}_{{\bf m}}^{u}\oplus \mathfrak{k}_{n',n''}^{u}$. Le membre de droite de (8) (en $x=1$) s'\'ecrit
$$ \int_{T_{w}(F)^{u}\backslash( K_{{\bf m}}^{u}\times K_{n',n''}^{u})}\int_{  \varpi\mathfrak{t}_{w}(\mathfrak{o})}\hat{\phi}_{Lie}(\iota(y,Z)+X_{w} )\,dZ\,dy.$$
D'apr\`es les propri\'et\'es de $\iota$, c'est aussi
$$\int_{\mathfrak{k}_{{\bf m}}^{u}\oplus \mathfrak{k}_{n',n''}^{u}}\hat{\phi}_{Lie}(X_{w}+Y)\,dY.$$
Mais ceci est la d\'efinition de $(\hat{\phi}_{Lie})_{red}(X_{w})$. Cela d\'emontre (8).    
   
   Utilisons (8) pour transformer (7). L'int\'egrale en $K_{{\bf m}}^{u}\times K_{n',n''}^{u}$ est absorb\'ee par celle en $K_{{\bf m}}\times K_{n',n''}$ mais introduit un facteur $mes(K_{{\bf m}}^{u}\times K_{n',n''}^{u})$ qui compense l'inverse de cette mesure intervenant dans (7). On obtient
   $$I_{w}=sgn(w)q^{n/2}\vert {\bf T}_{w}({\mathbb F}_{q})\vert ^{-1}mes({\cal X}_{w})^{-1}\int_{K_{{\bf m}}\times K_{n',n''}^{\pm}}\int_{{\cal X}_{w}}\hat{\phi}_{Lie}(x^{-1}Y_{w}x)\,dY_{w}\,dx.$$
Rappelons que l'on a suppos\'e $sgn_{CD}(w'')=1$ si $\sharp=iso$ et $sgn_{CD}(w'')=-1$ si $\sharp=an$.   Notons  $W({\bf m},n',n'')_{ell,\sharp}$ le sous-ensemble des \'el\'ements de $W({\bf m},n',n'')_{ell}$ dont la composante $w''$  v\'erifie cette condition. En revenant \`a (2), on obtient
$$(9) \qquad I=\vert W({\bf m},n',n'')\vert ^{-1}\sum_{w\in W({\bf m},n',n'')_{ell,\sharp}}res_{{\bf m}}(\kappa_{\pi})_{\gamma_{n',n''}}(w) sgn(w)q^{n/2}\vert {\bf T}_{w}({\mathbb F}_{q})\vert ^{-1}$$
$$mes({\cal X}_{w})^{-1}\int_{K_{{\bf m}}\times K_{n',n''}^{\pm}}\int_{{\cal X}_{w}}\hat{\phi}_{Lie}(x^{-1}Y_{w}x)\,dY_{w}\,dx.$$
   
  Notons plus pr\'ecis\'ement $I_{n',n''}(\phi)$ cette expression.   En 2.1, on a d\'efini un terme $\Theta_{\pi,cusp}^M(\phi)$. On a l'\'egalit\'e
  $$\Theta_{\pi,cusp}^M(\phi)=\sum_{n',n''}mes((A_{M}(F)K_{{\bf m}}\times K_{n',n''}^{\pm})/A_{M}(F))^{-1}\int_{A_{M}(F)\backslash M(F)}I_{n',n''}(^m\phi)\,dm,$$
  o\`u $(n',n'')$ parcourt $D(n)$ avec la restriction $n''\geq1$ si $\sharp=an$ et o\`u on a not\'e $^m\phi$ la fonction $x\mapsto \phi(m^{-1}xm)$. On peut oublier la restriction sur $n''$: si $\sharp=an$ et $n''=0$, la formule (9) vaut $0$ car l'ensemble $W({\bf m},n',n'')_{ell,\sharp}$ est vide. On voit que l'image par transformation de Fourier de $(^m\phi)_{Lie}$ est $^m(\hat{\phi}_{Lie})$. Les int\'egrales sur $K_{{\bf m}}\times K_{n',n''}^{\pm}$ de la formule (9) sont absorb\'ees par l'int\'egrale sur $A_{M}(F)\backslash M(F)$, mais introduisent des facteurs $mes(K_{{\bf m}}\times K_{n',n''}^{\pm})$. Notons $A_{M}(F)^c$ le plus grand sous-groupe compact de $A_{M}(F)$. C'est aussi l'intersection de $A_{M}(F)$ et de $K_{{\bf m}}$. Donc 
  $$mes((A_{M}(F)K_{{\bf m}}\times K_{n',n''}^{\pm})/A_{M}(F))=mes(K_{{\bf m}}\times K_{n',n''}^{\pm})mes(A_{M}(F)^c)^{-1}.$$
  D'o\`u
   $$(10) \qquad \Theta_{\pi,cusp}^M(\phi)=\sum_{(n',n'')\in D(n)}mes(A_{M}(F)^c)\vert W({\bf m},n',n'')\vert ^{-1}\sum_{w\in W({\bf m},n',n'')_{ell,\sharp}}$$
   $$res_{{\bf m}}(\kappa_{\pi})_{\gamma_{n',n''}}(w) sgn(w)q^{n/2}\vert {\bf T}_{w}({\mathbb F}_{q})\vert ^{-1}mes({\cal X}_{w})^{-1}\int_{A_{M}(F)\backslash M(F)}\int_{{\cal X}_{w}}\hat{\phi}_{Lie}(m^{-1}Y_{w}m)\,dY_{w}\,dm.$$
On peut encore remplacer l'int\'egrale en $A_{M}(F)\backslash M(F)$ par une int\'egrale sur $T_{w}(F)\backslash M(F)$, \`a condition de multiplier par $mes(A_{M}(F)\backslash T_{w}(F))$. Notons $T_{w}(F)^c$ le plus grand sous-groupe compact de $T_{w}(F)$. Parce que $T_{w}$ est non ramifi\'e, on a $T_{w}(F)=A_{M}(F)T_{w}(F)^c$, d'o\`u $mes(A_{M}(F)\backslash T_{w}(F))=mes(A_{M}(F)^c)^{-1}mes(T_{w}(F)^c)$. Le premier facteur compense son inverse qui figure dans la formule ci-dessus.   On a introduit plus haut le sous-groupe $T_{w}(F)^{u}$ de $T_{w}(F)$ et on a $T_{w}(F)^c/T_{w}(F)^{u}\simeq {\bf T}_{w}({\mathbb F}_{q})$. De plus, $mes(T_{w}(F)^{u})=mes(\varpi\mathfrak{t}_{w}(\mathfrak{o}))$. On a fix\'e sur $\mathfrak{t}_{w}(F)$ la mesure autoduale. Puisque $\mathfrak{t}_{w}(\mathfrak{o})$ et $\varpi\mathfrak{t}_{w}(\mathfrak{o})$ sont duaux pour le bicaract\`ere $(X,Y)\mapsto \psi_{F}(trace(XY))$, on calcule  
$$mes(\varpi\mathfrak{t}_{w}(\mathfrak{o}))=[\mathfrak{t}_{w}(\mathfrak{o}):\varpi\mathfrak{t}_{w}(\mathfrak{o})]^{-1/2}=q^{-n/2}.$$
D'o\`u
$$mes(T_{w}(F)^c)=q^{-n/2}\vert {\bf T}_{w}({\mathbb F}_{q})\vert .$$
Ces termes compensent leurs inverses figurant dans la formule (10). Finalement
  $$(11) \qquad \Theta_{\pi,cusp}^M(\phi)=\sum_{(n',n'')\in D(n_{0})} \vert W({\bf m},n',n'')\vert ^{-1}\sum_{w\in W({\bf m},n',n'')_{ell,\sharp}}res_{{\bf m}}(\kappa_{\pi})_{\gamma_{n',n''}}(w) $$
   $$sgn(w) mes({\cal X}_{w})^{-1}\int_{T_{w}(F)\backslash M(F)}\int_{{\cal X}_{w}}\hat{\phi}_{Lie}(m^{-1}Y_{w}m)\,dY_{w}\,dm.$$
   
   Soit maintenant $f\in C_{c}^{\infty}(G_{\sharp}(F))$. On suppose que le support de $f$ est form\'e d'\'el\'ements topologiquement unipotents.  En 2.1, on a d\'efini le terme 
   $$\Theta_{\pi,{\bf m},cusp}(f)=\int_{P(F)\backslash G_{\sharp}(F)}\Theta^M_{\pi,cusp}((^gf)_{U})\,dg. $$
   On d\'efinit la fonction $f_{Lie}$, cf.  2.2.  Posons $\phi=f_{U} $. On voit que 
   $$\phi_{Lie}(X)=\int_{\mathfrak{u}(F)}f_{Lie}(X+Y)\,dY ,$$
o\`u   $dY$ est la mesure de Haar sur $\mathfrak{u}(F)$ telle que l'exponentielle de $\mathfrak{u}(F)$ sur $U(F)$ pr\'eserve les mesures. D'o\`u aussi
$$\hat{\phi}_{Lie}(X)=\int_{\mathfrak{u}(F)}\hat{f}_{Lie}(X+Y)\,dY.$$
Ou encore
$$\hat{\phi}_{Lie}(X)=D(X)^{1/2}\int_{U(F)}\hat{f}_{Lie}(u^{-1}Xu)\,du,$$
o\`u $D$ est le discriminant de Weyl. On peut remplacer $f$ par $^gf$. En posant $\phi=(^gf)_{U}$, on obtient
$$\hat{\phi}_{Lie}(X)=D(X)^{1/2}\int_{U(F)}\hat{f}_{Lie}(g^{-1}u^{-1}Xug)\,du.$$
Les \'el\'ements $Y_{w}$ intervenant dans (11) v\'erifient $D(Y_{w})=1$ car les valeurs propres des r\'eductions ${\bf X}_{w}$ sont toutes distinctes.   On en d\'eduit  \'egalit\'e
  $$\int_{P(F)\backslash G_{\sharp}(F)}  \int_{T_{w}(F)\backslash M(F)} \hat{\phi}_{Lie}(m^{-1}Y_{w}m)\,dm\,dg=\int_{T_{w}(F)\backslash G_{\sharp}(F)}\hat{f}_{Lie}(g^{-1}Y_{w}g)\,dg.$$
  On obtient
 $$\Theta_{\pi,{\bf m},cusp}(f)=\sum_{(n',n'')\in D(n_{0})} \vert W({\bf m},n',n'')\vert ^{-1}\sum_{w\in W({\bf m},n',n'')_{ell,\sharp}}res_{{\bf m}}(\kappa_{\pi})_{\gamma_{n',n''}}(w) $$
   $$sgn(w) mes({\cal X}_{w})^{-1}\int_{T_{w}(F)\backslash G(F)}\int_{{\cal X}_{w}}\hat{f}_{Lie}(g^{-1}Y_{w}g)\,dY_{w}\,dg.$$
   
   Dans \cite{W2} page 53, on a introduit la distribution
   $$\varphi\mapsto \int_{T_{w}(F)\backslash G_{\sharp}(F)}\varphi(g^{-1}Y_{w}g)\,dg$$
   sur $C_{c}^{\infty}(\mathfrak{g}_{\sharp}(F))$ (elle y est not\'ee $\phi_{\theta}(X_{T},\varphi)$). On a montr\'e en \cite{W2} corollaire III.5 que sa restriction \`a un certain sous-espace ${\cal H}\subset C_{c}^{\infty}(\mathfrak{g}_{\sharp}(F))$ ne d\'ependait pas de l'\'el\'ement $Y_{w}\in {\cal X}_{w}$ (elle ne d\'epend d'ailleurs pas non plus du choix de ${\bf X}_{w}$ mais cela r\'esulte d\'ej\`a de nos calculs ci-dessus). L'espace ${\cal H}$ est  d\'efini ainsi. Soit $B$ un  sous-groupe d'Iwahori  de $G_{\sharp}(F)$. Il lui correspond un sous-$\mathfrak{o}$-r\'eseau $\mathfrak{b}$ de $\mathfrak{g}_{\sharp}(F)$.  Notons $C_{c}^{\infty}(g_{\sharp}(F)/\mathfrak{b})$ le sous-espace des fonctions invariantes par $\mathfrak{b}$. Alors ${\cal H}$ est la somme de ces espaces   $C_{c}^{\infty}(g_{\sharp}(F)/\mathfrak{b})$ quand $\mathfrak{b}$ d\'ecrit tous les sous-groupes d'Iwahori de $G_{\sharp}$. On voit facilement que ${\cal H}$ est exactement le sous-espace des fonctions $\varphi$ telles que $\hat{\varphi}$ soit \`a support topologiquement nilpotent. En particulier, $\hat{f}_{Lie}$ appartient \`a ${\cal H}$. Donc l'int\'egrale  
   $$\int_{T_{w}(F)\backslash G(F)} \hat{f}_{Lie}(g^{-1}Y_{w}g)\,dg$$
   ne d\'epend pas du point $Y_{w}\in {\cal X}_{w}$. Int\'egrer cette formule en $Y_{w}$ revient \`a la multiplier par $mes({\cal X}_{w})$ et ce facteur compense son inverse figurant dans la formule plus haut. On obtient simplement
    $$(12) \qquad \Theta_{\pi,{\bf m},cusp}(f)=\sum_{(n',n'')\in D(n_{0})} \vert W({\bf m},n',n'')\vert ^{-1}\sum_{w\in W({\bf m},n',n'')_{ell,\sharp}}res_{{\bf m}}(\kappa_{\pi})_{\gamma_{n',n''}}(w) $$
   $$sgn(w)  \int_{T_{w}(F)\backslash G(F)} \hat{f}_{Lie}(g^{-1}X_{w}g) \,dg.$$

   Pour expliciter davantage la formule obtenue, introduisons l'\'el\'ement $\boldsymbol{\gamma}_{0}=(0,0,n)$ de $\boldsymbol{\Gamma}$  (cf. \cite{W3} 1.8) et la composante $\kappa_{\pi,0}$ de $\kappa_{\pi}$ dans la composante ${\cal R}(\boldsymbol{\gamma}_{0})$ de ${\cal R}$. Soit $(\alpha,\beta',\beta'')\in {\cal P}_{3}(n)$. On a d\'efini en \cite{W3} 1.8  la valeur $\kappa_{\pi,0}(w_{\alpha,\beta',\beta''})$. Associons \`a notre triplet de partitions la partition ${\bf m}=\alpha$ et les entiers $n'=S(\beta')$, $n''=S(\beta'')$, $n_{0}=n'+n''$. Soit $w=(w_{1},...,w_{t},w',w'')$ un \'el\'ement  de $W({\bf m},n',n'')_{ell}$ tel que $w'$ et $w''$ soient param\'etr\'es par les partitions  $(\emptyset,\beta')$, resp. $(\emptyset,\beta'')$. Les d\'efinitions entra\^{\i}nent que $res_{{\bf m}}(\kappa_{\pi})_{\gamma_{n',n''}}(w)=\kappa_{\pi,0}(w_{\alpha,\beta',\beta''})$. Posons $ sgn(w_{\alpha,\beta',\beta''})=sgn(w)$ et d\'efinissons une distribution  $\phi_{\alpha,\beta',\beta''}$ sur $C_{c}^{\infty}(\mathfrak{g}_{\sharp}(F))$ par
   
   si $\sharp=iso$ et $l(\beta'')$ est impair ou si $\sharp=an$ et $l(\beta'')$ est pair, $ \phi_{\alpha,\beta',\beta''}=0$;
   
    si $\sharp=iso$ et $l(\beta'')$ est pair ou si $\sharp=an$ et $l(\beta'')$ est impair,
   $$\phi_{\alpha,\beta',\beta''}(\varphi)= \int_{T_{w}(F)\backslash G(F)}  \varphi(g^{-1}X_{w}g) \,dg$$
   pour tout $\varphi\in C_{c}^{\infty}(\mathfrak{g}_{\sharp}(F))$. 
   
   La distinction entre les deux cas provient de ce que $X_{w}$ n'existe que si $w\in W({\bf m},n',n'')_{ell,\sharp}$, c'est-\`a-dire si $sgn(w'')$ vaut $1$ si $\sharp=iso$, $-1$ si $\sharp=an$, ce qui se traduit par les conditions indiqu\'ees.  Cette d\'efinition d\'epend des choix de $w$ dans sa classe de conjugaison et de l'\'el\'ement $X_{w}$. Mais nous n'appliquerons cette distribution qu'\`a des \'el\'ements de l'espace ${\cal H}$. Comme on l'a dit ci-dessus, cette restriction ne d\'epend pas de ces choix. 
     Dans la formule (12), l'int\'egrale devient $\phi_{\alpha,\beta',\beta''}(\hat{f}_{Lie})$. 
   Cette formule  devient une somme index\'ee par les triplets  $(\alpha,\beta',\beta'')$  tels que $\alpha={\bf m}$ de termes ne d\'ependant que de ces triplets. Chaque triplet $(\alpha,\beta',\beta'')$ intervient avec une certaine multiplicit\'e. Celle-ci est le produit de $\vert W({\bf m},n',n'')\vert ^{-1}$ et du nombre d'\'el\'ements $w=(w_{1},...,w_{t},w',w'')\in W({\bf m},n',n'')_{ell}$ tels que $w'$ et $w''$ soient param\'etr\'es par  $(\emptyset,\beta')$, resp. $(\emptyset,\beta'')$.  Pour toute partition $\lambda$, posons
   $$z(\lambda)=(\prod_{j=1,...,l(\lambda)}2\lambda_{j})\prod_{i\geq1}mult_{\lambda}(i)!,$$
   et posons
   $$z(\alpha,\beta',\beta'')=   z(\alpha)z(\beta')z(\beta'').$$
   On voit que la multiplicit\'e pr\'ec\'edente est \'egale \`a 
   $$2^{l(\alpha)}mult!_{\alpha}z(\alpha,\beta',\beta'')^{-1}.$$ 
   Alors (12) se r\'ecrit
   $$\Theta_{\pi,{\bf m},cusp}(f)=\sum_{(\alpha,\beta',\beta'')\in {\cal P}_{3}(n); \alpha={\bf m}}2^{l(\alpha)}mult!_{\alpha}z(\alpha,\beta',\beta'')^{-1}$$
   $$sgn(w_{\alpha,\beta',\beta''})\kappa_{\pi,0}(w_{\alpha,\beta',\beta''})\phi_{\alpha,\beta',\beta''}(\hat{f}_{Lie}).$$
   Le r\'esultat de 2.1 est que $\Theta_{\pi}(f)$ est la somme sur ${\bf m}$  des expressions ci-dessus, multipli\'ees par $2^{-l({\bf m})}mult!_{{\bf m}}^{-1}$.  D'o\`u
   $$(13) \qquad \Theta_{\pi}(f)=\sum_{(\alpha,\beta',\beta'')\in {\cal P}_{3}(n) } z(\alpha,\beta',\beta'')^{-1}sgn(w_{\mu,\beta',\beta''})\kappa_{\pi,0}(w_{\alpha,\beta',\beta''})\phi_{\alpha,\beta',\beta''}(\hat{f}_{Lie}).$$
   
   \bigskip
   
   \section{Fronts d'onde}
   
   \bigskip
   
   \subsection{Rappel sur les orbites unipotentes}
   Soit $\sharp=iso$ ou $an$. On appelle orbite nilpotente une  classe de conjugaison par $G_{\sharp}(F)$ d'\'el\'ements nilpotents dans $\mathfrak{g}_{\sharp}(F)$. On note $Nil_{\sharp}$ l'ensemble des orbites nilpotentes. Les orbites nilpotentes  sont classifi\'ees par des donn\'ees $(\mu,(q_{i})_{i\in Jord_{bp}(\mu)})$ o\`u:
   
    $\mu\in {\cal P}^{orth}(2n+1)$;
    
    pour tout $i\in Jord_{bp}(\mu)$, $q_{i}$ est une classe d'\'equivalence d'une forme quadratique non d\'eg\'en\'er\'ee sur un espace vectoriel sur $F$ de dimension $mult_{\mu}(i)$; 
    
   le noyau anisotrope de la forme quadratique $\oplus_{i\in Jord_{bp}(\mu)}$  est \'equivalent \`a celui de $Q_{\sharp}$.
   
   Pour une orbite nilpotente ${\cal O}$, on note $\mu({\cal O})$ la partition associ\'ee \`a ${\cal O}$.
   
   Une classification analogue vaut pour les groupes ${\bf SO}(2n+1)$ et ${\bf SO}(2n)_{\sharp}$ d\'efinis sur ${\mathbb F}_{q}$.  Il y a une petite perturbation dans le cas du groupe ${\bf SO}(2n)_{iso}$. La classification ci-dessus vaut pour les classes de conjugaison par ${\bf O}(2n)_{iso}({\mathbb F}_{q})$ et non pas par ${\bf SO}(2n)_{iso}({\mathbb F}_{q})$. Il peut y avoir des classes de conjugaison par  ${\bf O}(2n)_{iso}({\mathbb F}_{q})$ qui se coupent en deux classes de conjugaison par ${\bf SO}(2n)_{iso}({\mathbb F}_{q})$. A ces deux classes sont  associ\'ees les m\^emes donn\'ees $(\mu,(q_{i})_{i\in Jord_{bp}(\mu)})$.
   
   La d\'efinition suivante va nous \^etre utile. Consid\'erons deux espaces vectoriels $l_{1}$ et $l_{2}$ sur ${\mathbb F}_{q}$ de dimensions $d_{1}$, resp. $d_{2}$. Soient $q_{1}$, resp. $q_{2}$, des formes quadratiques non d\'eg\'en\'er\'ees  sur  ces espaces. A isomorphisme pr\`es, il existe un unique triplet $(V,Q,L)$ v\'erifiant les conditions suivantes:
   
   $V$ est un espace vectoriel sur $F$ de dimension $d_{1}+d_{2}$;
   
   $Q$ est une forme quadratique non d\'eg\'en\'er\'ee sur $V$;
   
   $L\subset V$ est un r\'eseau presque autodual, c'est-\`a-dire $L^*\subset L\subset \varpi L^*$;
   
   $(l_{1},q_{1})$ est isomorphe \`a $(l',Q')$ et $(l_{2},q_{2})$ est isomorphe \`a $(l'',Q'')$  (rappelons que $l'=L/\varpi L^*$, $l''=L^*/L$ et que $Q'$ et $Q''$ sont les formes sur ces espaces qui se d\'eduisent naturellement de $Q$, cf. \cite{W3} 1.1). 
   
   On note $Q_{q_{1},q_{2}}$ cette forme quadratique $Q$ dont la  classe d'\'equivalence est bien d\'etermin\'ee.

    Consid\'erons l'ensemble ${\bf Nil}_{\sharp}$ des paires $({\cal O}_{1},{\cal O}_{2})$ telles que
   
  il existe  $n_{1},n_{2}\in {\mathbb N}$, avec $n_{1}+n_{2}=n$ et $n_{2}\geq1$ si $\sharp=an$, de sorte que
   ${\cal  O}_{1}$ soit une orbite nilpotente dans $\boldsymbol{\mathfrak{so}}(2n_{1}+1)({\mathbb F}_{q})$ et ${\cal  O}_{2}$ est une orbite nilpotente dans $\boldsymbol{\mathfrak{so}}(2n_{2})_{\sharp}({\mathbb F}_{q})$.
   
   A un telle paire, on  va associer une orbite nilpotente ${\cal O}_{{\cal O}_{1},{\cal O}_{2}}$.  Notons $(\mu_{1},(q_{1,i})_{i\in Jord_{bp}(\mu_{1})})$  et $(\mu_{2},(q_{2,i})_{i\in Jord_{bp}(\mu_{2})})$ les param\`etres de ${\cal O}_{1}$ et ${\cal O}_{2}$. On pose $\mu=\mu_{1}\cup \mu_{2}$ et, pour tout $i\in Jord_{bp}(\mu)$, $q_{i}=Q_{q_{1,i},q_{2,i}}$ (avec $q_{1,i}$ ou $q_{2,i}=0$ si $i\not\in Jord_{bp}(\mu_{1})$ ou $i\not\in Jord_{bp}(\mu_{2})$). On v\'erifie que $(\mu,(q_{i})_{i\in Jord_{bp}(\mu)})$ classifie une orbite nilpotente dans $\mathfrak{g}_{\sharp}(F)$. Alors ${\cal O}_{{\cal O}_{1},{\cal O}_{2}}$ est cette orbite unipotente.L'application 
   $$\begin{array}{ccc}{\bf Nil}_{\sharp}&\to&Nil_{\sharp}\\ ({\cal O}_{1},{\cal O}_{2})&\mapsto&{\cal O}_{{\cal O}_{1},{\cal O}_{2}}\\ \end{array}$$
   est surjective. 
   
    Pour ${\cal O}\in Nil_{\sharp}$, on note $I_{{\cal O}}$ l'int\'egrale orbitale associ\'ee \`a ${\cal O}$. Pour la d\'efinir, il faut bien s\^ur fixer une mesure sur ${\cal O}$ invariante par conjugaison. La d\'efinition de cette mesure n'aura pas d'importance pour nous.
    
    Soit $({\cal O}_{1},{\cal O}_{2})\in {\bf Nil}_{\sharp}$. En \cite{W2} IX.2, on a d\'efini une fonction $h_{{\cal O}_{1},{\cal O}_{2}}\in C_{c}^{\infty}(\mathfrak{g}_{\sharp}(F))$ (dans cette r\'ef\'erence, les \'el\'ements de ${\bf Nil}_{\sharp}$ \'etaient not\'es $\check{N}$). Elle v\'erifie les propri\'et\'es suivantes:
    
    $h_{{\cal O}_{1},{\cal O}_{2}}\in {\cal H}$; l'espace ${\cal H}$ a \'et\'e d\'efini en 2.3; c'est celui des fonctions $\varphi$ dont la transform\'ee de Fourier est \`a support topologiquement nilpotent;
    
  (1)   pour ${\cal O}\in Nil_{\sharp}$ dont l'adh\'erence $\bar{{\cal O}}$ ne contient pas ${\cal O}_{{\cal O}_{1},{\cal O}_{2}}$, $I_{{\cal O}}(h_{{\cal O}_{1},{\cal O}_{2}})=0$;
    
  (2)   pour ${\cal O}={\cal O}_{{\cal O}_{1},{\cal O}_{2}}$, $I_{{\cal O}}(h_{{\cal O}_{1},{\cal O}_{2}})\not=0$.
    
    Cf. \cite{W2} lemme IX.4. On d\'efinit une fonction $f_{{\cal O}_{1},{\cal O}_{2}}\in C_{c}^{\infty}(G_{\sharp}(F))$ comme suit: c'est la fonction \`a support topologiquement unipotent telle que $(f_{{\cal O}_{1},{\cal O}_{2}})_{Lie}=\hat{h}_{{\cal O}_{1},{\cal O}_{2}}$. 
    
   \bigskip
   
   \subsection{Developpement des caract\`eres \`a l'origine}
   
   Soit $\pi$ une repr\'esentation lisse et irr\'eductible de $G_{\sharp}(F)$.  D'apr\`es Harish-Chandra, on sait qu'il existe une unique famille de nombres complexes $(c_{{\cal O}}(\pi))_{{\cal O}\in Nil_{\sharp}}$  et un voisinage $V(\pi)$ de $1$ dans $G_{\sharp}(F)$ de sorte que les propri\'et\'es suivantes soit v\'erifi\'ees. Le voisinage $V(\pi)$ est invariant par conjugaison par $G_{\sharp}(F)$ et est form\'e d'\'el\'ements topologiquement unipotents.   Soit $f\in C_{c}^{\infty}(G_{\sharp}(F))$. On suppose que le support de $f$ est contenu dans $V(\pi)$. En particulier, on peut associer \`a $f$ une fonction $ f_{Lie}$ sur $\mathfrak{g}_{\sharp}(F)$, \`a support  topologiquement nilpotent. Alors on a l'\'egalit\'e
   $$(1) \qquad \Theta_{\pi}(f)=\sum_{{\cal O}\in Nil_{\sharp}}c_{{\cal O}}(\pi)I_{{\cal O}}(\hat{f}_{Lie}).$$
   Remarquons que les coefficients $c_{{\cal O}}(\pi)$ ne sont pas tous nuls. En effet,  si $f$ est la fonction caract\'eristique d'un sous-groupe ouvert compact $H$ contenu dans $V(\pi)$, $\Theta_{\pi}(f)$ est \'egal au produit de la mesure de $H$ et de la dimension du sous-espace des invariants par $H$ dans l'espace de $\pi$. Ce terme est non nul si $H$ est assez petit. 
   On dit que $\pi$ admet un front d'onde s'il existe $\mu(\pi)\in {\cal P}^{orth}(2n+1)$ de sorte que
   
   pour tout ${\cal O}\in Nil_{\sharp}$ tel que $c_{{\cal O}}(\pi)\not=0$, on a $\mu({\cal O})\leq \mu(\pi)$;
   
   il existe ${\cal O}\in Nil_{\sharp}$ tel que $c_{{\cal O}}(\pi)\not=0$ et $\mu({\cal O})=\mu(\pi)$.
   
   Evidemment, $\mu(\pi)$ est unique si elle existe.   On conjecture que toute repr\'esentation lisse irr\'eductible admet un front d'onde. Supposons que $\pi$ admette un front d'onde. On montre que
   
   $\mu(\pi)$ est une partition sp\'eciale, cf. \cite{M} th\'eor\`eme 1.4.
   
   pour tout ${\cal O}\in Nil_{\sharp}$ tel que  $\mu({\cal O})=\mu(\pi)$, on a $c_{{\cal O}}(\pi)\geq 0$, cf. \cite{MW2} corollaire 1.17. 
   
  {\bf Remarque.} La construction d'Harish-Chandra utilise l'exponentielle et non pas notre exponentielle tronqu\'ee $E$. Mais le r\'esultat est le m\^eme, avec les m\^emes coefficients, que l'on utilise l'une ou l'autre de ces applications.  
  
  \bigskip
  
  Dans le cas o\`u $\pi\in Irr_{unip,\sharp}$, on peut prendre pour voisinage $V(\pi)$ l'ensemble tout entier des \'el\'ements topologiquement unipotents de $G_{\sharp}(F)$. En effet, pour $f$ \`a support topologiquement unipotent,  $\Theta_{\pi}(f)$ est calcul\'e par la formule 2.3(13). Or il r\'esulte de \cite{D} th\'eor\`eme 2.1.5 que le membre de droite de cette formule est de la m\^eme forme que celui de (1) ci-dessus. Ces deux expressions  doivent co\"{\i}ncider si le support de $f$ est dans un voisinage assez petit de l'origine. Les distributions $f\mapsto I_{{\cal O}}(\hat{f}_{Lie})$  sont lin\'eairement ind\'ependantes, m\^eme si on les restreint aux fonctions v\'erifiant cette condition de support. Cela implique  que les coefficients sont les m\^emes dans les deux expressions. Donc (1) est valable pour toute $f$ \`a support topologiquement unipotent. 
   
   \bigskip
   
   \subsection{Le th\'eor\`eme}
   Pour $\sharp=iso$ ou $an$, notons $Irr_{tunip,\sharp}$ le sous-ensemble des repr\'esentations admissibles irr\'eductibles de $G_{\sharp}(F)$ qui sont temp\'er\'ees et de r\'eduction unipotente. Notons $Irr_{tunip}$ la r\'eunion disjointe de $Irr_{tunip,iso}$ et $Irr_{tunip,an}$. Dans \cite{W3} 1.3, on a adapt\'e l'habituelle classification de Langlands: l'ensemble $Irr_{tunip}$ est param\'etr\'e par un ensemble $\mathfrak{Irr}_{tunip}$ de triplets $(\lambda,s,\epsilon)$. En particulier, le  terme $\lambda$ est un \'el\'ement de ${\cal P}^{symp}(2n)$.  Pour un tel triplet, on a not\'e $\pi(\lambda,s,\epsilon)$ la repr\'esentation qui lui est  associ\'ee par Lusztig  (elle est temp\'er\'ee).  On a  introduit l'involution $D$ de Zelevinsky-Aubert-Schneider-Stuhler en  \cite{W3} 1.7 et une dualit\'e $d$ entre partitions en 1.6  et 1.7 ci-dessus.   On pose 
    $\delta(\lambda,s,\epsilon)= D(\pi(\lambda,s,\epsilon))$.

   \ass{Th\'eor\`eme}{Soit  $(\lambda,s,\epsilon)\in \mathfrak{Irr}_{tunip}$.  Alors  $\delta(\lambda,s,\epsilon)$ admet un front d'onde et on a l'\'egalit\'e $\mu(\delta(\lambda,s,\epsilon))=d(\lambda)$.} 
 
 La fin de l'article est consacr\'e \`a la d\'emonstration du th\'eor\`eme.  
   \bigskip
   
   \subsection{Une premi\`ere r\'eduction}
   En \cite{W3} 1.3, on a introduit le sous-ensemble $\mathfrak{Irr}_{unip-quad}$ des triplets $(\lambda,s,\epsilon)\in \mathfrak{Irr}_{tunip}$ tels que $s^2=1$.   Supposons que le th\'eor\`eme soit prouv\'e pour les triplets $(\lambda,s,\epsilon)\in \mathfrak{Irr}_{unip-quad}$ tels que $\lambda$ n'ait que des termes pairs. Montrons que le th\'eor\`eme r\'esulte de ce cas particulier.

   Soit $\sharp=iso$ ou $an$ et soit $P$ un sous-groupe  parabolique de $G_{\sharp}$, de composante de Levi $M$. Soit $\pi^M$ une repr\'esentation lisse irr\'eductible de $M(F)$, notons $\pi=Ind_{P}^{G_{\sharp}}(\pi^M)$ son induite et supposons $\pi$ irr\'eductible. On a d\'efini en 3.2 la notion de front d'onde pour le groupe $G_{\sharp}$ mais on sait bien que la d\'efinition est g\'en\'erale et s'applique en particulier au Levi $M$. Supposons que $\pi^M$ admette un front d'onde. Si 
   $$M=GL(m_{1})\times...\times GL(m_{t})\times G_{\sharp,n_{0}},$$
   $\mu(\pi^M)$ est alors une famille $(\mu_{1},...,\mu_{t},\mu_{0})$ o\`u, pour $i=1,...,t$, $\mu_{i}\in {\cal P}(m_{i})$ et $\mu_{0}\in {\cal P}^{orth}(2n_{0}+1)$. On a d\'efini en 1.8 une op\'eration d'induction qui envoie $\mu(\pi^M)$ sur une partition $ind(\mu(\pi^M))\in {\cal P}^{orth}(2n+1)$.

   Sous ces hypoth\`eses, on a
   
   (1) $\pi$ admet un front d'onde et on a $\mu(\pi)=ind(\mu(\pi^M))$.
   
   Preuve. Pour $f\in C_{c}^{\infty}(G_{\sharp}(F))$, on a l'\'egalit\'e $\Theta_{\pi}(f)=\Theta_{\pi^M}(f_{P})$, o\`u $f_{P}$ est l'habituel "terme constant" de $f$. On d\'efinit facilement un voisinage $V(\pi)$ de $1$ dans $G_{\sharp}(F)$ invariant par conjugaison et form\'e d'\'el\'ements topologiquement unipotents, de sorte que, si $f$ est \`a support dans $V(\pi)$, $f_{P}$ soit \`a support dans $V(\pi^M)$. Pour une telle fonction $f$, on a alors
   $$\Theta_{\pi}(f)=\sum_{{\cal O}^M\in Nil^M}c_{{\cal O}^M}(\pi^M)I_{{\cal O}^M}(f_{P}),$$
   o\`u $Nil^M$ est l'analogue de $Nil_{\sharp}$ pour le groupe $M$. 
   Notons $\mathfrak{u}$ le radical nilpotent de l'alg\`ebre de Lie $\mathfrak{p}$. Pour ${\cal O}^M\in Nil^M$ et ${\cal O}\in Nil_{\sharp}$, on dit que ${\cal O}$ est induite de ${\cal O}^M$ si ${\cal O}$ coupe ${\cal O}^M+\mathfrak{u}(F)$ selon un ouvert non vide. On note cette relation ${\cal O}\subset ind({\cal O}^M)$.  Si on se pla\c{c}ait sur la cl\^oture alg\'ebrique, $\bar{F}$ il y aurait une et une seule orbite induite mais, parce que l'on travaille sur $F$, il y en a plusieurs en g\'en\'eral. Cette op\'eration d'induction d'orbites est reli\'ee \`a l'induction des partitions par la relation suivante:
   
   si ${\cal O}\subset ind({\cal O}^M)$, alors $\mu({\cal O})=ind(\mu({\cal O}^M))$.
      
   On a une \'egalit\'e
   $$I_{{\cal O}^M}(f_{P})=\sum_{{\cal O}\subset ind({\cal O}^M)}c_{{\cal O}^M,{\cal O}}I_{{\cal O}}(f),$$
   avec des coefficients $c_{{\cal O}^M,{\cal O}}>0$. D'o\`u
   $$\Theta_{\pi}(f)=\sum_{{\cal O}\in Nil_{\sharp}}c_{{\cal O}}(\pi)I_{{\cal O}}(f),$$
   o\`u, pour tout ${\cal O}\in Nil_{\sharp}$, on a
   $$(2) \qquad c_{{\cal O}}(\pi)=\sum_{{\cal O}^M\in Nil^M; {\cal O}\subset ind( {\cal O}^M)}c_{{\cal O}^M}(\pi^M)c_{{\cal O}^M,{\cal O}}.$$
   Si $c_{{\cal O}}(\pi)\not=0$, il  existe ${\cal O}^M$ tel que ${\cal O}\subset ind({\cal O}^M)$ et $c_{{\cal O}^M}(\pi^M)\not=0$. On a alors $\mu({\cal O})=ind(\mu({\cal O}^M))$ et ${\cal O}^M\leq \mu(\pi^M)$.   D'o\`u $\mu({\cal O})\leq ind(\mu(\pi^M))$ car l'op\'eration d'induction est croissante. Inversement, soit ${\cal O}_{0}^M\in Nil^M$ tel que $c_{{\cal O}_{0}^M}(\pi^M)\not=0$ et $\mu({\cal O}_{0}^M)=\mu(\pi^M)$. Soit ${\cal O}\subset ind({\cal O}_{0}^M)$. On a $\mu({\cal O})=ind(\mu(\pi^M))$. Montrons que $c_{{\cal O}}(\pi)\not=0$. Soit ${\cal O}^M$ intervenant de fa\c{c}on non nulle dans la formule (2). On a $c_{{\cal O}^M}(\pi^M)\not=0$ donc $\mu({\cal O}^M)\leq \mu(\pi^M)$. Si cette relation n'est pas une \'egalit\'e, on a $ind(\mu({\cal O}^M))< ind(\mu(\pi^M))$ car l'op\'eration d'induction est strictement croissante. Cela contredit la relation $ind(\mu({\cal O}^M))=\mu({\cal O})=ind(\mu(\pi^M))$. Donc $\mu({\cal O}^M)=\mu(\pi^M)$. Alors le coefficient $c_{{\cal O}^M}(\pi^M)c_{{\cal O}^M,{\cal O}}$ est strictement
 positif. Par construction, il y a au moins une telle orbite $ {\cal O}^M$, \`a savoir ${\cal O}_{0}^M$. Le coefficient   $c_{{\cal O}}(\pi)$ est une somme non vide de termes strictement positifs, donc $c_{{\cal O}}(\pi)>0$, ce qui ach\`eve la d\'emonstration de (1). 
 
 Soit  $ (\lambda,s,\epsilon)\in \mathfrak{Irr}_{tunip,\sharp}$.  En reprenant les consid\'erations de \cite{W3} 1.3, on voit qu'il existe 
 
 - un sous-groupe parabolique $P$ de $G_{\sharp}$ de composante de Levi 
  $$M=GL(m_{1})\times...\times GL(m_{t})\times G_{\sharp,n_{0}};$$
  
 - pour tout $j=1,...,t$, un caract\`ere non ramifi\'e $\chi_{j}$ de $F^{\times}$;
  
  -  un \'el\'ement $ (\lambda_{0},s_{0},\epsilon_{0})\in \mathfrak{Irr}_{unip-quad,n_{0},\sharp}$ tel que $\lambda_{0}$ n'ait que des termes pairs; 
  
  de sorte que les propri\'et\'es suivantes soient v\'erifi\'ees:
  
  $$(3) \qquad \pi(\lambda,s,\epsilon)=Ind_{P}^{G_{\sharp}}(st_{m_{1}}(\chi_{1}\circ det)\otimes...\otimes st_{m_{t}}(\chi_{t}\circ det)\otimes \pi(\lambda_{0},s_{0},\epsilon_{0} )),$$
   o\`u $st_{m_{j}}$ est la repr\'esentation de Steinberg de $GL(m_{j};F)$;
   
 $$(4) \qquad \lambda =(m_{1} ,m_{1})\cup...\cup(m_{t} ,m_{t})\cup \lambda_{0} .$$
 
 En appliquant l'involution  $ D$, on d\'eduit de (3) l'\'egalit\'e

   $$\delta( \lambda,s,\epsilon)=Ind_{P}^{G_{\sharp}}(\delta^M),$$
   o\`u
   $$\delta^M= (\chi_{1}\circ det)\otimes...\otimes  (\chi_{t}\circ det)\otimes \delta(\lambda_{0} ,s_{0},\epsilon_{0} )).$$
   Pusiqu'on suppose connue le th\'eor\`eme pour $\delta(\lambda_{0} ,s_{0},\epsilon_{0})$ (et que les fronts d'onde des repr\'esentations des groupes $GL(m_{j})$ sont bien connus), il r\'esulte de (1) que $\delta( \lambda,s,\epsilon)$ admet un front d'onde et que 
   $$\mu(\delta( \lambda,s,\epsilon))=ind(\mu(\delta^M)).$$
    Le front d'onde d'un caract\`ere de $GL(m_{j};F)$ est $ ^t(m_{j})$, c'est-\`a-dire la partition compos\'ee de $m_{j}$ fois le nombre $1$. On a donc
$$ \mu(\delta ^M)=({^t(m}_{1}),...,{^t(m}_{t}),d(\lambda_{0} )).$$
Posons $\boldsymbol{\lambda}=((m_{1}),...,(m_{t}),\lambda_{0} )$. Avec les d\'efinitions de 1.8, cette derni\`ere relation  s'\'ecrit $\mu(\delta^M)=d(\boldsymbol{\lambda})$ tandis que l'\'egalit\'e (4) s'\'ecrit $\lambda=cup(\boldsymbol{\lambda})$. En appliquant le lemme 1.8, on obtient $ind(\mu(\delta^M))=d(\lambda )$. 
 D'o\`u $\mu(\delta( \lambda,s,\epsilon))=d(\lambda )$, ce qui d\'emontre le th\'eor\`eme.  
 
 \bigskip
 
 \subsection{Traduction de ce que l'on veut d\'emontrer}
 On  fixe d\'esormais un \'el\'ement $(\lambda,s,\epsilon)\in \mathfrak{Irr}_{tunip}$ et on pose $\delta=\delta(\lambda,s,\epsilon)$. On note $\sharp$ l'indice tel que $\delta$ soit une repr\'esentation de $G_{\sharp}(F)$.   On veut prouver que $\delta$ admet un front d'onde et que $\mu(\delta)=d(\lambda)$. Montrons qu'il suffit de prouver
 
 (1) pour tout $({\cal O}_{1},{\cal O}_{2})\in {\bf Nil}_{\sharp}$, la relation $\Theta_{\delta}(f_{{\cal O}_{1},{\cal O}_{2}})\not=0$ entra\^{\i}ne $\mu({\cal O}_{1})\cup\mu({\cal O}_{2})\leq d(\lambda)$;
 
 (2) il existe $({\cal O}_{1},{\cal O}_{2})\in {\bf Nil}_{\sharp}$ tel que $\Theta_{\delta}(f_{{\cal O}_{1},{\cal O}_{2}})\not=0$ et $\mu({\cal O}_{1})\cup\mu({\cal O}_{2})= d(\lambda)$.
   
   Comme on l'a dit en 3.2, on  peut prendre pour voisinage $V(\delta)$ l'ensemble tout entier des \'el\'ements topologiquement unipotents.  En particulier, le d\'eveloppement  de 3.2 vaut pour toute fonction $f_{{\cal O}_{1},{\cal O}_{2}}$. D'apr\`es la d\'efinition de cette fonction, on a
  $$\Theta_{\delta}(f_{{\cal O}_{1},{\cal O}_{2}})=\sum_{{\cal O}\in Nil_{\sharp}}c_{{\cal O}}(\delta)I_{{\cal O}}(h_{{\cal O}_{1},{\cal O}_{2}}).$$
   Soit ${\cal O}_{0}$ un \'el\'ement  maximal dans l'ensemble des ${\cal O}\in Nil_{\sharp}$ pour lesquels 
  $c_{{\cal O}}(\delta)\not=0$. Appliquons l'\'egalit\'e (3) \`a une paire $({\cal O}_{1},{\cal O}_{2})$ telle que ${\cal O}_{{\cal O}_{1},{\cal O}_{2}}={\cal O}_{0}$.   En vertu de 3.1(1), il ne reste dans la somme que des ${\cal O}$ pour lesquels $\bar{{\cal O}}$ contient ${\cal O}_{0}$. Par maximalit\'e de ${\cal O}_{0}$, il ne reste donc que ${\cal O}_{0}$. Le coefficient $c_{{\cal O}_{0}}(\delta)$ est non nul par hypoth\`ese et l'int\'egrale orbitale $I_{{\cal O}_{0}}( h_{{\cal O}_{1},{\cal O}_{2}})$ ne l'est pas par 3.1(2). Donc $\Theta_{\delta}(f_{{\cal O}_{1},{\cal O}_{2}})\not=0$. D'o\`u $\mu({\cal O}_{0})\leq d(\lambda)$ d'apr\`es (1). Ceci \'etant vrai pour tout \'el\'ement maximal ${\cal O}_{0}$, c'est vrai pour tout \'el\'ement: pour tout ${\cal O}$ tel que $c_{{\cal O}}(\delta)\not=0$, on a $\mu({\cal O})\leq d(\lambda)$. En appliquant maintenant (3) pour une paire $({\cal O}_{1},{\cal O}_{2})$ v\'erifiant (2), le m\^eme calcul montre que $c_{{\cal O}}(\delta)\not=0$ pour l'orbite ${\cal O}={\cal O}_{{\cal O}_{1},{\cal O}_{2}}$. Pour cette orbite, on a $\mu({\cal O})=d(\lambda)$. Cela v\'erifie les propri\'et\'es requises pour que $\delta$ admette un front d'onde et que l'on ait $\mu(\delta)=d(\lambda)$.

  \bigskip
  
   \subsection{D\'ebut du calcul}
  Fixons un couple $({\cal O}_{1},{\cal O}_{2})\in {\bf Nil}_{\sharp}$, posons $\mu_{1}=\mu({\cal O}_{1})$, $\mu_{2}=\mu({\cal O}_{2})$, $S(\mu_{1})=2n_{1}+1$, $S(\mu_{2})=2n_{2}$.  Si $\sharp=iso$, on pose $W_{n_{2},iso}=W_{n_{2}}^D$. Si $\sharp=an$, auquel cas $n_{2}>0$, on note $W_{n_{2},an}=\{w\in W_{n_{2}};sgn_{CD}(w)=-1\}$. On fixe des \'el\'ements nilpotents $Y_{1}\in {\cal O}_{1}$ et $Y_{2}\in {\cal O}_{2}$.
  
  La formule 2.3 (13) calcule $\Theta_{\delta}(f_{{\cal O}_{1},{\cal O}_{2}})$ en fonction de termes $\phi_{\alpha,\beta_{1},\beta_{2}}(h_{{\cal O}_{1},{\cal O}_{2}})$ pour $(\alpha,\beta_{1},\beta_{2})\in {\cal P}_{3}(n)$. On a calcul\'e ce terme en \cite{W2} proposition 3.5. On va rappeler ce r\'esultat en modifiant quelque peu ses notations. Pour $w_{1}\in W_{n_{1}}$, on d\'efinit une certaine fonction $Q_{w_{1}}^{\natural}$ sur l'ensemble des \'el\'ements nilpotents de ${\bf SO}(2n_{1}+1)({\mathbb F}_{q})$, cf. \cite{W2} VIII.13. Elle est invariante par conjugaison par ${\bf SO}(2n_{1}+1)({\mathbb F}_{q})$ et ne d\'epend que de la classe de conjugaison de $w_{1}$. Pour $w_{2}\in W_{n_{2},\sharp}$, on d\'efinit de m\^eme une fonction $Q_{w_{2}}^{\natural}$ sur l'ensemble des \'el\'ements nilpotents de ${\bf SO}(2n_{2})_{\sharp}({\mathbb F}_{q})$, cf. \cite{W2} VIII.13. Elle est invariante par conjugaison par ${\bf SO}(2n_{2})_{\sharp}({\mathbb F}_{q})$ et ne d\'epend que de la classe de conjugaison par $W_{n_{2}}^D$ de $w_{2}$. 
  
  {\bf Remarque.} Dans le cas o\`u $\sharp=an$, la construction de \cite{W2} \'etait un peu diff\'erente. On y avait fix\'e une certaine sym\'etrie \'el\'ementaire $w_{\phi}\in W_{n_{2},an}$ et d\'efini une fonction $Q_{w_{2}}^{\natural}$ index\'ee non pas par un \'el\'ement $w_{2}\in W_{n_{2},an}$, mais par un \'el\'ement $w_{2}\in W_{n_{2}}^D$. Cette fonction ne d\'ependait que de la classe de $w_{\phi}$-conjugaison de $w_{2}$. Notre pr\'esente fonction $Q_{w_{2}}^{\natural}$ est la fonction $Q_{w_{\phi}w_{2}}^{\natural}$ de \cite{W2}.
  
  \bigskip
  Notons $W(\alpha,\beta_{1},\beta_{2})$ l'ensemble des paires $(w_{1},w_{2})\in W_{n_{1}}\times W_{n_{2},\sharp}$ v\'erifiant la condition suivante. Notons $(\alpha_{1},\beta'_{1})$ la paire de partitions param\'etrant la classe de conjugaison de $w_{1}$ et $(\alpha_{2},\beta'_{2})$ celle qui param\`etre la classe de conjugaison par $W_{n_{2}}$ de $w_{2}$. Alors
  $$\alpha=\alpha_{1}\cup \alpha_{2},\,\, \beta'_{1}=\beta_{1},\,\, \beta'_{2}=\beta_{2}.$$
  Pour une telle paire $(w_{1},w_{2})$, posons
  
  $$[w_{1},w_{2}]=\frac{z(\alpha)}{z(\alpha_{1})z(\alpha_{2})},$$
  cf. 2.3 pour la d\'efinition de ces termes;
  
  $\eta(w_{1},w_{2})=2$ si $n_{2}\geq1$ et la classe de conjugaison de $w_{2}$ par $W_{n_{2}}$ co\"{\i}ncide avec sa classe de conjugaison par $W_{n_{2}}^D$;
  
  $\eta(w_{1},w_{2})=1$ si $n_{2}=0$ ou si $n_{2}\geq1$ et la classe de conjugaison de $w_{2}$ par $W_{n_{2}}$ se coupe en deux classes de conjugaison par $W_{n_{2}}^D$.
  
  Fixons    un ensemble de repr\'esentants  ${\cal W}(\alpha,\beta_{1},\beta_{2})$ des classes de conjugaison par $W_{n_{1}}\times W_{n_{2}}^D$ dans $W(\alpha,\beta_{1},\beta_{2})$. 
  
  La proposition 3.5 de \cite{W2} affirme alors l'existence d'un demi-entier $d_{{\cal O}_{1},{\cal O}_{2}}$ (ne d\'ependant que de $({\cal O}_{1},{\cal O}_{2})$) de sorte que
  $$\phi_{\alpha,\beta_{1},\beta_{2}}(h_{{\cal O}_{1},{\cal O}_{2}})=q^{d_{{\cal O}_{1},{\cal O}_{2}}}\sum_{(w_{1},w_{2})\in {\cal W}(\alpha,\beta_{1},\beta_{2})} \eta(w_{1},w_{2})[w_{1},w_{2}] Q_{w_{1}}^{\natural}(Y_{1})Q_{w_{2}}^{\natural}(Y_{2}).$$
  Cette formule peut se simplifier. Pour $w_{1}\in W_{n_{1}}$, notons $Z(w_{1})$ son centralisateur dans $W_{n_{1}}$. Pour $w_{2}\in W_{n_{2}}$, notons $Z^D(w_{2})$ son centralisateur dans $W_{n_{2}}^D$. Le terme $Q_{w_{1}}^{\natural}(Y_{1})Q_{w_{2}}^{\natural}(Y_{2})$ ne d\'ependant que de la classe de conjugaison de $(w_{1},w_{2})$ par $W_{n_{1}}\times W_{n_{2}}^D$, on peut remplacer la somme sur le syst\`eme de repr\'esentants ${\cal W}(\alpha,\beta_{1},\beta_{2})$ par une somme sur $W(\alpha,\beta_{1},\beta_{2})$, \`a condition de multiplier chaque terme index\'e par $(w_{1},w_{2})$ par l'inverse du nombre d'\'el\'ements de sa classe de conjugaison par  $W_{n_{1}}\times W_{n_{2}}^D$. Cet inverse est \'egal \`a
  $$\vert Z(w_{1})\vert \vert Z^D(w_{2})\vert \vert W_{n_{1}}\vert ^{-1}\vert W_{n_{2}}^D\vert ^{-1}.$$
  D'autre part,
  soit $(w_{1}  ,w_{2})\in W(\alpha,\beta_{1},\beta_{2})$. On v\'erifie l'\'egalit\'e
  $$z(\alpha,\beta_{1},\beta_{2})^{-1}\eta(w_{1},w_{2})[w_{1},w_{2}]=\vert Z(w_{1})\vert^{-1} \vert Z^D(w_{2})\vert ^{-1}.$$
  La formule ci-dessus se r\'ecrit
  $$z(\alpha,\beta_{1},\beta_{2})^{-1}\phi_{\alpha,\beta_{1},\beta_{2}}(h_{{\cal O}_{1},{\cal O}_{2}})=q^{d_{{\cal O}_{1},{\cal O}_{2}}}\vert W_{n_{1}}\vert ^{-1}\vert W_{n_{2}}^D\vert ^{-1}\sum_{(w_{1},w_{2})\in  W(\alpha,\beta_{1},\beta_{2})}Q_{w_{1}}^{\natural}(Y_{1})Q_{w_{2}}^{\natural}(Y_{2}).$$
  Reportons cette \'egalit\'e dans la formule 2.3 (13). On obtient
  $$\Theta_{\delta}(f_{{\cal O}_{1},{\cal O}_{2}})=q^{d_{{\cal O}_{1},{\cal O}_{2}}}\vert W_{n_{1}}\vert ^{-1}\vert W_{n_{2}}^D\vert ^{-1}\sum_{(\alpha,\beta_{1},\beta_{2})\in {\cal P}_{3}(n)}sgn(w_{\alpha,\beta_{1},\beta_{2}})\kappa_{\delta,0}(w_{\alpha,\beta_{1},\beta_{2}})$$
  $$\sum_{(w_{1},w_{2})\in  W(\alpha,\beta_{1},\beta_{2})}Q_{w_{1}}^{\natural}(Y_{1})Q_{w_{2}}^{\natural}(Y_{2}).$$
  Sommer en $(\alpha,\beta_{1},\beta_{2})$ puis en $(w_{1},w_{2})\in W(\alpha,\beta_{1},\beta_{2})$ revient \`a sommer sur tout $(w_{1},w_{2})\in W_{n_{1}}\times W_{n_{2},\sharp}$. D'o\`u
   $$ \Theta_{\delta}(f_{{\cal O}_{1},{\cal O}_{2}})=q^{d_{{\cal O}_{1},{\cal O}_{2}}}\vert W_{n_{1}}\vert ^{-1}\vert W_{n_{2}}^D\vert ^{-1}$$
   $$\sum_{(w_{1},w_{2})\in W_{n_{1}}\times W_{n_{2},\sharp}}sgn(w_{1})sgn(w_{2})\kappa_{\delta,0}(w_{1}\times w_{2})Q_{w_{1}}^{\natural}(Y_{1})Q_{w_{2}}^{\natural}(Y_{2}).$$
   Pour toute paire $(\rho_{1},\rho_{2})\in \hat{W}_{n_{1}}\times \hat{W}_{n_{2}}$, posons
   $$m_{\delta}(\rho_{1},\rho_{2})=\vert W_{n_{1}}\vert ^{-1}\vert W_{n_{2}}\vert ^{-1}\sum_{(w_{1},w_{n_{2}})\in W_{n_{1}}\times W_{n_{2}}}k_{\delta,0}(w_{1},w_{2})\rho_{1}(w_{1})\rho_{2}(w_{2}).$$
   Posons aussi
   $$\chi^{\natural}_{\rho_{1}}=\vert W_{n_{2}}\vert ^{-1}\sum_{w_{1}\in W_{n_{1}}}\rho_{1}(w_{1})Q_{w_{1}}^{\natural}$$
   et
   $$\chi^{\natural}_{\rho_{2},\sharp}=\vert W_{n_{2}}^D\vert ^{-1}\sum_{w_{2}\in W_{n_{2},\sharp}}\rho_{2}(w_{2})Q_{w_{2}}^{\natural}.$$
   Pour $(w_{1},w_{2})\in W_{n_{1}}\times W_{n_{2}}$, on a l'\'egalit\'e
   $$sgn(w_{1})sgn(w_{2})\kappa_{\delta,0}(w_{1}\times w_{2})=\sum_{(\rho_{1},\rho_{2})\in \hat{W}_{n_{1}}\times \hat{W}_{n_{2}}}m_{\delta}(\rho_{1}\otimes sgn,\rho_{2}\otimes sgn)\rho_{1}(w_{1})\rho_{2}(w_{2}).$$
   On obtient alors
    $$(1) \qquad  \Theta_{\delta}(f_{{\cal O}_{1},{\cal O}_{2}})=q^{d_{{\cal O}_{1},{\cal O}_{2}}}\sum_{(\rho_{1},\rho_{2})\in \hat{W}_{n_{1}}\times \hat{W}_{n_{2}}}m_{\delta}(\rho_{1}\otimes sgn,\rho_{2}\otimes sgn)\chi^{\natural}_{\rho_{1}}(Y_{1})\chi^{\natural}_{\rho_{2},\sharp}(Y_{2}).$$
    
   Soit $(\mu_{1},\eta_{1})\in \boldsymbol{{\cal P}^{orth}}(2n_{1}+1)$. Supposons $k_{\mu_{1},\eta_{1}}=1$. On a d\'efini en \cite{W2} VIII.13 une fonction $\chi_{\mu_{1},\eta_{1}}^{\natural}$ sur l'ensemble des \'el\'ements nilpotents de $\boldsymbol{\mathfrak{so}}(2n)({\mathbb F}_{q})$. Il existe un demi-entier $d_{\mu_{1},\eta_{1}}$ tel que $\chi_{\mu_{1},\eta_{1}}^{\natural}=q^{d_{\mu_{1},\eta_{1}}}\chi^{\natural}_{\rho_{\mu_{1},\eta_{1}}}$. On note $\boldsymbol{{\cal P}^{orth}}(2n_{1}+1;k=1)$ le sous-ensemble des $(\mu_{1},\eta_{1})\in \boldsymbol{{\cal P}^{orth}}(2n_{1}+1)$ tels que $k_{\mu_{1},\eta_{1}}=1$. Rappelons que l'application
   $$\begin{array}{ccc}\boldsymbol{{\cal P}^{orth}}(2n_{1}+1;k=1)&\to&\hat{W}_{n_{2}}\\ (\mu_{1},\eta_{1})&\mapsto& \rho_{\mu_{1},\eta_{1}}\\ \end{array}$$
   est bijective.
   
   Soit $(\underline{\mu}_{2},\eta_{2})\in \boldsymbol{\underline{{\cal P}}^{orth}}(2n_{2})$, cf. 1.5. Supposons $k_{\mu_{2},\eta_{2}}=0$. On a d\'efini en \cite{W2} VIII.13 une fonction $\chi_{\underline{\mu}_{2},\eta_{2},\sharp}^{\natural}$ sur l'ensemble des \'el\'ements nilpotents de $\boldsymbol{\mathfrak{so}}(2n)_{\sharp}({\mathbb F}_{q})$ (on a ajout\'e un indice $\sharp$ \`a la notation de \cite{W2}). Soit  $(\mu_{2},\eta_{2})\in \boldsymbol{{\cal P}^{orth}}(2n_{2})$ tel que  $k_{\mu_{2},\eta_{2}}=0$. Si  $\mu_{2}$ n'est pas exceptionnel, $(\mu_{2},\eta_{2})$ est aussi un \'el\'ement de $\boldsymbol{\underline{{\cal P}}^{orth}}(2n_{2})$; la repr\'esentation $\rho_{\mu_{2},\epsilon_{2}}$ est d\'efinie ainsi que ses prolongements $\rho^+_{\mu_{2},\epsilon_{2}}$ et $\rho^-_{\mu_{2},\epsilon_{2}}$, cf. 1.12. 
   La fonction $\chi^{\natural}_{\mu_{2},\eta_{2},\sharp}$ est aussi d\'efinie. Si  $\mu_{2}$ est exceptionnel, auquel cas $\eta_{2}=1$, $\mu_{2}$ se rel\`eve en les deux \'el\'ements $(\mu_{2},+)$ et $(\mu_{2},-)$ de $\underline{{\cal P}}^{orth}(2n_{2})$. Les repr\'esentations $\rho_{\mu_{2},+,\eta_{2}}$ et $\rho_{\mu_{2},-,\eta_{2}}$ sont d\'efinies. Les repr\'esentations $\rho^+_{\mu_{2},+,\eta_{2}}$, $\rho^-_{\mu_{2},+,\eta_{2}}$, $\rho^+_{\mu_{2},-,\eta_{2}}$ et $\rho^-_{\mu_{2},-,\eta_{2}}$ sont toutes \'egales. On note $\rho^+_{\mu_{2},\eta_{2}}=\rho^-_{\mu_{2},\eta_{2}}$ ce prolongement.
   On pose
   $$\chi^{\natural}_{\mu_{2},\eta_{2},\sharp}=\chi^{\natural}_{\mu_{2},+,\eta_{2},\sharp}+\chi^{\natural}_{\mu_{2},-,\eta_{2},\sharp}.$$
    En tout cas, il existe un demi-entier $d_{\mu_{2},\eta_{2}}$ tel que les \'egalit\'es suivantes soient v\'erifi\'ees:

   si $\sharp=iso$,
   $$\chi_{\mu_{2},\eta_{2},iso}^{\natural}=q^{d_{\mu_{2},\eta_{2}}}\chi^{\natural}_{\rho^+_{\mu_{2},\eta_{2}},iso}=q^{d_{\mu_{2},\eta_{2}}}\chi^{\natural}_{\rho^-_{\mu_{2},\eta_{2}},iso};$$

   si $\sharp=an$, 
   $$\chi_{\mu_{2},\eta_{2},an}^{\natural}=q^{d_{\mu_{2},\eta_{2}}}\chi^{\natural}_{\rho^+_{\mu_{2},\eta_{2}},an}=-q^{d_{\mu_{2},\eta_{2}}}\chi^{\natural}_{\rho^-_{\mu_{2},\eta_{2}},an}.$$
   
   Remarquons que, quand $\mu_{2}$ est exceptionnel, on a $\chi^{\natural}_{\mu_{2},\eta_{2},an}=0$. 
   
   On note $\boldsymbol{{\cal P}^{orth}}(2n_{2};k=0)$ l'ensemble des $(\mu_{2},\eta_{2})\in \boldsymbol{{\cal P}^{orth}}(2n_{2})$ tels que $k_{\mu_{2},\eta_{2}}=0$. Rappelons que $\hat{W}_{n_{2}}$ est r\'eunion disjointe des ensembles $\{\rho^+_{\mu_{2},\eta_{2}},\rho^-_{\mu_{2},\eta_{2}}\}$ pour $(\mu_{2},\eta_{2})\in \boldsymbol{{\cal P}^{orth}}(2n_{2};k=0)$ avec $\mu_{2}$ non exceptionnel et des ensembles $\{\rho^+_{\mu_{2},+,\eta_{2}}\}=\{\rho^-_{\mu_{2},+,\eta_{2}}\}=\{\rho^+_{\mu_{2},-,\eta_{2}}\}=\{\rho^-_{\mu_{2},-,\eta_{2}}\}$  pour $(\mu_{2},\eta_{2})\in \boldsymbol{{\cal P}^{orth}}(2n_{2};k=0)$ avec $\mu_{2}$  exceptionnel.
   
 Pour $(\mu_{1},\eta_{1};\mu_{2},\eta_{2})\in \boldsymbol{{\cal P}^{orth}}(2n_{1}+1;k=1)\times \boldsymbol{{\cal P}^{orth}}(2n_{2};k=0)$, posons
 
 si $\sharp=iso$, $n_{2}\not=0$ et $\mu_{2}$ n'est pas exceptionnel,
 $$m_{\delta,iso}(\mu_{1},\eta_{1};\mu_{2},\eta_{2})= m_{\delta}(\rho_{\mu_{1},\eta_{1}}\otimes sgn,\rho^+_{\mu_{2},\eta_{2}} \otimes sgn)+ m_{\delta}(\rho_{\mu_{1},\eta_{1}}\otimes sgn,\rho^-_{\mu_{2},\eta_{2}} \otimes sgn);$$
 
  si $\sharp=an$, $n_{2}\not=0$ et $\mu_{2}$ n'est pas exceptionnel,
   $$m_{\delta,an}(\mu_{1},\eta_{1};\mu_{2},\eta_{2})= m_{\delta}(\rho_{\mu_{1},\eta_{1}}\otimes sgn,\rho^+_{\mu_{2},\eta_{2}} \otimes sgn)- m_{\delta}(\rho_{\mu_{1},\eta_{1}}\otimes sgn,\rho^-_{\mu_{2},\eta_{2}} \otimes sgn);$$
   
   si $\sharp=iso$ et  si $n_{2}=0$ ou $\mu_{2}$ est exceptionnel
   $$m_{\delta,iso}(\mu_{1},\eta_{1};\mu_{2},\eta_{2})=\frac{ m_{\delta}(\rho_{\mu_{1},\eta_{1}}\otimes sgn,\rho^+_{\mu_{2},\eta_{2}} \otimes sgn)+ m_{\delta}(\rho_{\mu_{1},\eta_{1}}\otimes sgn,\rho^-_{\mu_{2},\eta_{2}} \otimes sgn)}{2};$$
   
   si $\sharp=an$ et  $\mu_{2}$ est exceptionnel
    $$m_{\delta,an}(\mu_{1},\eta_{1};\mu_{2},\eta_{2})=0.$$
 
 On voit alors que la formule (1) se r\'ecrit
 $$(2)\qquad \Theta_{\delta}(f_{{\cal O}_{1},{\cal O}_{2}})=q^{d_{{\cal O}_{1},{\cal O}_{2}}}\sum_{(\mu_{1},\eta_{1};\mu_{2},\eta_{2})\in \boldsymbol{{\cal P}^{orth}}(2n_{1}+1;k=1)\times \boldsymbol{{\cal P}^{orth}}(2n_{2})}q^{-d_{\mu_{1},\eta_{1}}-d_{\mu_{2},\eta_{2}}}$$
 $$m_{\delta,\sharp}(\mu_{1},\eta_{1};\mu_{2},\eta_{2})\chi^{\natural}_{\mu_{1},\eta_{1}}(Y_{1})\chi^{\natural}_{\mu_{2},\eta_{2},\sharp}(Y_{2}).$$
 
 \bigskip
 \subsection{Traduction des conditions en termes de repr\'esentations de groupes de Weyl}
 Montrons qu'il nous suffit de prouver les deux assertions suivantes:
 
 (1) soient $n_{1},n_{2}\in {\mathbb N}$ avec $n_{1}+n_{2}=n$ et $n_{2}\geq1$ si $\sharp=an$; soient $(\mu_{1},\eta_{1};\mu_{2},\eta_{2})\in \boldsymbol{{\cal P}^{orth}}(2n_{1}+1;k=1)\times \boldsymbol{{\cal P}^{orth}}(2n_{2};k=0)$; supposons $m_{\delta,\sharp}(\mu_{1},\eta_{1};\mu_{2},\eta_{2})\not=0$; alors $\mu_{1}\cup\mu_{2}\leq d(\lambda)$;
 
 (2) il existe $n_{1},n_{2}\in {\mathbb N}$ avec $n_{1}+n_{2}=n$ et $n_{2}\geq1$ si $\sharp=an$ et il existe $(\mu_{1},\eta_{1};\mu_{2},\eta_{2})\in \boldsymbol{{\cal P}^{orth}}(2n_{1}+1;k=1)\times \boldsymbol{{\cal P}^{orth}}(2n_{2};k=0)$ tels que $m_{\delta,\sharp}(\mu_{1},\eta_{1};\mu_{2},\eta_{2})\not=0$ et $\mu_{1}\cup\mu_{2}=d(\lambda)$.

Soit $({\cal O}_{1},{\cal O}_{2})\in {\bf Nil}_{\sharp}$. On en d\'eduit des entiers $n_{1}$, $n_{2}$ comme en 3.6. Supposons $\Theta_{\delta}(f_{{\cal O}_{1},{\cal O}_{2}})\not=0$.  La formule 3.6(2) implique qu'il existe $(\mu_{1},\eta_{1};\mu_{2},\eta_{2})\in \boldsymbol{{\cal P}^{orth}}(2n_{1}+1;k=1)\times \boldsymbol{{\cal P}^{orth}}(2n_{2};k=0)$ tel que $m_{\delta,iso}(\mu_{1},\eta_{1};\mu_{2},\eta_{2})\not=0$ et $\chi^{\natural}_{\mu_{1},\eta_{1}}(Y_{1})\chi^{\natural}_{\mu_{2},\eta_{2},\sharp}(Y_{2})\not=0$. Or la fonction $\chi^{\natural}_{\mu_{1},\eta_{1}}$ n'est non nulle que sur les orbites nilpotentes ${\cal O}'_{1}$ v\'erifiant $\mu({\cal O}'_{1})\leq \mu_{1}$, cf. \cite{W2} VIII.13. Puisque $Y_{1}\in {\cal O}_{1}$, cela entra\^{\i}ne  $\mu({\cal O}_{1})\leq \mu_{1}$. La fonction $\chi^{\natural}_{\mu_{2},\eta_{2},\sharp}$ v\'erifie une propri\'et\'e analogue. Donc $\mu({\cal O}_{2})\leq \mu_{2}$. Gr\^ace \`a (1), on a aussi $\mu_{1}\cup\mu_{2}\leq d(\lambda)$. Donc $\mu({\cal O}_{1})\cup \mu({\cal O}_{2})\leq d(\lambda)$, ce qui v\'erifie la propri\'et\'e (1) de 3.5. 

Fixons des donn\'ees v\'erifiant (2). Consid\'erons la somme
$$\Psi=\sum_{\eta'_{1},\eta'_{2}}q^{-d_{\mu_{1},\eta_{1}}-d_{\mu_{2},\eta_{2}}}m_{\delta,\sharp}(\mu_{1},\eta'_{1};\mu_{2},\eta'_{2})\chi_{\mu_{1},\eta'_{1}}^{\natural}\chi^{\natural}_{\mu_{2},\eta'_{2},\sharp},$$
o\`u $\eta'_{1}$ parcourt les \'el\'ements de $\{\pm 1\}^{Jord_{bp}(\mu_{1})}/\{\pm 1\}$ tels que $k(\mu_{1},\eta'_{1})=1$ et $\eta'_{2}$ parcourt les \'el\'ements de $\{\pm 1\}^{Jord_{bp}(\mu_{2})}/\{\pm 1\}$ tels que $k(\mu_{2},\eta'_{2})=1$.  C'est une fonction invariante par conjugaison sur le produit des ensembles d'\'el\'ements nilpotents de $\boldsymbol{\mathfrak{so}}(2n_{1}+1)({\mathbb F}_{q})$ et de $\boldsymbol{\mathfrak{so}}(2n_{2})_{\sharp}({\mathbb F}_{q})$. Notons ${\cal U}(\mu_{1})$ la r\'eunion des orbites nilpotentes ${\cal O}_{1}$ dans $\boldsymbol{\mathfrak{so}}(2n_{1}+1)({\mathbb F}_{q})$ telles que $\mu({\cal O}_{1})=\mu_{1}$. Notons ${\cal U}(\mu_{2})$ la r\'eunion des orbites nilpotentes ${\cal O}_{2}$ dans $\boldsymbol{\mathfrak{so}}(2n_{2})_{\sharp}({\mathbb F}_{q})$ telles que $\mu({\cal O}_{2})=\mu_{2}$. 
Quand $\eta'_{1}$ d\'ecrit les \'el\'ements ci-dessus, les restrictions \`a ${\cal U}(\mu_{1})$ des fonctions $\chi_{\mu_{1},\eta'_{1}}^{\natural}$ sont lin\'eairement ind\'ependantes, cf. \cite{W2} VIII.13. Quand $\eta'_{2}$ d\'ecrit les \'el\'ements ci-dessus, les restrictions \`a ${\cal U}(\mu_{2})$ des fonctions $\chi_{\mu_{2},\eta'_{2},\sharp}^{\natural}$ sont elles aussi lin\'eairement ind\'ependantes (on doit remarquer que, si $\sharp=an$, l'hypoth\`ese $m_{\delta,\sharp}(\mu_{1},\eta_{1};\mu_{2},\eta_{2})\not=0$ implique que $\mu_{2}$ n'est pas exceptionnel).  L'hypoth\`ese $m_{\delta,\sharp}(\mu_{1},\eta_{1};\mu_{2},\eta_{2})\not=0$ implique donc que la restriction de $\Psi$  \`a ${\cal U}(\mu_{1})\times {\cal U}(\mu_{2})$ est non nulle. Fixons donc des orbites ${\cal O}_{1}\subset {\cal U}_{1}$ et ${\cal O}_{2}\subset {\cal U}_{2}$ telles que $\Psi$ soit non nulle sur ${\cal O}_{1}\times {\cal O}_{2}$. On a $\mu({\cal O}_{1})\cup\mu({\cal O}_{2})=\mu_{1}\cup \mu_{2}=d(\lambda)$. Appliquons la formule 3.6(2) au couple $({\cal O}_{1},{\cal O}_{2})$. Notons plut\^ot $(\mu'_{1},\eta'_{1};\mu'_{2},\eta'_{2})$ les termes indexant la somme de cette formule. Je dis que, si $(\mu'_{1},\mu'_{2})\not=(\mu_{1},\mu_{2})$, le terme 
$$m_{\delta,\sharp}(\mu'_{1},\eta'_{1};\mu'_{2},\eta'_{2})\chi_{\mu'_{1},\eta'_{1}}^{\natural}(Y_{1})\chi^{\natural}_{\mu'_{2},\eta'_{2},\sharp}(Y_{2})$$
est nul. En effet, la non-nullit\'e des deux derniers termes entra\^{\i}ne comme plus les in\'egalit\'es $\mu({\cal O}_{1})\leq \mu'_{1}$ et $\mu({\cal O}_{2})\leq \mu'_{2}$, c'est-\`a-dire $\mu_{1}\leq \mu'_{1}$ et $\mu_{2}\leq \mu'_{2}$. La non-nullit\'e du premier terme entra\^{\i}ne $\mu'_{1}\cup \mu'_{2}\leq d(\lambda)$ d'apr\`es (1). Puisque $\mu_{1}\cup \mu_{2}=d(\lambda)$, les in\'egalit\'es pr\'ec\'edentes sont forc\'ement des \'egalit\'es. Donc $\mu'_{1}=\mu_{1}$ et $\mu'_{2}=\mu_{2}$, contrairement \`a l'hypoth\`ese. Dans la somme de 3.6(2) ne restent donc que les $(\mu'_{1},\eta'_{1};\mu'_{2},\eta'_{2})$ pour lesquels $\mu'_{1}=\mu_{1}$ et $\mu'_{2}=\mu_{2}
$. C'est-\`a-dire que l'on obtient
$$\Theta_{\delta}(f_{{\cal O}_{1},{\cal O}_{2}})=q^{d_{{\cal O}_{1},{\cal O}_{2}}}\Psi(Y_{1},Y_{2}),$$
o\`u $(Y_{1},Y_{2})\in {\cal O}_{1}\times {\cal O}_{2}$. 
D'o\`u $\Theta_{\delta}(f_{{\cal O}_{1},{\cal O}_{2}})\not=0$. Alors $({\cal O}_{1},{\cal O}_{2})$ v\'erifie 3.5(2). 

 \bigskip

\subsection{Une  description de $Res(\delta)$}
On a fix\'e $(\lambda,s,\epsilon)\in \mathfrak{Irr}_{tunip}$ en 3.5. Nous supposons d\'esormais que c'est un \'el\'ement de $\mathfrak{Irr}_{unip-quad}$, ce qui nous suffit d'apr\`es 3.4. On a montr\'e en \cite{W3} 2.2 que cet ensemble s'identifiait \`a $\boldsymbol{{\cal P}}_{2}^{\boldsymbol{symp}}(2n)$. On note $(\lambda^+,\epsilon^+,\lambda^-,\epsilon^-)$ l'\'el\'ement de cet ensemble auquel s'identifie $(\lambda,s,\epsilon)$. On a $\lambda=\lambda^+\cup \lambda^-$ et cette d\'ecomposition est d\'etermin\'ee par l'\'el\'ement $s$. 
 
 Notons $\mathfrak{n}(\lambda^+,\lambda^-)$ l'ensemble des couples de partitions symplectiques $(\nu^+,\nu^-)$ v\'erifiant les conditions suivantes:

(1)(a) $\nu^+\cup\nu^-=\lambda$;

(1)(b) pour tout entier $i\in {\mathbb N}$ impair, $mult_{\nu^-}(i)=0$;

(1)(c) pour tout $i\in Jord_{bp}(\lambda^+)\cap Jord_{bp}(\lambda^-)$, $mult_{\nu^-}(i)\leq 2$;

(1)(d) pour tout $i\in Jord_{bp}(\lambda)$ tel que      $mult_{\lambda^+}(i)=0$ ou $mult_{\lambda^-}(i)=0$, $mult_{\nu^-}(i)\leq 1$. 

Notons $\mathfrak{N}(\lambda_{1},\lambda_{2})$ l'ensemble des quadruplets $(\nu^+,\xi^+,\nu^-,\xi^-)\in \mathfrak{Irr}_{unip-quad}$ tels que $(\nu^+,\nu^-)\in \mathfrak{n}(\lambda^+,\lambda^-)$. 
Fixons un tel quadruplet $(\nu^+,\xi^+,\nu^-,\xi^-)$. Soit $i\in Jord_{bp}(\lambda)$. On d\'efinit un nombre complexe $e(i)$ par les formules suivantes, o\`u on pose par convention $\xi^+(i)=1$ si $i\not\in Jord_{bp}(\nu^+)$ et  $\xi^-(i)=1$ si $i\not\in Jord_{bp}(\nu^-)$:

(2)(a) si $mult_{\nu^-}(i)=0$, $e(i)=\xi^+(i)^{mult_{\lambda^-}(i)}$;

(2)(b) si $mult_{\nu^-}(i)=2$, $e(i)=\epsilon^+(i)\epsilon^-(i)\xi^-(i)\xi^+(i)^{mult_{\lambda^-}(i)-1}$;

(2)(c) si $mult_{\nu^-}(i)=1$ et $mult_{\lambda^-}(i)=0$, $e(i)=\epsilon^+(i)$;

(2)(d) si $mult_{\nu^-}(i)=1$ et $mult_{\lambda^+}(i)=0$, $e(i)=\epsilon^-(i)\xi^-(i)\xi^+(i)^{mult_{\lambda^-}(i)-1}$;

(2)(e) si $mult_{\nu^-}(i)=1$, $mult_{\lambda^+}(i)\geq1$ et $mult_{\lambda^-}(i)\geq1$, $e(i)=\epsilon^-(i)\xi^-(i)\xi^+(i)^{mult_{\lambda^-}(i)-1}+\epsilon^+(i)\xi^+(i)^{mult_{\lambda^-}(i)}$.

On pose
$$e(\lambda^+,\epsilon^+,\lambda^-,\epsilon^-;\nu^+,\xi^+,\nu^-,\xi^-)=2^{-\vert Jord_{bp}(\lambda^+)\vert -\vert Jord_{bp}(\lambda^-)\vert}\prod_{i\in Jord_{bp}(\lambda)}e(i).$$

D'autre part, en \cite{W3} 1.11, on a associ\'e \`a $(\nu^+,\xi^+,\nu^-,\xi^-)$ un \'el\'ement de ${\cal R}$ que l'on a not\'e $j(\boldsymbol{\rho}_{\nu^+,\xi^+}\otimes \boldsymbol{\rho}_{\nu^-,\xi^-})$. Le terme $j$ \'etait un isomorphisme entre ${\cal R}$ et un autre espace. Distinguer ces deux espaces nous \'etait alors utile. Ce ne l'est plus, on identifie l'espace en question \`a ${\cal R}$ gr\^ace \`a l'isomorphisme $j$ et on fait dispara\^{\i}tre ce $j$ de la notation. 

\ass{Lemme}{On a l'\'egalit\'e
$$\kappa_{\delta }=\sum_{(\nu^+,\xi^+,\nu^-,\xi^-)\in \mathfrak{N}(\lambda^+,\lambda^-)} e(\lambda^+,\epsilon^+,\lambda^-,\epsilon^-;\nu^+,\xi^+,\nu^-,\xi^-)\rho\iota(\boldsymbol{\rho}_{\nu^+,\xi^+}\otimes\boldsymbol{\rho}_{\nu^-,\xi^-}).$$}

Preuve. Rappelons quelques  notations de \cite{W3} 1.3. On fixe  un homomorphisme  $\rho_{\lambda}:SL(2;{\mathbb C}) \to Sp(2n;{\mathbb C})$  param\'etr\'e par $\lambda$ et on note $Z(\lambda)$ le commutant dans $Sp(2n;{\mathbb C})$ de son image. Le terme $s$ appartient \`a $Z(\lambda)$ et v\'erifie $s^2=1$. On note $Z(\lambda,s)$ le commutant de $s$ dans $Z(\lambda)$, ${\bf Z}(\lambda,s)$ son groupe des composantes connexes et ${\bf Z}(\lambda,s)^{\vee}$ le groupe des caract\`eres de ${\bf Z}(\lambda,s)$. On a $\epsilon\in {\bf Z}(\lambda,s)^{\vee}$. 
Consid\'erons un sous-ensemble $H\subset Z(\lambda,s)$ qui s'envoie bijectivement sur ${\bf Z}(\lambda,s)$ et est form\'e d'\'el\'ements $h$ v\'erfiant $h^2=1$. Pour $h\in H$, le triplet $(\lambda,s,h)$ appartient \`a l'ensemble $\mathfrak{Endo}_{unip-quad}$ de \cite{W3} 2.2. Il lui est associ\'e une repr\'esentation virtuelle
$$(3) \qquad \Pi(\lambda,s,h)=\sum_{\epsilon'\in {\bf Z}(\lambda,s)^{\vee}}\pi(\lambda,s,\epsilon')\epsilon'(h).$$
Par inversion de Fourier dans le groupe ${\bf Z}(\lambda,s)$, on a
$$\pi(\lambda,s,\epsilon)=\vert {\bf Z}(\lambda,s)\vert ^{-1}\sum_{h\in H}\Pi(\lambda,s,h)\epsilon(h).$$
On applique $Res\circ D$ \`a cette \'egalit\'e:
$$ Res(\delta)=\vert {\bf Z}(\lambda,s)\vert ^{-1}\sum_{h\in H}Res\circ D\circ \Pi(\lambda,s,h)\epsilon(h).$$
Utilisons  l'involution ${\cal F}^{par}$ de \cite{W3} 1.9. Puisque, justement, c'est une involution, on peut composer \`a gauche le membre de droite ci-dessus par ${\cal F}^{par}\circ {\cal F}^{par}$. Le th\'eor\`eme 2.7 de \cite{W3}
n'est plus conditionnel puisqu'on a d\'emontr\'e en \cite{W4} le th\'eor\`eme 2.1 de \cite{W3}. Il nous dit que ${\cal F}^{par}\circ Res\circ D\circ\Pi(\lambda,s,h)=Res\circ D\circ\Pi(\lambda,h,s)$. D'o\`u
$$  Res(\delta)=\vert {\bf Z}(\lambda,s)\vert ^{-1}\sum_{h\in H}{\cal F}^{par }\circ Res\circ D\circ \Pi(\lambda,h,s)\epsilon(h).$$
On peut encore d\'evelopper le membre de droite en utilisant (3) o\`u l'on \'echange $s$ et $h$:
$$(4) \qquad   Res(\delta)=\vert {\bf Z}(\lambda,s)\vert ^{-1}\sum_{h\in H}\sum_{\xi\in {\bf Z}(\lambda,h)^{\vee}}{\cal F}^{par}\circ Res\circ D(\pi(\lambda,h,\xi))\xi(s)\epsilon(h).$$
L'ensemble $\mathfrak{Endo}_{unip-quad}$ s'identifie \`a ${\cal P}^{symp}_{4}(2n)$. Plus pr\'ecis\'ement, l'ensemble des \'el\'ements de $\mathfrak{Endo}_{unip-quad}$ de la forme $(\lambda,s,h)$ (c'est-\`a-dire dont les deux premiers termes sont nos \'el\'ements fix\'es $\lambda$ et $s$) s'identifie aux quadruplets $(\lambda^{++},\lambda^{-+},\lambda^{+-},\lambda^{--})\in {\cal P}^{symp}_{4}(2n)$ tels que $\lambda^{++}\cup \lambda^{+-}=\lambda^+$ et $\lambda^{-+}\cup \lambda^{--}=\lambda^{-}$. Si $(\lambda,s,h)$ correspond ainsi \`a $(\lambda^{++},\lambda^{-+},\lambda^{+-},\lambda^{--})$, l'image de $h$ dans 
$${\bf Z}(\lambda,s)\simeq ({\mathbb Z}/2{\mathbb Z})^{Jord_{bp}(\lambda^+)}\times ({\mathbb Z}/2{\mathbb Z})^{Jord_{bp}(\lambda^-)}$$ est
$$(5) \qquad (mult_{\lambda^{+-}(i)})_{i\in Jord_{bp}(\lambda^+)}\times (mult_{\lambda^{--}(i)})_{i\in Jord_{bp}(\lambda^-)},$$
o\`u il s'agit en fait des images des multiplicit\'es dans ${\mathbb Z}/2{\mathbb Z}$. 
On voit que l'on peut choisir  $H$ de sorte que l'ensemble des $(\lambda,s,h)$ pour $h\in H$ s'identifie \`a l'ensemble des quadruplets satisfaisant les conditions ci-dessus et de plus: $mult_{\lambda^{+-}}(i)\leq 1$ et $mult_{\lambda^{--}}(i)\leq 1$ pour tout $i$. On peut \'evidemment renforcer des in\'egalit\'es en $mult_{\lambda^{+-}}(i)\leq inf( 1,mult_{\lambda^+}(i))$ et $mult_{\lambda^{--}}(i)\leq inf( 1,mult_{\lambda^-}(i))$ (par exemple, $mult_{\lambda^{+-}}(i)\leq mult_{\lambda^+}(i)$ puisque $\lambda^+=\lambda^{++}\cup \lambda^{+-}$). 
 On choisit ainsi l'ensemble $H$. Pour $h\in H$, continuons \`a noter $(\lambda^{++},\lambda^{-+},\lambda^{+-},\lambda^{--})$ le quadruplet associ\'e \`a $(\lambda,s,h)$ et posons $\nu^+=\lambda^{++}\cup \lambda^{-+}$, $\nu^-=\lambda^{+-}\cup \lambda^{--}$. Le couple $(\nu^+,\nu^-)$  appartient \`a notre ensemble $\mathfrak{n}(\lambda^+,\lambda^-)$.  Le groupe ${\bf Z}(\lambda,h)$ s'identifie \`a 
$$({\mathbb Z}/2{\mathbb Z})^{Jord_{bp}(\nu^+)}\times ({\mathbb Z}/2{\mathbb Z})^{Jord_{bp}(\nu^-)}$$
et un \'el\'ement $\xi\in {\bf Z}(\lambda,h)^{\vee}$ s'identifie \`a un couple $(\xi^+,\xi^-)$. Le quadruplet $(\nu^+,\xi^+,\nu^-,\xi^-)$ appartient \`a $ \mathfrak{N}(\lambda^+,\lambda^-)$. Mais l'application $h\mapsto (\nu^+,\nu^-)$ n'est pas injective ($h$ est une classe de conjugaison par $Z(\lambda,s)$ et $(\nu^+,\nu^-)$ param\`etre sa classe de conjugaison par $Z(\lambda)$). L'\'egalit\'e (4) se r\'ecrit
$$ Res(\delta)=\sum_{(\nu^+,\xi^+,\nu^-,\xi^-)\in \mathfrak{N}(\lambda^+,\lambda^-)}f(\nu^+,\xi^+,\nu^-,\xi^-){\cal F}^{par}\circ Res\circ D(\pi(\nu^+,\xi^+,\nu^-,\xi^-)),$$
o\`u:

pour $(\lambda,h,\xi)$ correspondant \`a $(\nu^+,\xi^+,\nu^-,\xi^-)$, on a not\'e $\pi(\nu^+,\xi^+,\nu^-,\xi^-)=\pi(\lambda,h,\xi)$;

$f(\nu^+,\xi^+,\nu^-,\xi^-)$ est la somme des $\vert {\bf Z}(\lambda,s)\vert ^{-1}\epsilon(h)\xi(s)$ sur les $h\in H$ d'image $(\nu^+,\nu^-)$.

\noindent Pour $(\nu^+,\xi^+,\nu^-,\xi^-)\in \mathfrak{N}(\lambda^+,\lambda^-)$, on a 
$$Res\circ D(\pi(\nu^+,\xi^+,\nu^-,\xi^-))=Rep\circ\rho\iota(\boldsymbol{\rho}_{\nu^+,\xi^+}\otimes \boldsymbol{\rho}_{\nu^-,\xi^-})$$
d'apr\`es la proposition 1.11 de \cite{W3}. On a aussi ${\cal F}^{par}\circ Rep=k$, cf. \cite{W3} 1.9. Donc
$$ Res(\delta)=\sum_{(\nu^+,\xi^+,\nu^-,\xi^-)\in \mathfrak{N}(\lambda^+,\lambda^-)}f(\nu^+,\xi^+,\nu^-,\xi^-) k\circ \rho\iota(\boldsymbol{\rho}_{\nu^+,\xi^+}\otimes \boldsymbol{\rho}_{\nu^-,\xi^-}).$$
Puisque $\kappa_{\delta}$ est l'\'el\'ement de ${\cal R}$ tel que $Res(\delta)=k(\kappa_{\delta})$, on obtient la formule de l'\'enonc\'e, \`a condition de prouver l'\'egalit\'e:

 $f(\nu^+,\xi^+,\nu^-,\xi^-)=e(\lambda^+,\epsilon^+,\lambda^-,\epsilon^-;\nu^+,\xi^+,\nu^-,\xi^-)$ pour tout $(\nu^+,\xi^+,\nu^-,\xi^-)\in \mathfrak{N}(\lambda^+,\lambda^-)$.

Fixons donc $(\nu^+,\xi^+,\nu^-,\xi^-)\in \mathfrak{N}(\lambda^+,\lambda^-)$. On compare tout de suite le facteur $\vert {\bf Z}(\lambda,s)\vert ^{-1}$ figurant dans la d\'efinition de $f(\nu^+,\xi^+,\nu^-,\xi^-)$ avec le facteur $2^{-\vert Jord_{bp}(\lambda^+)\vert -\vert Jord_{bp}(\lambda^-)\vert}$ figurant dans celle de $e(\lambda^+,\epsilon^+,\lambda^-,\epsilon^-;\nu^+,\xi^+,\nu^-,\xi^-)$: ils sont \'egaux. On note $f$ le terme $f(\nu^+,\xi^+,\nu^-,\xi^-)$ priv\'e de ce facteur et on doit prouver que $f=\prod_{i\in Jord_{bp}(\lambda)}e(i)$, avec les notations pr\'ec\'edant l'\'enonc\'e. 
Soit $h\in H$ d'image $(\nu^+,\nu^-)$. Notons encore 
 $(\lambda^{++},\lambda^{-+},\lambda^{+-},\lambda^{--})$ le quadruplet associ\'e \`a $(\lambda,s,h)$.  D'apr\`es  (5), on a
$$\epsilon(h)=(\prod_{i\in Jord_{bp}(\lambda^+)}\epsilon^+(i)^{mult_{\lambda^{+-}}(i)})(\prod_{i\in Jord_{bp}(\lambda^-)}\epsilon^-(i)^{mult_{\lambda^{--}}(i)}).$$
On simplifie cette \'egalit\'e en
$$\epsilon(h)=\prod_{i\in Jord_{bp}(\lambda)}\epsilon^+(i)^{mult_{\lambda^{+-}}(i)}\epsilon^-(i)^{mult_{\lambda^{--}}(i)},$$
avec la convention $\epsilon^+(i)=1$ si $i\not\in Jord_{bp}(\lambda^+)$ et $\epsilon^-(i)=1$ pour $i\not\in Jord_{bp}(\lambda^-)$. 
Le quadruplet associ\'e \`a $(\lambda,h,s)$ est $(\lambda^{++},\lambda^{+-},\lambda^{-+},\lambda^{--})$. Le couple $(\xi^+,\xi^-)$ s'identifie \`a un \'el\'ement de ${\bf Z}(\lambda,h)^{\vee}$ et on a la formule similaire (avec une convention analogue):
 $$\xi(s)=\prod_{i\in Jord_{bp}( \lambda)}\xi^+(i)^{mult_{\lambda^{-+}}(i)} \xi^-(i)^{mult_{\lambda^{--}}(i)}.$$
 Posons simplement $m^+(i)=mult_{\lambda^{+-}}(i)$ et $m^-(i)=mult_{\lambda^{--}}(i)$. On a $mult_{\lambda^{-+}}(i)=mult_{\lambda^-}(i)-m^-(i)$ et on obtient
 $$\epsilon(h)\xi(s)=\prod_{i\in Jord_{bp}( \lambda)}\epsilon^+(i)^{m^+(i)}\epsilon^-(i)^{m^-(i)}\xi^+(i)^{mult_{\lambda^-}(i)-m^-(i)}\xi^-(i)^{m^-(i)}.$$
 Pour $i\in Jord_{bp}(\lambda)$, notons $E(i)$ l'ensemble des couples $(m^+,m^-)\in \{0,1\}^2$ tels que
 
 $m^+\leq inf(1,mult_{\lambda^+}(i))$ et  $m^-\leq inf(1,mult_{\lambda^-}(i))$;
 
 $m^++m^-=mult_{\nu^-}(i)$.
 
 L'application 
 $$h\mapsto ( m^+(i),m^-(i))_{i\in Jord_{bp}(\lambda)}$$
 identifie l'ensemble des $h\in H$ d'image $(\nu^+,\nu^-)$ avec $\prod_{i\in Jord_{bp}(\lambda)}E(i)$.  On voit alors que $f=\prod_{i\in Jord_{bp}(\lambda)}f_{i}$, o\`u, pour $i\in Jord_{bp}(\lambda)$, on a pos\'e
 $$f_{i}=\sum_{(m^+,m^-)\in E(i)}\epsilon^+(i)^{m^+}\epsilon^-(i)^{m^-}\xi^+(i)^{mult_{\lambda^-}(i)-m^-}\xi^-(i)^{m^-}.$$
 Il reste \`a d\'emontrer l'\'egalit\'e $f_{i}=e_{i}$ pour tout $i\in Jord_{bp}(\lambda)$. C'est un calcul \'el\'ementaire que l'on effectue en distinguant chacun des cas (2)(a) \`a (2)(e). On le laisse au lecteur. $\square$

\bigskip

\subsection{Preuve de 3.7(1) }
On fixe $n_{1},n_{2}\in {\mathbb N}$ avec $n_{1}+n_{2}=n$ et $n_{2}\geq1$ si $\sharp=an$. On fixe $(\mu_{1},\eta_{1};\mu_{2},\eta_{2})\in \boldsymbol{{\cal P}^{orth}}(2n_{1}+1;k=1)\times \boldsymbol{{\cal P}^{orth}}(2n_{2};k=0)$ et on suppose $m_{\delta,\sharp}(\mu_{1},\eta_{1};\mu_{2},\eta_{2})\not=0$.   D'apr\`es la d\'efinition de $m_{\delta,\sharp}(\mu_{1},\eta_{1};\mu_{2},\eta_{2})\not=0$, on peut fixer $\zeta'=\pm$ tel que $m_{\delta}(\rho_{\mu_{1},\eta_{1}}\otimes sgn,\rho_{\mu_{2},\eta_{2}}^{\zeta'}\otimes sgn)\not=0$.

On a d\'efini en 1.4 et 1.5 les partitions sp\'eciales $sp(\mu_{1},\eta_{1})$ et $sp(\mu_{2},\eta_{2})$. D'apr\`es les r\'esultats de ces paragraphes, on a 
$$(1) \qquad \mu_{1}\leq sp(\mu_{1},\eta_{1})\,\,\mu_{2}\leq sp(\mu_{2},\eta_{2}).$$
Le symbole associ\'e \`a $\rho_{\mu_{1},\eta_{1}}$ appartient \`a la famille de $sp(\mu_{1},\eta_{1})$. Posons $\lambda_{1}=d(sp(\mu_{1},\eta_{1}))$ et $\rho_{1}=\rho_{\mu_{1},\eta_{1}}\otimes sgn$. D'apr\`es 1.6, le symbole associ\'e \`a $\rho_{1}$ appartient \`a la famille de $\lambda_{1}$. Posons $\lambda_{2}=d(sp(\mu_{2},\eta_{2}))$ et $\rho_{2}=\rho_{\mu_{2},\eta_{2}}\otimes sgn$ si $\mu_{2}$ n'est pas exceptionnel.  Si $\mu_{2}$ est exceptionnel, on rel\`eve $\mu_{2}$ en un \'el\'ement $\underline{\mu}_{2}$ de $\underline{{\cal P}}^{orth}(2n_{2})$ et on pose $\rho_{2}=\rho_{\underline{\mu}_{2},\eta_{2}}\otimes sgn$. On a de m\^eme: le symbole associ\'e \`a $\rho_{2}$ appartient \`a la famille de $\lambda_{2}$. La repr\'esentation $\rho_{\mu_{2},\eta_{2}}^{\zeta'}\otimes sgn$ est l'un des prolongements de $\rho_{2}$ \`a $W_{n_{2}}$ donc est de la forme $\rho_{2}^{\zeta}$ pour un $\zeta=\pm$ (on n'a pas en g\'en\'eral $\zeta=\zeta'$ mais peu importe). On a donc $m_{\delta}(\rho_{1},\rho_{2}^{\zeta})\not=0$. Rappelons que ce terme est la "multiplicit\'e" de $\rho_{1}\otimes \rho_{2}^{\zeta}$ dans $\kappa_{\delta,0}$. La fonction $\kappa_{\delta}$ est calcul\'ee par le lemme 3.8 (rappelons que l'on suppose $(\lambda,s,\epsilon)\in \mathfrak{Irr}_{unip-quad}$). La fonction $\kappa_{\delta,0}$ est calcul\'ee par une formule analogue, o\`u l'on se restreint aux $(\nu^+,\xi^+,\nu^-,\xi^-)\in \mathfrak{N}(\lambda^+,\lambda^-)$ tels que $k_{\nu^+,\xi^+}=k_{\nu^-,\xi^-}=0$.  Notons ce sous-ensemble $\mathfrak{N}_{0}(\lambda^+,\lambda^-)$. D'apr\`es ce lemme, on peut fixer un \'el\'ement  $(\nu^+,\xi^+,\nu^-,\xi^-)\in \mathfrak{N}_{0}(\lambda^+,\lambda^-)$ tel que la multiplicit\'e de $\rho_{1}\otimes \rho_{2}^{\zeta}$ dans $\rho\iota(\boldsymbol{\rho}_{\nu^+,\xi^+}\otimes \boldsymbol{\rho}_{\nu^-,\xi^-})$ est non nulle.  
On sait que la repr\'esentation $\boldsymbol{\rho}_{\nu^+,\xi^+}$ est de la forme 
$$\boldsymbol{\rho}_{\nu^+,\xi^+}=\sum_{\dot{\nu}^+,\dot{\xi}^+}c(\nu^+,\xi^+;\dot{\nu}^+,\dot{\xi}^+)\rho_{\dot{\nu}^+,\dot{\xi}^+},$$
o\`u $(\dot{\nu}^+,\dot{\xi}^+)$ parcourt les \'el\'ements de $\boldsymbol{{\cal P}^{symp}}(S(\nu^+))$ tels que $k_{\dot{\nu}^+,\dot{\xi}^+}=0$. On note  

\noindent $\boldsymbol{{\cal P}^{symp}}(S(\nu^+);k=0)$ l'ensemble de ces \'el\'ements. 
Le coefficient $c(\nu^+,\xi^+;\dot{\nu}^+,\dot{\xi}^+)$ n'est non nul que si $\nu^+\leq \dot{\nu}^+$. De m\^emes propri\'et\'es valent pour $\boldsymbol{\rho}_{\nu^-,\xi^-}$. On peut donc fixer $(\dot{\nu}^+,\dot{\xi}^+) \in\boldsymbol{{\cal P}^{symp}}(S(\nu^+);k=0)$ et $(\dot{\nu}^-,\dot{\xi}^-) \in\boldsymbol{{\cal P}^{symp}}(S(\nu^-);k=0)$ tels que  

(2) $\nu^+\leq \dot{\nu}^+,\,\,\nu^-\leq \dot{\nu}^-$,

 \noindent et la multiplicit\'e de  $\rho_{1}\otimes \rho_{2}^{\zeta}$ dans $\rho_{\dot{\nu}^+,\dot{\xi}^+}\otimes \rho_{\dot{\nu}^-,\dot{\xi}^-}$ soit non nulle. Cette multiplicit\'e est exactement le terme $m(\rho_{1},\rho_{2}^{\zeta};\rho_{\dot{\nu}^+,\dot{\xi}^+},\rho_{\dot{\nu}^-,\dot{\xi}^-})$ d\'efini en 1.12. D'apr\`es la proposition de ce paragraphe, la non-nullit\'e de cette multiplicit\'e entra\^{\i}ne

(3) $\dot{\nu}^+\cup \dot{\nu}^-\leq ind(\lambda_{1},\lambda_{2})$.

Par d\'efinition de $\mathfrak{N}(\lambda^+,\lambda^-)$, on a $\nu^+\cup\nu^-=\lambda$. Les in\'egalit\'es (2) et (3) entra\^{\i}nent
$\lambda\leq ind(\lambda_{1},\lambda_{2})$ d'o\`u aussi $d(ind(\lambda_{1},\lambda_{2}))\leq d(\lambda)$ puisque la dualit\'e est d\'ecroissante.  On applique la proposition 1.9: $d(\lambda_{1})\cup d(\lambda_{2})\leq d(ind(\lambda_{1},\lambda_{2}))$, d'o\`u aussi $d(\lambda_{1})\cup d(\lambda_{2})\leq d(\lambda)$. Or $d(\lambda_{1})=sp(\mu_{1},\eta_{1})$ et $d(\lambda_{2})=sp(\mu_{2},\eta_{2})$ par d\'efinition. En utilisant les in\'egalit\'es (1), on en d\'eduit $\mu_{1}\cup \mu_{2}\leq d(\lambda)$, ce qui d\'emontre 3.7(1).

\bigskip

\subsection{Preuve de 3.7(2)}
On a suppos\'e $(\lambda,s,\epsilon)\in \mathfrak{Irr}_{unip-quad}$. Maintenant, on suppose de plus que les termes de $\lambda$ sont tous pairs. C'est loisible d'apr\`es 3.4.

On  fixe une fonction $\tau:Jord_{bp}(\lambda)\to {\mathbb Z}/2{\mathbb Z}$  v\'erifiant les conditions suivantes:

(1)(a) pour $i\in Jord_{bp}(\lambda)$ tel que $mult_{\lambda}(i)=1$, $\tau(i)=0$;

(1)(b) pour $i\in Jord_{bp}(\lambda^+)\cap Jord_{bp}(\lambda^-)$, $(-1)^{\tau(i)}=\epsilon^+(i)\epsilon^-(i)$.

Remarquons que les deux cas sont exclusifs: dans le cas (b), on a $mult_{\lambda}(i)\geq2$.

Appliquant la proposition 1.11, on introduit des entiers $n_{1},n_{2}\in {\mathbb N}$ tels que $n_{1}+n_{2}=n$ et des partitions $\lambda_{1}\in {\cal P}^{symp,sp}(2n_{1})$, $\lambda_{2}\in {\cal P}^{orth,sp}(2n_{2})$ telles que $\lambda_{1}$ et $\lambda_{2}$ induisent r\'eguli\`erement $\lambda$, $d(\lambda_{1})\cup d(\lambda_{2})=d(\lambda)$ et $\tau_{\lambda_{1},\lambda_{2}}=\tau$. On fixe des couples $(\tau_{1},\delta_{1})$ et $(\tau_{2},\delta_{2})$  param\'etrant des symboles $(X_{1},Y_{1})$ dans la famille de $\lambda_{1}$ et $(X_{2},,Y_{2})$ dans la famille de $\lambda_{2}$, avec $\delta_{1}=0$ et $\delta_{2}=0$. On pose $\mu_{1}=d(\lambda_{1})$, $\mu_{2}=d(\lambda_{2})$. Ce sont des partitions sp\'eciales. Le symbole $d(X_{1},Y_{1})$ appartient \`a la famille de $\mu_{1}$ et est param\'etr\'e par un couple $(\tau'_{1},\delta'_{1})$ tel que $\delta'_{1}=0$. D'apr\`es le lemme 1.4, c'est le symbole de la repr\'esentation $\rho_{\mu_{1},\eta_{1}}$ pour un couple $(\mu_{1},\eta_{1})\in \boldsymbol{{\cal P}^{orth}}(2n_{1}+1;k=1)$. Le symbole $d(X_{2},Y_{2})$ appartient \`a la famille de $\mu_{2}$ et est param\'etr\'e par un couple $(\tau'_{2},\delta'_{2})$ tel que $\delta'_{2}=0$. Supposons $\mu_{2}$ non exceptionnel. D'apr\`es le lemme 1.5, $d(X_{2},Y_{2})$ est le symbole de la repr\'esentation $\rho_{\mu_{2},\eta_{2}}$ pour un couple $(\mu_{2},\eta_{2})\in \boldsymbol{{\cal P}^{orth}}(2n_{2};k=0)$. Dans le cas o\`u $\mu_{2}$ est exceptionnel, on a de m\^eme un couple $(\mu_{2},\eta_{2})$ mais la repr\'esentation $\rho_{\mu_{2},\eta_{2}}$ doit \^etre remplac\'ee par $\rho_{\underline{\mu}_{2},\eta_{2}}$, o\`u $\underline{\mu}_{2}$ est l'un des rel\`evements de $\mu_{2}$ dans $\underline{{\cal P}}^{orth}(2n_{2})$. On va montrer que le quadruplet $(\mu_{1},\eta_{1};\mu_{2},\eta_{2})$
 v\'erifie la condition  3.7(2). 
 
 Tout d'abord, on a $\mu_{1}\cup \mu_{2}=d(\lambda_{1})\cup d(\lambda_{2})=d(\lambda)$.
 
 On veut prouver que $m_{\delta,\sharp}(\mu_{1},\eta_{1};\mu_{2},\eta_{2})\not=0$. Posons $sgn_{\sharp}=1$ si $\sharp=iso$, $sgn_{\sharp}=-1$ si $\sharp=an$. Le terme $m_{\delta,\sharp}(\mu_{1},\eta_{1};\mu_{2},\eta_{2})$ est \'egal  \`a
 $$(2) \qquad m_{\delta}(\rho_{\mu_{1},\eta_{1}}\otimes sgn,\rho^+_{\mu_{2},\eta_{2}}\otimes sgn)+sgn_{\sharp}m_{\delta}(\rho_{\mu_{1},\eta_{1}}\otimes sgn,\rho^-_{\mu_{2},\eta_{2}}\otimes sgn),$$
 \'eventuellement divis\'e par $1/2$. La repr\'esentation $\rho_{\mu_{1},\eta_{1}}\otimes sgn$ n'est autre que la repr\'esentation $\rho_{1}$ de $W_{n_{1}}$ dont le symbole est $(X_{1},Y_{1})$. Les repr\'esentations 
$\rho^+_{\mu_{2},\eta_{2}}\otimes sgn$ et $\rho^-_{\mu_{2},\eta_{2}}\otimes sgn$ sont les prolongements (\'eventuellement \'egaux) \`a $W_{n_{2}}$ d'une repr\'esentation $\rho_{2}\in W_{n_{2}}^D$ dont le symbole est $(X_{2},Y_{2})$. Le terme (2) est \'egal, au signe pr\`es, \`a
$$(3) \qquad m_{\delta}(\rho_{1},\rho_{2}^+)+s_{\sharp}m_{\delta}(\rho_{1},\rho_{2}^-).$$
Soit $\zeta=\pm$. En reprenant la preuve du paragraphe pr\'ec\'edent, on calcule
$$m_{\delta}(\rho_{1},\rho_{2}^{\zeta})=\sum_{(\nu^+,\xi^+,\nu^-,\xi^-)\in {\cal N}_{0}(\lambda^+,\lambda^-)}e(\lambda^+,\epsilon^+,\lambda^-,\epsilon^-;\nu^+,\xi^+,\nu^-,\xi^-)$$
$$\sum_{(\dot{\nu}^+,\dot{\xi}^+)\in \boldsymbol{{\cal P}^{symp}}(S(\nu^+);k=0)}\sum_{(\dot{\nu}^-,\dot{\xi}^-)\in \boldsymbol{{\cal P}^{symp}}(S(\nu^-);k=0)}c(\nu^+,\xi^+;\dot{\nu}^+,\dot{\xi}^+)c(\nu^-,\xi^-;\dot{\nu}^-,\dot{\xi}^-)$$
$$m(\rho_{1},\rho_{2}^{\zeta};\rho_{\dot{\nu}^+,\dot{\xi}^+},\rho_{\dot{\nu}^-,\dot{\xi}^-}).$$
Comme dans le paragraphe pr\'ec\'edent, la non-nullit\'e du terme que l'on somme entra\^{\i}ne les in\'egalit\'es  (2) et (3) de ce paragraphe. On a $\nu^+\cup\nu^-=\lambda$ et, ici, on sait par hypoth\`ese que $ind(\lambda_{1},\lambda_{2})=\lambda$. Ces in\'egalit\'es (2) et (3) sont donc des \'egalit\'es.  Maintenant que $\nu^+=\dot{\nu}^+$, on sait que  la relation $c(\nu^+,\xi^+;\dot{\nu}^+,\dot{\xi}^+)\not=0$ \'equivaut \`a $\dot{\xi}^+=\xi^+$ et que, si  elle est v\'erifi\'ee, on a $c(\nu^+,\xi^+;\dot{\nu}^+,\dot{\xi}^+)=1$. De m\^eme bien s\^ur pour les objets associ\'es \`a $\nu^-$. Cela nous d\'ebarrasse des sommes en $(\dot{\nu}^+,\dot{\xi}^+)$ et $(\dot{\nu}^-,\dot{\xi}^-)$ et des coefficients $c(\nu^+,\xi^+;\dot{\nu}^+,\dot{\xi}^+)$ et $c(\nu^-,\xi^-;\dot{\nu}^-,\dot{\xi}^-)$. On a simplement
 $$m_{\delta}(\rho_{1},\rho_{2}^{\zeta})=\sum_{(\nu^+,\xi^+,\nu^-,\xi^-)\in {\cal N}_{0}(\lambda^+,\lambda^-)}e(\lambda^+,\epsilon^+,\lambda^-,\epsilon^-;\nu^+,\xi^+,\nu^-,\xi^-)$$
$$m(\rho_{1},\rho_{2}^{\zeta};\rho_{\nu^+,\xi^+},\rho_{\nu^-,\xi^-}).$$ 
Les \'el\'ements $(\nu^+,\xi^+,\nu^-,\xi^-)\in {\cal N}_{0}(\lambda^+,\lambda^-)$ pour lesquels $m(\rho_{1},\rho_{2}^{\zeta};\rho_{\nu^+,\xi^+},\rho_{\nu^-,\xi^-})\not=0$ sont exactement les \'el\'ements de ${\cal N}(\lambda^+,\lambda^-)$ qui v\'erifient l'hypoth\`ese $(B)^{\zeta}$ de 1.13. D'apr\`es la proposition de ce paragraphe, ce sont aussi les \'el\'ements de ${\cal N}(\lambda^+,\lambda^-)$ qui v\'erifient $(A)^{\zeta}$ et, pour ces \'el\'ements, on a $m(\rho_{1},\rho_{2}^{\zeta};\rho_{\nu^+,\xi^+},\rho_{\nu^-,\xi^-})=1$. La condition $(A)^{\zeta}$ se d\'ecompose en deux: 

(4) pour $i\in Jord_{bp}(\lambda)$, $mult_{\nu^-}(i)\equiv c^{\zeta}(i)$,

\noindent o\`u on a pos\'e $c^{\zeta}(i)=\delta^{-\zeta}(i)-\delta^{-\zeta}(i^+)$;

(5) pour $i\in Jord_{bp}(\nu^+)$, $\xi^+(i)=(-1)^{\tau^{\zeta}(i)}$; pour $i\in Jord_{bp}(\nu^-)$, $\xi^-(i)=(-1)^{\tau^{-\zeta}(i)}$.

  Notons $\mathfrak{n}^{\zeta}(\lambda_{1},\lambda_{2})$ l'ensemble des  $(\nu^+,\nu^-)\in \mathfrak{n}(\lambda_{1},\lambda_{2})$ v\'erifiant la condition (4). Dans la suite du calcul, pour $(\nu^+,\nu^-)\in \mathfrak{n}^{\zeta}(\lambda_{1},\lambda_{2})$, notons $(\xi^+,\xi^-)$ le couple d\'etermin\'e par la condition (5). On obtient
 $$(6) \qquad m_{\delta}(\rho_{1},\rho_{2}^{\zeta})=\sum_{(\nu^+,\nu^-)\in \mathfrak{n}^{\zeta}(\lambda^+,\lambda^-)}e(\lambda^+,\epsilon^+,\lambda^-,\epsilon^-;\nu^+,\xi^+,\nu^-,\xi^-).$$
 
 Soit $i\in Jord_{bp}(\lambda)$. D\'efinissons un ensemble $M^{\zeta}(i)$ par les \'egalit\'es suivantes:
 
  si $c^{\zeta}(i)=1$, $M^{\zeta}(i)=\{1\}$;

si $c^{\zeta}(i)=0$ et $mult_{\lambda^+}(i)mult_{\lambda^-}(i)=0$, $M^{\zeta}(i)=\{0\}$;

 si $c^{\zeta}(i)=0$ et $mult_{\lambda^+}(i)mult_{\lambda^-}(i)\not=0$, $M^{\zeta}(i)=\{0,2\}$.
 
  On v\'erifie \`a l'aide de la d\'efinition de 3.8 et de la relation (4) ci-dessus que l'application $(\nu^+,\nu^-)\mapsto (mult_{\nu^-}(i))_{i\in Jord_{bp}(\lambda)}$ est une bijection de $\mathfrak{n}^{\zeta}(\lambda^+,\lambda^-)$  sur $\prod_{i\in Jord_{bp}(\lambda)}M^{\zeta}_{i}$.
  
  Soit $i\in Jord_{bp}(\lambda)$ et $m\in M^{\zeta}(i)$. D\'efinissons 
 un nombre $e^{\zeta}(i,m)$ par les \'egalit\'es suivantes:
 
 $e^{\zeta}(i,0)=(-1)^{\tau^{\zeta}(i)mult_{\lambda^-}(i)}$;
 
si $mult_{\lambda^+}(i)mult_{\lambda^-}(i)\not=0$,  $e^{\zeta}(i,1)=\epsilon^+(i)(-1)^{\tau^{\zeta}(i)mult_{\lambda^-}(i)}+\epsilon^-(i)(-1)^{\tau^{-\zeta}(i)+\tau^{\zeta}(i)(mult_{\lambda^-}(i)-1)}$;

si $mult_{\lambda^+}(i)=0$,  $e^{\zeta}(i,1)= \epsilon^-(i)(-1)^{\tau^{-\zeta}(i)+\tau^{\zeta}(i)(mult_{\lambda^-}(i)-1)}$;

si $mult_{\lambda^-}(i)=0$, $e^{\zeta}(i,1)=\epsilon^+(i)  $;
 
 $e^{\zeta}(i,2)=\epsilon^+(i)\epsilon^-(i)(-1)^{\tau^{-\zeta}(i)+\tau^{\zeta}(i)(mult_{\lambda^-}(i)-1)}$.

  Pour $(\nu^+,\nu^-)\in \mathfrak{n}^{\zeta}(\lambda^+,\lambda^-)$, on  v\'erifie \`a l'aide de la d\'efinition de 3.8 et de la relation (5) ci-dessus que l'on a l'\'egalit\'e
  $$e(\lambda^+,\epsilon^+,\lambda^-,\epsilon^-;\nu^+,\xi^+,\nu^-,\xi^-) =2^{-\vert Jord_{bp}(\lambda^{+})\vert-\vert Jord_{bp}(\lambda^-)\vert }\prod_{i\in Jord_{bp}(\lambda)}e^{\zeta}(i,mult_{\nu^-}(i)).$$
La relation (6) se transforme en   
  $$(7) \qquad m_{\delta}(\rho_{1},\rho_{2}^{\zeta})=2^{-\vert Jord_{bp}(\lambda^{+})\vert-\vert Jord_{bp}(\lambda^-)\vert }\prod_{i\in Jord_{bp}(\lambda)}S^{\zeta}(i),$$
o\`u on a pos\'e
$$S^{\zeta}(i)=\sum_{m\in M^{\zeta}(i)}e^{\zeta}(i,m).$$
Fixons $i\in Jord_{bp}(\lambda)$ et calculons
$S^{\zeta}(i) $. Remarquons qu'en vertu de la condition (1)(b) ci-dessus,  de l'\'egalit\'e $\tau=\tau_{\lambda_{1},\lambda_{2}}$ et de la relation 1.13(1), on a l'\'egalit\'e

(8) $(-1)^{\tau^+(i)+\tau^-(i)}=\epsilon^+(i)\epsilon^-(i)$ si $mult_{\lambda^+}(i)mult_{\lambda^-}(i)\not=0$;

\noindent La d\'efinition des termes $e^{\zeta}(i,m)$ se simplifie alors dans les deux cas suivants:

si $mult_{\lambda^+}(i)mult_{\lambda^-}(i)\not=0$,  $e^{\zeta}(i,1)=2\epsilon^+(i)(-1)^{\tau^{\zeta}(i)mult_{\lambda^-}(i)}$;

$e^{\zeta}(i,2)=(-1)^{\tau^{\zeta}(i)mult_{\lambda^-}(i)}$.

On calcule alors

(9)(a) si $c^{\zeta}(i)=1$ et $mult_{\lambda^-}(i)=0$, $S^{\zeta}(i)=\epsilon^+(i) $;

(9)(b) si $c^{\zeta}(i)=1$ et $mult_{\lambda^+}(i)=0$, $S^{\zeta}(i)=\epsilon^-(i)(-1)^{\tau^{-\zeta}(i)+\tau^{\zeta}(i)(mult_{\lambda^-}(i)-1)}$;

(9)(c) si $c^{\zeta}(i)=1$ et $mult_{\lambda^+}(i)mult_{\lambda^-}(i)\not=0$, $S^{\zeta}(i)=2\epsilon^+(i)(-1)^{mult_{\lambda^-}(i)\tau^{\zeta}(i)} $;

(9)(d) si $c^{\zeta}(i)=0$ et $mult_{\lambda^+}(i)mult_{\lambda^-}(i)=0$, $S^{\zeta}(i)=(-1)^{mult_{\lambda^-}(i)\tau^{\zeta}(i)}$;

(9)(e) si $c^{\zeta}(i)=0$ et $mult_{\lambda^+}(i)mult_{\lambda^-}(i)\not=0$, $S^{\zeta}(i)=2(-1)^{mult_{\lambda^-}(i)\tau^{\zeta}(i)}$.

On voit  que $S^{\zeta}(i)$ est non nul pour tout $i$. On d\'eduit de (7) que $m_{\delta}(\rho_{1},\rho_{2}^{\zeta})\not=0$. Mais cela ne suffit pas \`a prouver que l'expression (3) est non nulle. Pour cela, montrons que, pour tout $i\in Jord_{bp}(\lambda)$, on a l'\'egalit\'e

(10) $S^-(i)=\epsilon^+(i)^{mult_{\lambda^+}(i)}\epsilon^-(i)^{mult_{\lambda^-}(i)}S^+(i)$,

\noindent o\`u, pour simplifier les notations, on a pos\'e $\epsilon^+(i)=1$ si $i\not\in Jord_{bp}(\lambda^+)$ et $\epsilon^-(i)=1$ si $i\not\in Jord_{bp}(\lambda^-)$.   La v\'erification de cette assertion se fait cas par cas. Traitons seulement le cas o\`u 
$mult_{\lambda^+}(i)mult_{\lambda^-}(i)\not=0$. D'apr\`es la d\'efinition de $c^{\zeta}(i)$ et la remarque 1.13(2), on a la relation $c^+(i)+c^-(i)\equiv mult_{\lambda}(i)\,\,mod\,\, 2{\mathbb Z}$. Supposons d'abord $mult_{\lambda}(i)$ pair. Alors $c^+(i)=c^-(i)$. Si ces deux nombres valent $1$, on a d'apr\`es (9)(c)
$$S^+(i)=2\epsilon^+(i)(-1)^{mult_{\lambda^-}(i)\tau^+(i)},\,\, S^-(i)=2\epsilon^+(i)(-1)^{mult_{\lambda^-}(i)\tau^-(i)}.$$
D'o\`u $S^-(i)=(-1)^{(\tau^+(i)+\tau^-(i))mult_{\lambda^-}(i)}$. 
En vertu de (8), cela \'equivaut  
$$S^-(i)=\epsilon^+(i)^{mult_{\lambda^-}(i)}\epsilon^-(i)^{mult_{\lambda^-}(i)}S^+(i).$$
 Mais, puisque $mult_{\lambda}(i)$ est pair, $mult_{\lambda^+}(i)$ et $mult_{\lambda^-}(i)$ sont de m\^eme parit\'e et l'\'egalit\'e pr\'ec\'edente co\"{\i}ncide avec (10). Si $c^+(i)=c^-(i)=0$, on a d'apr\`es (9)(e)
$$S^+(i)=2(-1)^{mult_{\lambda^-}(i)\tau^+(i)},\,\,S^-(i)=2(-1)^{mult_{\lambda^-}(i)\tau^-(i)},$$
d'o\`u encore $S^-(i)=\epsilon^+(i)^{mult_{\lambda^-}(i)}\epsilon^-(i)^{mult_{\lambda^-}(i)}S^+(i)$ et la m\^eme conclusion. Supposons maintenant $mult_{\lambda}(i)$ impair. Alors $c^+(i)\not=c^-(i)$. Soit $\zeta=\pm$ tel que $c^{\zeta}(i)=1$ et $c^{-\zeta}(i)=0$. D'apr\`es (9)(c) et (9)(e), on a
$$S^{\zeta}(i)=2\epsilon^+(i)(-1)^{mult_{\lambda^-}(i)\tau^{\zeta}(i)},\,\, S^{-\zeta}(i)=2(-1)^{mult_{\lambda^-}(i)\tau^{-\zeta}(i)}.$$
D'o\`u
$$S^{\zeta}(i)=\epsilon^+(i)(-1)^{mult_{\lambda^-}(i)(\tau^+(i)+\tau^-(i))}S^{-\zeta}(i),$$
ou encore, d'apr\`es (8):
$$S^{\zeta}(i)=\epsilon^+(i) (\epsilon^+(i)\epsilon^-(i))^{mult_{\lambda^-}(i)}S^{-\zeta}(i)=\epsilon^+(i)^{1+mult_{\lambda^-}(i)}\epsilon^-(i)^{mult_{\lambda^-}(i)}S^{-\zeta}(i).$$
Mais $mult_{\lambda}(i)$ est impair donc $1+ mult_{\lambda^-}(i)$ est de la m\^eme parit\'e que $mult_{\lambda^+}(i)$. L'\'egalit\'e pr\'ec\'edente co\"{\i}ncide avec (10). Cela d\'emontre (10) dans  le cas o\`u 
$mult_{\lambda^+}(i)mult_{\lambda^-}(i)\not=0$. On laisse les autres cas au lecteur.

En vertu de (7) et (10), on a l'\'egalit\'e
$$(11) \qquad m_{\delta}(\rho_{1},\rho_{2}^+)=cm_{\delta}(\rho_{1},\rho^-_{2}),$$
o\`u
$$c=(\prod_{i\in Jord_{bp}(\lambda^+)}\epsilon^+(i)^{mult_{\lambda^+}(i)})(\prod_{i\in Jord_{bp}(\lambda^-)}\epsilon^-(i)^{mult_{\lambda^-}(i)}).$$
Mais on a vu en \cite{W3} 1.3(1)  que ce produit d\'eterminait l'indice $\sharp$: celui-ci est $iso$ si $c=1$, $an$ si $c=-1$. Autrement dit, $c=sgn_{\sharp}$. L'\'egalit\'e (11) et la non nullit\'e de ses deux membres  entra\^{\i}nent la non-nullit\'e de l'expression (3). Cela ach\`eve de prouver que $m_{\delta,\sharp}(\mu_{1},\eta_{1};\mu_{2},\eta_{2})\not=0$.
 
Il reste une derni\`ere condition \`a prouver, \`a savoir que $n_{2}>0$ si $\sharp=an$. Mais, si $n_{2}=0$, les repr\'esentations $\rho_{2}^+$ et $\rho_{2}^-$ sont les m\^emes: ce sont l'unique repr\'esentation du groupe $W_{0}=\{1\}$. Donc $m_{\delta}(\rho_{1},\rho_{2}^+)=m_{\delta}(\rho_{1},\rho^-_{2})$ et le calcul ci-dessus entra\^{\i}ne que $\sharp=iso$. Cela ach\`eve la v\'erification de la condition 3.7(2) et en m\^eme temps la preuve du th\'eor\`eme 3.3.

{\bf Remarque.}  On peut v\'erifier directement que, si $\sharp=an$, le couple $(n_{1},n_{2})$ fourni par la proposition 1.11, pour notre fonction $\tau$,  v\'erifie $n_{2}\not=0$. En effet, supposons par l'absurde que $n_{2}=0$. A fortiori, l'ensemble d'intervalles de $\lambda_{2}$ est vide donc $j$ et $j+1$ ne sont $2$-li\'es pour aucun $j\geq1$. Pour $j$ impair, la condition 1.11(2) entra\^{\i}ne $\lambda_{j}=\lambda_{j+1}$. Donc les multiplicit\'es $mult_{\lambda}(i)$ sont toutes paires. Puisque $\sharp=an$, on a
$$(\prod_{i\in Jord_{bp}(\lambda^+)}\epsilon^+(i)^{mult_{\lambda^+}(i)})(\prod_{i\in Jord_{bp}(\lambda^-)}\epsilon^-(i)^{mult_{\lambda^-}(i)})=-1.$$
Un $i\in Jord_{bp}(\lambda)$ tel que $mult_{\lambda^+}(i)mult_{\lambda^-}(i)=0$ n'intervient pas dans ce produit. Par exemple, si $mult_{\lambda^-}(i)=0$, il n'intervient \'evidemment pas dans le second produit. Il intervient dans le premier par $\epsilon^+(i)^{mult_{\lambda^+}(i)}$. Mais $mult_{\lambda^+}(i)=mult_{\lambda}(i)$ est pair et cette contribution vaut $1$. En supprimant ces termes et en utilisant que $mult_{\lambda^+}(i)$ et $mult_{\lambda^-}(i)$ sont de m\^eme parit\'e, on obtient
$$\prod_{i\in Jord_{bp}(\lambda^+)\cap Jord_{bp}(\lambda^-)}(\epsilon^+(i)\epsilon^-(i))^{mult_{\lambda^-}(i)}=-1.$$
On peut donc fixer $i\in Jord_{bp}(\lambda^+)\cap Jord_{bp}(\lambda^-)$ tel que $\epsilon^+(i)\epsilon^-(i)=-1$. D'apr\`es  (8), on a $\tau^+(i)\not=\tau^-(i)$. Par construction de ces fonctions, cela entra\^{\i}ne qu'il existe un intervalle $\Delta_{2}\in Int(\lambda_{2})$  tel que $J(\{i\})\subset J(\Delta_{2})$.  A fortiori, $Int(\lambda_{2})$ est non vide, ce qui contredit notre hypoth\`ese $n_{2}=0$. 

\bigskip

{\bf Index des notations}

${\mathbb C}[X]$ 1.1; $cup$ 1.8; $c_{{\cal O}}(\pi)$ 3.2; $d$ 1.2, 1.6, 1.7; $\Delta_{min}$ 1.3, 1.4, 1.5; $\Delta_{max}$ 1.4, 1.5; $\delta^+$, $\delta^-$ 1.13; $\delta(\lambda,s,\epsilon)$ 3.3; $E$ 2.2;  $fam$ 1.3, 1.4, 1.5; $f_{Lie}$ 2.2; $f_{red}$ 2.2; $\phi_{\alpha,\beta',\beta''}$ 2.3; $f_{{\cal O}_{1},{\cal O}_{2}}$ 3.1;  ${\cal H}$ 2.3; $h_{{\cal O}_{1},{\cal O}_{2}}$ 3.1; $Int(\lambda)$ 1.3, 1.4, 1.5;  $ind$ 1.8; $ind(\lambda_{1},\lambda_{2})$ 1.9; $Int_{\lambda_{1},\lambda_{2}}$ 1.10; $Jord(\lambda)$ 1.1; $Jord_{bp}(\lambda)$ 1.3, 1.4, 1.5; $Jord_{bp}^k(\lambda)$ 1.3, 1.4, 1.5; $J(\Delta)$ 1.6, 1.7; $j_{min}(\Delta)$ 1.6, 1.7; $j_{max}(\Delta)$ 1.6, 1.7; $k_{\lambda,\epsilon}$ 1.3, 1.4, 1.5; $k(r',r'';w)$ 2.3; $k(w)$ 2.3; $\kappa_{\pi,0}$ 2.3;  $l(\lambda)$ 1.1; $mult_{\lambda}(i)$ 1.1; $mult_{\lambda}(\geq i)$ 1.1; $m(\rho_{1},\rho_{2}^{\zeta}; \rho_{\lambda',\epsilon'},\rho_{\lambda'',\epsilon''})$ 1.12; $mult!_{{\bf m}}$ 2.1; $\mu({\cal O})$ 3.1; $\mu(\pi)$ 3.2; $Nil_{\sharp}$ 3.1; $\boldsymbol{Nil}_{\sharp}$ 3.1; $\mathfrak{n}(\lambda^+,\lambda^-)$ 3.8; $\mathfrak{N}(\lambda^+,\lambda^-)$ 3.8; ${\cal O}_{{\cal O}_{1},{\cal O}_{2}}$ 3.1; ${\cal P}(N)$ 1.1; ${\cal P}_{k}(N)$ 1.1; ${\cal P}^{symp}(2n)$ 1.3; $\boldsymbol{{\cal P}^{symp}}(2n)$ 1.3; ${\cal P}^{symp,sp}(2n)$ 1.3; 
${\cal P}^{orth}(2n+1)$ 1.4; $\boldsymbol{{\cal P}^{orth}}(2n+1)$ 1.4; ${\cal P}^{orth,sp}(2n+1)$ 1.4; ${\cal P}^{orth}(2n)$ 1.5; $\boldsymbol{{\cal P}^{orth}}(2n)$ 1.5; ${\cal P}^{orth,sp}(2n)$ 1.5; $\underline{{\cal P}}^{orth}(2n)$ 1.5; $\boldsymbol{\underline{{\cal P}}^{orth}}(2n)$ 1.5; ${\cal P}^{orth}({\bf n})$ 1.8; $\rho(\alpha,\beta)$ 1.1; $\rho^D(\alpha,\beta)$ 1.1; $\rho_{\lambda,\epsilon}$ 1.3, 1.4, 1.5; 
$\rho_{2}^+$, $\rho_{2}^-$ 1.12; $S(\lambda)$ 1.1; $S_{k}(\lambda)$ 1.1; $\mathfrak{S}_{N}$ 1.1; $sgn$ 1.1; $sgn_{CD}$ 1.1; ${\cal S}_{N,D}$ 1.2; $symb$ 1.2; $sp(\lambda)$ 1.3, 1.4, 1.5; $sp(\lambda,\epsilon)$ 1.3, 1.4, 1.5;  $\tau_{\lambda_{1},\lambda_{2}}$ 1.10;  $\tau^+$, $\tau^-$ 1.13; $\Theta_{\pi}$ 2.1; $\Theta_{\pi,cusp}$ 2.1; $\Theta_{\pi,{\bf m},cusp}$ 2.1; $\Theta_{\pi,cusp}^M$ 2.1;  $V(\pi)$ 3.2; $W_{N}$ 1.1; $W_{N}^D$ 1.1; $W({\bf m},N',N'')$ 2.3; $W({\bf m},N',N'')_{ell}$ 2.3; $W_{n_{2},iso}$, $W_{n_{2},an}$ 3.6; $\xi$ 1.9; $\zeta(\lambda)$ 1.6, 1.7.

\bigskip

 {\bf Index des notations de \cite{W3}}
 
${\mathbb C}[X]$ 1.4; $C'_{n'}$ 1.5; $C^{''\pm}_{n'',\sharp}$ 1.5; $C''_{n''}$ 1.5; $C^{GL(m)}$ 1.5; ${\mathbb C}[\hat{W}_{N}]_{cusp}$ 1.8; $D(n)$ 1.2; $D_{iso}(n)$ 1.2;  $D_{an}(n)$ 1.2; $D$ 1.7; $D^{par}$ 1.7; $\eta(Q)$ 1.1; $\eta^+(Q)$ 1.1; $\eta^-(Q)$ 1.1; $Ell_{unip}$ 1.4; $\mathfrak{Ell}_{unip}$ 1.4; $\mathfrak{Endo}_{tunip}$ 2.1; $\mathfrak{Endo}_{unip-quad}$ 2.2; $\mathfrak{Endo}_{unip-quad}^{red}$ 2.2; $\mathfrak{Endo}_{unip,disc}$ 2.4;  ${\cal F}^L$ 1.9; ${\cal F}^{par}$ 1.9; ${\cal F}$ 2.3: $\mathfrak{F}^{par}$ 2.3; $G_{iso}$ 1.1; $G_{an}$ 1.1; $\Gamma$ 1.8; $\boldsymbol{\Gamma}$ 1.8; $\tilde{GL}(2n)$ 2.1; $Irr_{tunip}$ 1.3; $\mathfrak{Irr}_{tunip}$ 1.3; $Irr_{unip-quad}$1.3; $\mathfrak{Irr}_{unip-quad}$ 1.3; $Jord(\lambda)$ 1.3; $Jord_{bp}(\lambda)$ 1.3; $Jord_{bp}^{k}(\lambda)$ 1.4; $K_{n',n''}^{\pm}$ 1.2; $k$ 1.9; $L^*$ 1.1; $L_{n',n''}$ 1.2; $l(\lambda)$ 1.3; $mult_{\lambda}$ 1.3; $\mathfrak{o}$ 1.1; $O^+(Q)$ 1.1; $O^-(Q)$ 1.1; $\varpi$ 1.1; $\pi_{n',n''}$ 1.3; ${\cal P}(N)$ 1.3; ${\cal P}^{symp}(2N)$ 1.3; $\boldsymbol{{\cal P}^{symp}}(2N)$ 1.3; $\pi(\lambda,s,\epsilon)$ 1.3; $\pi(\lambda^+,\epsilon^+,\lambda^-,\epsilon^-)$ 1.3; $\pi_{ell}(\lambda^+,\epsilon^+,\lambda^-,\epsilon^-)$ 1.4; $proj_{cusp}$ 1.5; ${\cal P}(\leq n)$ 1.5; ${\cal P}_{k}(N)$ 1.8; $\Pi(\lambda,s,h)$ 2.1; $\Pi^{st}(\lambda^+,\lambda^-)$ 2.4; ${\cal P}^{symp,disc}(2n)$ 2.4; $Q_{iso}$ 1.1; $Q_{an}$ 1.1; $\rho_{\lambda}$ 1.3;  ${\cal R}^{par}$ 1.5; ${\cal R}^{par,glob}$ 1.5; ${\cal R}^{par}_{cusp}$ 1.5; ${\cal R}^{par,glob}_{{\bf m}}$ 1.5; ${\cal R}^{par}_{{\bf m},cusp}$ 1.5; $res'_{m}$ 1.5; $res''_{m}$ 1.5; $res_{m}$ 1.5 et 1.8; $res_{{\bf m}}$ 1.5; ${\cal R}$ 1.8; ${\cal R}(\gamma)$ 1.8; ${\cal R}(\boldsymbol{\gamma})$ 1.8; ${\cal R}^{glob}$ 1.8; ${\cal R}_{cusp}$ 1.8; $Rep$ 1.9; $\rho\iota$ 1.10; $S(\lambda)$ 1.3; $\mathfrak{S}_{N}$ 1.8; $\hat{\mathfrak{S}}_{N}$ 1.8; $sgn$ 1.8; $sgn_{CD}$ 1.8; ${\cal S}_{n}$ 1.11; $\mathfrak{St}_{tunip}$ 2.1; $\mathfrak{St}_{unip-quad}$ 2.4; $\mathfrak{St}_{unip,disc}$ 2.4;  $sgn_{iso}$ 2.6; $sgn_{an}$ 2.6; $val_{F}$ 1.1; $V_{iso}$ 1.1; $V_{an}$ 1.1; $W_{N}$ 1.8; $\hat{W}_{N}$ 1.8; $w_{\alpha}$ 1.8; $w_{\alpha,\beta}$ 1.8; $w_{\alpha,\beta',\beta''}$ 1.8; $Z(\lambda)$ 1.3; $Z(\lambda,s)$ 1.3; ${\bf Z}(\lambda,s)$ 1.3; ${\bf Z}(\lambda,s)^{\vee}$ 1.3; $\vert .\vert _{F}$ 1.1.

CNRS IMJ-PRG

4 place Jussieu

75005 Paris

jean-loup.waldspurger@imj-prg.fr
 
\end{document}